\documentclass[12pt,leqno]{article}

\usepackage{amsthm}
\usepackage[usenames]{color}

\usepackage{amsmath}
\usepackage{amscd}
\usepackage{amsmath}
\usepackage{amssymb}
\usepackage{latexsym}
\makeatletter

\@addtoreset{equation}{section}
\makeatother

\newtheorem{thm}{Theorem}[subsection]
\newtheorem{pr}[thm]{Proposition}
\newtheorem{df}[thm]{Definition}
\newtheorem{lm}[thm]{Lemma}
\newtheorem{cor}[thm]{Corollary}

\newcommand{\sm}{\raisebox{2.33pt}{~\rule{6.4pt}{1.3pt}~}}

\input xy \xyoption{all} \CompileMatrices

\begin{document}

\title{On singular supports in mixed characteristic}
\author{Takeshi SAITO}

\maketitle

\begin{abstract}
We fix an excellent regular noetherian scheme
$S$ over ${\mathbf Z}_{(p)}$ 
satisfying a certain finiteness condition.
For a constructible \'etale sheaf
${\cal F}$ on a regular scheme $X$ of
finite type over $S$,
we introduce a variant of the singular support
relatively to $S$
and prove the existence of
a saturated relative variant of the singular support
by adopting the method of Beilinson in \cite{SS}
using the Radon transform.
We may deduce the
existence of the singular support itself,
if we admit an expected property on
the micro support of tensor product
and if the scheme $X$ is sufficiently ramified
over the base $S$.
\end{abstract}

For a constructible \'etale sheaf on 
a smooth scheme over a field,
Beilinson proved in \cite{SS} the existence of
the singular support and its fundamental properties.
A key tool used in the proof of the existence is the 
Radon transform studied in \cite{Radon}.
In mixed characteristic, 
the existence of singular support 
remains an open question
although
a framework to study the singular support 
was prepared in \cite{ANT} using
the Frobenius--Witt cotangent bundle $FT^*X$
supported on the characteristic $p$ fiber
$X_{{\mathbf F}_p}$.
We prove in Theorem \ref{thmSS} the existence of
a saturated relative variant of the singular support
in mixed characteristic by adopting the method of 
Beilinson.

Using Radon transform requires
the introduction of two ingredients.
The first is to work in a relative setting
with a fixed base scheme $S$.
The second is to work with 
pairs $(h,f)$ of morphisms
$h\colon W\to X$ and $f\colon W\to Y$
of schemes over $S$.
The Radon transform
is defined for a projective space
${\mathbf P}={\mathbf P}_S$
over a fixed scheme $S$
using the pair $(q,q^\vee)$ 
of the universal family 
$q^\vee \colon Q\to {\mathbf P}^\vee$
of hyperplanes in ${\mathbf P}$ 
and the canonical morphism
$q\colon Q\to {\mathbf P}$.
Thus, we need to work in a relative setting over $S$
and also with pairs $(h,f)$ of morphisms
while the definition of singular support
in \cite{ANT} is formulated only with morphisms $h$
without the second morphisms $f$.

The format of the definitions of various variants of
singular supports is summarized as follows.
First, we introduce a certain condition
on morphisms of schemes
with respect to a given closed conical
subset $C$ in the cotangent bundle and a corresponding
condition with respect to an \'etale sheaf ${\cal F}$.
Next, we define to say that
${\cal F}$ is micro supported on $C$
if every morphism satisfying the condition for $C$
satisfies the corresponding condition for ${\cal F}$.
The singular support $SS{\cal F}$ is then
defined to be the smallest $C$ on which
${\cal F}$ is micro supported.
The existence of singular support is non-trivial 
because
the condition that 
${\cal F}$ is micro supported on $C$ and on $C'$
does not a priori imply
that 
${\cal F}$ is micro supported on
the intersection $C\cap C'$.

In the original definition by Beilinson in \cite{SS},
he considers pairs $(h,f)$ of morphisms
and the condition on ${\cal F}$ is
given by the local acyclicity of $f$ relatively to $h^*{\cal F}$.
Since a naive imitation does not give a
good definition of singular support
in mixed characteristic,
we consider in \cite{ANT} morphisms $h$ without $f$
and the ${\cal F}$-transversality
as a condition on $h$ relatively to ${\cal F}$:
A separated morphism $h$ of finite type
is ${\cal F}$-transversal if the canonical morphism
$$c_{{\cal F},h}\colon h^*{\cal F}\otimes h^!\Lambda
\to h^!{\cal F}$$
is an isomorphism. Here and in the following
the letter $R$ indicating the derived functor is omitted.

We fix a base scheme $S$.
To define in Definition \ref{dfSss}
the $S$-singular support
of a sheaf ${\cal F}$ on a regular scheme
$X$ of finite type over $S$ along the format above,
we introduce in Definition \ref{dfFacyc}
a condition
for a pair of morphisms $h\colon W\to X$
and $f\colon W\to Y$ over $S$
to be ${\cal F}$-acyclic over $S$,
by modifying the local acyclicity of
$f$ relatively to $h^*{\cal F}$ used by Beilinson
in his definition of the singular support.
For the modification,
we use a criterion of the local acyclicity
due to Braverman--Gaitsgory \cite[Appendix Theorem B.2]{Eis}
recalled in Proposition \ref{prBG}:
A morphism $f\colon X\to Y$
to a smooth scheme $Y$
over a field $k$
is locally acyclic relatively to ${\cal F}$ if
and only if the morphism
$(1,f)\colon X\to X\times_kY$ is 
${\cal F}\boxtimes {\cal G}$-transversal
for every sheaf ${\cal G}$ on $Y$.

If $S$ is a scheme over a field $k$
and if $Y=Y_0\times_kS$ is the base change of 
a smooth scheme $Y_0$ over $k$,
then we may consider the condition that
$(h,f)\colon W\to X\times_SY=
X\times_kY_0$ is 
${\cal F}\boxtimes {\cal G}$-transversal
for the pull-back ${\cal G}$ to $Y$
of every sheaf ${\cal G}_0$ on $Y_0$.
In mixed characteristic case
where a base field cannot exist,
we need to find a replacement of
the condition for a sheaf ${\cal G}$
on a smooth scheme $Y$ over $S$
to `come from $Y_0$ as pull-back'.
We will formulate such a condition
called the $S$-acyclicity
in Definition \ref{dfFf0}
by using micro support defined
in an absolute setting without 
a base scheme and use it
to define the ${\cal F}$-acyclicity over $S$
in Definition \ref{dfFacyc}.

The ${\cal F}$-acyclicity over $S$ enables us
to define the $S$-singular support
$SS_S{\cal F}$ using the format above
and to ask its existence
as a closed conical subset
of the FW cotangent bundle $FT^*X$.
A key reason why
the use of Radon transform 
is effective in the proof of the
existence of singular support is
that the universal family $Q$
of hyperplanes in ${\mathbf P}$ over $S$
is
canonically identified
with the projective space bundle
${\mathbf P}(T^*{\mathbf P}/S)$
associated to the relative cotangent bundle
$T^*{\mathbf P}/S$.
To fill the gap between 
the FW cotangent bundle $FT^*X$
and the relative cotangent bundle
$T^*{\mathbf P}/S$,
we need to introduce
the $S$-saturation
$SS_S^{\rm sat}{\cal F}$
defined as the smallest closed
conical subset stable
under the action of
$FT^*S\times_SX$ and
containing
$SS_S{\cal F}$ as a subset.

With this variant, the method of Beilinson
allows us to prove the existence
of $SS_S^{\rm sat}{\cal F}$ in
Theorem \ref{thmSS}.
If $X$ is smooth over $S$,
the $S$-saturated singular support
$SS_S^{\rm sat}{\cal F}$
defines a closed conical subset
of the relative cotangent bundle
$T^*X/S|_{X_{{\mathbf F}_p}}$
restricted to the characteristic $p$ fiber.
In the case where $S={\rm Spec}\, k$
for a perfect field $k$ of characteristic $p>0$,
we have $T^*X/S|_{X_{{\mathbf F}_p}}=T^*X$ and
recover the singular support defined by
Beilinson in \cite{SS}. If $k$ is not perfect,
ours is different from that in \cite{SS}
as we see in Example after Definition \ref{dfSss}.
If we admit a certain property on 
micro support of tensor product,
we would have
$SS_S{\cal F}
\subset
SS{\cal F}
\subset
SS_S^{\rm sat}{\cal F}$
by Lemma \ref{lmS}.
The inclusions may be strict in general.
If we further assume
that the morphism
$FT^*S\times_SX\to FT^*X$
is 0,
then we have
$SS_S{\cal F}
=SS_S^{\rm sat}{\cal F}$.

In Section \ref{sms},
we recall the definition of singular support
from \cite{ANT} and give various complements.
We recall the definition of ${\cal F}$-transversality
and of the $C$-transversality
and prepare their basic properties in Sections \ref{ssFtr}
and \ref{ssCtr} respectively.
We study the relation between the 
${\cal F}$-transversality and the local acyclicity
in Section \ref{ssFacyc}
and close the section with a proof of
the criterion of the local acyclicity
due to Braverman--Gaitsgory in Proposition \ref{prBG}.
After recalling the definition of micro support
in Section \ref{ssms},
we recall the definition of the singular support in
Section \ref{ssSS} following the general format.
We show in Corollary \ref{corSS}
that the existence of singular support
is reduced to the case of projective spaces.
Finally, we introduce the $S$-acyclicity as
a preparation of the definition of
the ${\cal F}$-acyclicity over $S$.

In Section \ref{sSms},
we introduce the $S$-singular support
and prove the existence of
its $S$-saturation.
First, in Sections \ref{ssCacycS} and \ref{ssFacycS},
we define the $C$-acyclicity over $S$
and the ${\cal F}$-acyclicity over $S$
and study their basic properties
respectively.
Following the general format,
we define the notion of 
$S$-micro support in Section \ref{ssSms}
and the $S$-singular support in Section \ref{ssSSS}
respectively.
After recalling basic properties of
the Legendre transform and the Radon transform,
we prove the existence of
$S$-saturated singular support in
Theorem \ref{thmSS} adopting the method of
Beilinson.

The research was partially supported by 
Kakenhi Grant-in-Aid (C) 24K06683.
The author thanks greatly Daichi Takeuchi for
proposing an essential improvement on the main result
and correcting some errors in earlier versions.
He also thanks an anonymous referee
for extremely careful reading and
for pointing out several inaccuracies and
numerous typographical errors.

\tableofcontents

\bigskip

Let $\Lambda$ be a finite field
of characteristic $\ell$.
Schemes are over ${\mathbf Z}[1/\ell]$.
An object of
the derived category
$D^+(X,\Lambda)$
of bounded below
complexes on the \'etale site of 
a scheme $X$
will be called a {\em sheaf} on $X$
by abuse of terminology.
If a sheaf belongs to the
full sub category
$D^b_c(X,\Lambda)$ consisting of
bounded constructible complexes,
we say it is {\em constructible}.
A constructible
sheaf ${\cal F}$ is said to be
locally constant if
every cohomology sheaf
${\cal H}^q{\cal F}$ is locally constant.
We drop $R$ in the notation for
derived functors
and $\Lambda$
in $\otimes_\Lambda$.

\section{Micro support}\label{sms}

\subsection{${\cal F}$-transversality}\label{ssFtr}

Let $h\colon W\to X$ be a separated morphism
of finite type of noetherian schemes
and let ${\cal F}$ and ${\cal G}$
be sheaves on $X$.
We define a canonical morphism
$$c_{h,{\cal F},{\cal G}}\colon
h^*{\cal F}\otimes h^!{\cal G}
\to
h^!({\cal F}\otimes{\cal G})$$
to be the adjoint of
the composition
$$\begin{CD}
h_!(h^*{\cal F}\otimes h^!{\cal G})
@>>>
{\cal F}\otimes h_!h^!{\cal G}
@>{1_{\cal F}\otimes{\rm adj}}>>
{\cal F}\otimes{\cal G}
\end{CD}$$
where the first arrow is
the isomorphism of projection formula
and the second is induced
by the adjunction.
There is another way to view
the morphism $c_{h,{\cal F},{\cal G}}$ 
under the assumption that ${\cal F}$ is bounded.
By taking the composition 
of the canonical morphism
$h^!{\cal G}\to h^!{\cal H}om({\cal F},
{\cal F}\otimes {\cal G})$
with the inverse of
the isomorphism
${\cal H}om(h^*{\cal F},
h^!({\cal F}\otimes {\cal G}))
\to h^!{\cal H}om({\cal F},
{\cal F}\otimes {\cal G})$
of
\cite[Corollaire 3.1.12.2]{PD},
\cite[Theorem 8.4.7]{Fu},
we obtain a morphism
$h^!{\cal G}\to 
{\cal H}om(h^*{\cal F},
h^!({\cal F}\otimes {\cal G}))$.
Taking the adjoint,
we obtain
$c_{h,{\cal F},{\cal G}}\colon
h^*{\cal F}\otimes h^!{\cal G}
\to h^!({\cal F}\otimes {\cal G})$.

For ${\cal G}=\Lambda$,
$$c_{h,{\cal F}}\colon
h^*{\cal F}\otimes h^!\Lambda
\to
h^!{\cal F}$$
is defined as
$c_{h,{\cal F},\Lambda}$.

\begin{lm}\label{lmcGL}
The following diagram is commutative:
$$\begin{CD}
h^*{\cal F}\otimes h^!{\cal G}
@>{c_{h,{\cal F},{\cal G}}}>>
h^!({\cal F}\otimes{\cal G})
\\
@A{1_{h^*{\cal F}}\otimes
c_{h,{\cal G}}}AA
@AA{c_{h,{\cal F}\otimes{\cal G}}}A
\\
h^*{\cal F}\otimes 
h^*{\cal G}\otimes h^!\Lambda
@>{{\rm can}\otimes 1}>>
h^*({\cal F}\otimes 
{\cal G})\otimes h^!\Lambda.
\end{CD}
$$
\end{lm}

\proof{
By taking the adjoint and applying the projection formula,
it suffices to show that
the diagram
$$\begin{CD}
{\cal F}\otimes h_!h^!{\cal G}
@>{1_{\cal F}\otimes{\rm adj}}>>
{\cal F}\otimes{\cal G}
\\
@A{1_{\cal F}\otimes c_{h,{\cal G}}}AA
@AA{1_{\cal F}\otimes1_{\cal G}
\otimes{\rm adj}
}A
\\
{\cal F}\otimes 
h_!(h^*{\cal G}\otimes h^!\Lambda)
@>{1_{\cal F}\otimes{\rm proj}}>>
{\cal F}\otimes 
{\cal G}\otimes h_!h^!\Lambda
\end{CD}
$$
is commutative.
It suffices to show the commutativity 
before tensoring ${\cal F}$.
Then, it follows from the definition
of $c_{h,{\cal G}}$.
\qed

}
\medskip

We prove the transitivity of
$c_{h,{\cal F}}$.

\begin{lm}\label{lmchtr}
Let $h\colon W\to X'$ and $g\colon X'\to X$
be separated morphisms
of finite type of noetherian schemes.
Then, the following diagram is commutative:
\begin{equation}
\xymatrix{
(gh)^*{\cal F}\otimes (gh)^!\Lambda
\ar[rr]^-{c_{gh,{\cal F}}}
\ar[rrd]^{c_{h,g^*{\cal F},g^!\Lambda}}
&&
(gh)^!{\cal F}
\\
h^*(g^*{\cal F}\otimes g^!\Lambda)
\otimes h^!\Lambda
\ar[u]^{1\otimes c_{h,g^!\Lambda}}
\ar[rr]^-{c_{h,g^*{\cal F}\otimes g^!\Lambda}}
&&
h^!(g^*{\cal F}\otimes g^!\Lambda)
\ar[u]_{h^!(c_{g,{\cal F}})}
}
\label{eqchtr}
\end{equation}
\end{lm}

\proof{
The lower left triangle is commutative
by Lemma \ref{lmcGL}.
To prove that
the upper right triangle is commutative,
we consider the diagram
$$\xymatrix{
(gh)_!((gh)^*{\cal F}\otimes (gh)^!\Lambda)
\ar[rr]^{\rm proj}\ar[d]&&
{\cal F}\otimes (gh)_!(gh)^!\Lambda
\ar[d]
\ar[ddr]^{\rm adj}&
\\
g_!h_!(h^*g^*{\cal F}\otimes h^!g^!\Lambda)
\ar[r]^{\rm proj}\ar[dr]&
g_!(g^*{\cal F}\otimes h_!h^!g^!\Lambda)
\ar[r]^{\rm proj}\ar[d]^{\rm adj}&
{\cal F}\otimes g_!h_!h^!g^!\Lambda
\ar[d]^{\rm adj}
&
\\
&
g_!(g^*{\cal F}\otimes g^!\Lambda)
\ar[r]^{\rm proj}&
{\cal F}\otimes g_!g^!\Lambda
\ar[r]^{\rm adj}
&{\cal F}.
}
$$
The morphism $c_{gh,{\cal F}}$ is 
induced by the composition 
via the upper right,
$c_{g,{\cal F}}$ is 
induced by the composition of the lower line
and
$c_{h,g^*{\cal F},g^!\Lambda}$
is induced by the
lower left slant arrow
without $g_!$.
Hence, it is reduced to showing
the commutativity of the diagram.

By the transitivity of the isomorphisms
of projection formula,
the upper rectangle of the diagram
is commutative.
The lower square is commutative
by the functoriality 
of the isomorphism of projection formula
and the right triangle is commutative
by the transitivity of the adjunction.
Hence the diagram is commutative
as required.
\qed

}
\medskip

We prove the compatibility with
base change.
We consider a cartesian diagram 
\begin{equation}
\begin{CD}
W'@>{h'}>>X'\\
@VgVV@VVfV\\
W@>h>>X,
\end{CD}
\label{eqbc}
\end{equation}
of noetherian schemes
where the horizontal arrows are
separated morphisms
of finite type.
The base change morphism
\begin{equation}
g^*h^!\to h'^!f^*
\label{eqbcdef}
\end{equation}
\cite[(3.1.14.2)]{PD}
is defined as the adjoint
of the composition
$h'_!g^*h^!\to 
f^*h_!h^!\to f^*$
of the inverse of the isomorphism
of proper base change theorem
and the adjunction.
As the adjoint of the isomorphism
$f^*h_!\to h'_!g^*$ of
proper base change theorem,
we obtain an isomorphism
\begin{equation}
h^!f_*\to g_*h'^!
\label{eqbcdef2}
\end{equation}
\cite[Corollaire 3.1.12.3]{PD},
\cite[Proposition 8.4.9]{Fu}.
The morphism (\ref{eqbcdef2})
is the adjoint of the composition 
of $g^*h^!f_*\to h'^!f^*f_*\to h'^!$
where the first arrow is induced by
(\ref{eqbcdef}) and the second is
induced by the adjunction.

\begin{lm}\label{lmbc}
Let {\rm (\ref{eqbc})}
be a cartesian diagram of
noetherian schemes
and assume that
the horizontal arrows are
separated morphisms
of finite type.

{\rm 1.}
For a sheaf ${\cal F}$ on $X$,
the diagram
\begin{equation}
\begin{CD}
g^*(h^*{\cal F}\otimes h^!\Lambda)
@>{g^*c_{h,{\cal F}}}>>
g^*h^!{\cal F}\\
@V{\rm can\otimes (\ref{eqbcdef})}VV@VV{\rm (\ref{eqbcdef})}V\\
h'^*f^*{\cal F}\otimes h'^!\Lambda
@>{c_{h',f^*{\cal F}}}>>
h'^!f^*{\cal F}
\end{CD}
\label{eqcbc}
\end{equation}
is commutative.

{\rm 2.}
Let ${\cal F}'$ be a sheaf  on
$X'$.
Then, 
the following diagram is commutative:
\begin{equation}
\begin{CD}
h^*f_*{\cal F}'\otimes h^!\Lambda
@>{c_{h,f_*{\cal F}'}}>>
h^!f_*{\cal F}'\\
@VVV@VV{\rm (\ref{eqbcdef2})}V\\
g_*(h'^*{\cal F}'\otimes h'^!\Lambda)
@>{g_*(c_{h',{\cal F}'})}>>
g_*h'^!{\cal F}'
\end{CD}
\label{eqprbc}
\end{equation}
where the left vertical arrow is
the adjoint of
$$\begin{CD}
g^*(h^*f_*{\cal F}'\otimes h^!\Lambda)
=
g^*h^*f_*{\cal F}'\otimes g^*h^!\Lambda
@>{{\rm can}\otimes1}>>
h'^*f^*f_*{\cal F}'\otimes g^*h^!\Lambda
@>{h'^*{\rm adj}\otimes {\rm 
(\ref{eqbcdef})}}>> 
h'^*{\cal F}'\otimes h'^!\Lambda.
\end{CD}$$
\end{lm}

\proof{
1.
By the definition of (\ref{eqbcdef}),
it suffices to show that
the diagram
$$\begin{CD}
f^*h_!(h^*{\cal F}\otimes h^!\Lambda)
@>>>
f^*({\cal F}\otimes h_!h^!\Lambda)
@>{f^*(1\otimes{\rm adj})}>>
f^*{\cal F}\\
@VVV@VV{{\rm can
\otimes (\ref{eqbcdef})}}V@|\\
h'_!(h'^*f^*{\cal F}\otimes h'^!\Lambda)
@>>>
f^*{\cal F}\otimes h'_!h'^!\Lambda
@>{1\otimes{\rm adj}}>>
f^*{\cal F}
\end{CD}$$
is commutative.
The vertical arrows are induced by the base
change morphisms.
The left square is the compatibility
of the projection formula with base change.
By the definition of (\ref{eqbcdef}),
the diagram
$$\begin{CD}
h'_!g^*h^!@>>> h'_!h'^!f^*\\
@VVV@VV{\rm adj}V\\
f^*h_!h^!@>{\rm adj}>>f^*
\end{CD}$$
is commutative.
This implies the commutativity of 
the right square.

2. 
By the commutative diagram (\ref{eqcbc})
in 1 for ${\cal F}=f_*{\cal F}'$
and the functoriality of
$c_{h',{\cal F}}$ for
$f^*f_*{\cal F}\to {\cal F}$,
we obtain a commutative diagram
$$\begin{CD}
g^*(h^*f_*{\cal F}'\otimes h^!\Lambda)
@>{g^*c_{h,f_*{\cal F}'}}>>
g^*h^!f_*{\cal F}'\\
@VVV@VVV\\
h'^*{\cal F}'\otimes h'^!\Lambda
@>{c_{h',{\cal F}'}}>>
h'^!{\cal F}'.
\end{CD}$$
Taking the adjoint,
we obtain the commutative diagram
since (\ref{eqbcdef2})
is the adjoint of the composition 
of the left vertical arrow.
\qed

}

\begin{df}\label{dfXtr}
Let $h\colon W\to X$ be a morphism
of finite type of regular noetherian schemes.

{\rm 1.}
If, locally on $W$,
$h$ is the composition of
a smooth morphism $P\to X$ of
relative dimension $n$
and a regular immersion
$W\to P$ of codimension $c$,
we say that $W$ is of relative dimension $n-c$
over $X$.

{\rm 2.}
Let $g\colon X'\to X$ be a morphism
of regular noetherian schemes
and
let $U\subset W'=W\times_XX'$ be
an open subscheme.
We say that $h$ is transversal to $g$
on $U$ if $U$ is regular
and if the relative dimension of $U$ over $X'$
equals the relative dimension of $W$ over $X$.
\end{df}

\begin{lm}\label{lmpurity}
Let $h\colon W\to X$
be a separated morphism of finite type
of regular noetherian schemes.

{\rm 1.}
Let $d$ be the relative dimension of $h$.
Then, there is a canonical isomorphism
$h^!\Lambda\to \Lambda(d)[2d]$.

{\rm 2.}
Let $g\colon X'\to X$ be a morphism
of regular noetherian schemes and let
$$
\begin{CD}
W'@>{h'}>>X'\\
@V{g'}VV@VV{g}V\\
W@>h>> X
\end{CD}$$
be a cartesian diagram of
noetherian schemes.
Assume that $h$ is transversal to $g$
on an open subscheme $U'\subset W'$.
Then 
the base change morphism
$g'^*h^!\Lambda\to h'^!\Lambda$ 
{\rm (\ref{eqbcdef})} is
an isomorphism on $U'$.
\end{lm}

\proof{
Since the assertion is local on $W$,
we may factor $h$ into the composition of
a smooth morphism and an immersion.
The smooth case is the Poincar\'e duality 
\cite[Th\'eor\`eme 3.2.5]{PD} and
the immersion case 
follows from the purity
\cite[Expos\'e XVI, Th\'eor\`eme 3.1.1]{Gabber}.
\qed

}

\begin{df}\label{dfFtrans}
Let $h\colon W\to X$ be a separated morphism
of finite type of noetherian schemes
and let ${\cal F}$ be a sheaf on $X$.

We say that $h$ is ${\cal F}$-transversal
if $c_{h,{\cal F}}\colon
h^*{\cal F}\otimes h^!\Lambda
\to
h^!{\cal F}$ is an isomorphism.
\end{df}

\begin{lm}\label{lmFtrans}
Let $h\colon W\to X$ be a separated morphism
of finite type of noetherian schemes
and let ${\cal F}$ be a sheaf on $X$.

{\rm 1.}
If $h$ is smooth,
then $h$ is ${\cal F}$-transversal.

{\rm 2.}
If ${\cal F}'\to {\cal F}\to {\cal F}''\to$
is a distinguished triangle
of sheaves on $X$,
and if $h$ is transversal
for two of
${\cal F}',{\cal F}$ and ${\cal F}''$.
Then,
$h$ is transversal
for the rest.

{\rm 3.}
If ${\cal F}$ is locally constant,
then every morphism $h$ is ${\cal F}$-transversal.

{\rm 4.}
Assume that
$h$ is a closed immersion and that
${\cal F}$ is a constructible sheaf
supported on the image $h(W)\subset X$.
If $h$ is ${\cal F}$-transversal
and if
$h^!\Lambda$ is isomorphic
to $\Lambda(-c)[-2c]$
for an integer $c\neq 0$,
then ${\cal F}=0$.

{\rm 5.}
Assume that $h$ is a closed immersion
and let $j\colon U=X\sm W\to X$
be the open immersion of the complement.
Then, the following conditions
are equivalent:

{\rm (1)}
$h$ is ${\cal F}$-transversal.

{\rm (2)}
The canonical morphism
${\cal F}\otimes j_*\Lambda
\to j_*j^*{\cal F}$ is an isomorphism.
\end{lm}

We show a partial converse of 3
for constructible sheaves
on regular schemes
in Proposition \ref{prlcc} below.

\proof{
1.
Poincar\'e duality \cite[Th\'eor\`eme 3.2.5]{PD}.

2.
It suffices to consider the morphism
$$\begin{CD}
h^*{\cal F}'\otimes h^!\Lambda
@>>>
h^*{\cal F}\otimes h^!\Lambda
@>>>
h^*{\cal F}''\otimes h^!\Lambda
@>>>\\
@V{c_{h,{\cal F}'}}VV@V{c_{h,{\cal F}}}VV@V{c_{h,{\cal F}''}}VV@.\\
h^!{\cal F}'
@>>>
h^!{\cal F}
@>>>
h^!{\cal F}''
@>>>
\end{CD}$$
of distinguished triangles.

3.
Since the question is \'etale local,
we may assume that ${\cal F}$
is constant and the assertion follows.

4. Since
${\cal F}$ is supported on $h(W)$,
the morphism
$h^!{\cal F}\to h^*{\cal F}$
is an isomorphism.
Further if $
h^*{\cal F}(-c)[-2c]
\to h^!{\cal F}$
is an isomorphism for
$c\neq 0$,
we have $h^*{\cal F}=0$
and ${\cal F}=0$.

5.
In the commutative diagram
$$\begin{CD}
{\cal F}\otimes
h_*h^!\Lambda
@>>>
{\cal F}
@>>>
{\cal F}\otimes
j_*\Lambda
@>>>\\
@VVV
@|
@VVV@.
\\
h_*h^!{\cal F}
@>>>
{\cal F}
@>>>
j_*j^*{\cal F}@>>>
\end{CD}$$
of distinguished triangles,
the left vertical arrow 
is an isomorphism if and only if
so is the right vertical arrow.
By the projection formula,
the left vertical arrow 
is identified with
$h_*c_{h,{\cal F}}$
and the assertion follows.
\qed

}

\medskip

We show the transitivity of the transversality.

\begin{lm}\label{lmFtrtr}
Let $h\colon W\to X'$
and $g\colon X'\to X$ be separated morphisms
of finite type of noetherian schemes
and let ${\cal F}$ be a sheaf on $X$.
Assume that $g^!\Lambda$ is 
isomorphic to $\Lambda(d)[2d]$
for some integer $d$
and that $g$ is ${\cal F}$-transversal.
Then, the following conditions are equivalent:

{\rm (1)} $gh$ is ${\cal F}$-transversal.

{\rm (2)} $h$ is $g^*{\cal F}$-transversal.
\end{lm}

The assumption
that $g^!\Lambda$ is 
isomorphic to $\Lambda(d)[2d]$
is satisfied if $g$ is smooth of relative dimension
$d$.

\proof{
We consider the commutative diagram
(\ref{eqchtr}).
Since $g^!\Lambda$ is locally constant,
the left vertical arrow is an isomorphism
by Lemma \ref{lmFtrans}.3.
The assumption that
$g$ is ${\cal F}$-transversal
implies that
the right vertical arrow is an isomorphism.
Since 
$g^!\Lambda$ is isomorphic to 
$\Lambda(d)[2d]$,
condition (2) is equivalent to the condition that
the lower horizontal arrow is an isomorphism.
Since condition (1) is equivalent to the condition that
the upper horizontal arrow is an isomorphism,
the assertion follows.
\qed

}

\begin{lm}\label{lmlccA}
{\rm (\cite[Proposition 2.11]{Artin})}
Let $X$ be a scheme and
${\cal F}$ be a constructible sheaf.
Let $x\to X$ be a geometric point.
Then, the following conditions are equivalent:

{\rm (1)}
${\cal F}$ is locally constant on
a neighborhood of the image of $x$.

{\rm (2)}
For every geometric point
$t$ of the strict localization
$X_x$,
the morphism
${\cal F}_x\to {\cal F}_t$ of
specialization is an isomorphism.
\end{lm}


\proof{
(1)$\Rightarrow$(2):
Since we may replace $X$ by
an \'etale neighborhood of $x$,
we may assume that ${\cal F}$ is
constant. Then the morphism
${\cal F}_x\to {\cal F}_t$ is an isomorphism.

(2)$\Rightarrow$(1):
By taking the
cohomology sheaves, we may assume that
${\cal F}$ is concentrated in degree 0.
Since ${\cal F}_x$ is finite,
after replacing $X$ by
an \'etale neighborhood of $x$,
we may assume that there is a morphism
${\cal C}\to {\cal F}$
from a constant sheaf such that
${\cal C}_x\to {\cal F}_x$
is an isomorphism.
Then, for every geometric point
$t$ of the strict localization
$X_x$,
the morphism
${\cal C}_t\to {\cal F}_t$
is also an isomorphism.
Since ${\cal F}$ and ${\cal C}$
are constructible,
the subset $Z$ of $X$
consisting of the image
of geometric points $z$
such that ${\cal C}_z\to {\cal F}_z$
is an isomorphism
is constructible as in
the proof of 
\cite[Proposition 2.11]{Artin}.
By this and (2),
$Z$ contains an open neighborhood $U$
of the image of $x$
and the restriction 
${\cal F}|_U$ is constant.
\qed

}

\begin{lm}\label{lmRjS}
Let $X={\rm Spec}\, {\cal O}_K$
for a henselian discrete valuation
ring and $j\colon U={\rm Spec}\, K
\to X$ be the open immersion.
Let ${\cal F}$ be a constructible sheaf on $X$
and assume that
for every integer $q$,
the action of the inertia subgroup
$I\subset G_K={\rm Gal}(\overline K/K)$
on ${\cal H}^q{\cal F}_{\overline K}$ is
trivial.
Then the following conditions are equivalent:

{\rm (1)}
${\cal F}$ is locally constant.

{\rm (2)}
The morphism
${\cal F}\otimes j_*\Lambda
\to j_*j^*{\cal F}$
is an isomorphism.
\end{lm}

\proof{
Since the implication (1)$\Rightarrow$(2)
is trivial, we show (2)$\Rightarrow$(1).
Assume that ${\cal F}$ is acyclic outside
$[a,b]$ and we show it by induction on $b$.
Let $\overline s$ and $\overline \eta$
be the closed and generic geometric points of
$X$ and let $I$ denote the inertia group.
Then, the stalk
$(R^1j_*j^*{\cal H}^q{\cal F})_{\overline s}$
is canonically isomorphic to the Tate twist
${\cal H}^q{\cal F}_{\overline \eta I}(-1)$
of the coinvariant
with respect to the action of the inertia group $I$.
Since the cohomological dimension of
$Rj_*$ is $1$
and the inertia action is assumed to be trivial,
we have an isomorphism
${\cal H}^b({\cal F})_{\overline s}(-1)
\to 
{\cal H}^b({\cal F})_{\overline \eta}(-1)$.
Thus ${\cal H}^b({\cal F})$ is locally constant
and the assertion follows by induction.
\qed

}

\begin{pr}\label{prlcc}
Let $X$ be a regular excellent
noetherian scheme
and $Z\subset X$ be a closed subset.
For a constructible sheaf ${\cal F}$ on $X$,
the following conditions are equivalent:

{\rm (1)}
${\cal F}$ is locally constant
on a neighborhood of $Z$.

{\rm (2)}
Every separated morphism $h\colon
W\to X$ of finite type
of regular noetherian schemes
is ${\cal F}$-transversal
on a neighborhood of the inverse image of $Z$.
\end{pr}

\proof{
Since the implication (1)$\Rightarrow$(2)
is trivial, we show (2)$\Rightarrow$(1).
Let $x$ be a geometric point of
$Z$ and let $t$ be a geometric point of
the strict localization $X_x$.
By Lemma \ref{lmlccA},
it suffices to show that
the morphism
${\cal F}_x\to {\cal F}_t$
of specialization is an isomorphism.
We may assume that the image  $t_0\in X$
of
$t$ is different from the image $x_0$
of $x$.
Let $T\subset X$ be the closure of 
$t_0$ regarded as a reduced closed 
subscheme and
let $T'\to T$ be the blow-up
at the closure $S$ of $x_0$.

Let $W$ be the normalization of
$T'$ in a finite extension of
the function field such that
the monodromy action on the stalk
of the pull-back is trivial
and $D$ be the reduced inverse image
of $S\subset T$.
By replacing $W$
by an open neighborhood of
the generic point $s$ of
an irreducible component $D$,
we may assume that
$W$ is regular and
$D\subset W$ is a regular divisor.
Then, 
the closed immersion
$i\colon D\to W$ 
is ${\cal F}_W$-transversal
for the pull-back ${\cal F}_W$ of ${\cal F}$
by Lemma \ref{lmFtrtr}.
Hence for the open immersion
$j\colon U=W\sm D\to W$,
the morphism
${\cal F}_W\otimes Rj_*\Lambda
\to Rj_*j^*{\cal F}_W$
is an isomorphism
by Lemma \ref{lmFtrans}.5.
Thus the pull-back of 
${\cal F}_W$ on the local ring
${\rm Spec}\, {\cal O}_{W,s}$
is locally constant by Lemma \ref{lmRjS} and
${\cal F}_x\to {\cal F}_t$
is an isomorphism.
\qed

\medskip
Let $g\colon X'\to X$ be a morphism
of noetherian schemes and $Z'\subset X'$ 
be a closed subset.
We say that $g$ is proper (resp.~finite)
on $Z$, if the restriction of $g$ to $Z$
is proper (resp.~finite) with respect
to a closed subscheme structure on $Z$.

\begin{lm}\label{lmFtrpr}
Let $$\begin{CD}
W'@>{h'}>>X'\\
@V{g'}VV@VVgV\\
W@>h>>X
\end{CD}$$
be a cartesian
diagram of separated morphisms
of finite type of
noetherian schemes
and let ${\cal F}'$ be a constructible sheaf on $X'$.
Assume that $g$ is proper 
on the support of ${\cal F}$ and that
the base change morphism
$g'^*h^!\Lambda\to h'^!\Lambda$ 
{\rm (\ref{eqbcdef})} is
an isomorphism.
We consider the following conditions:

{\rm (1)} $h$ is $g_*{\cal F}'$-transversal.

{\rm (2)} $h'$ is ${\cal F}'$-transversal.

We have the implication {\rm (2)}$\Rightarrow${\rm (1)}.
If $g$ is finite
on the support of ${\cal F}'$, 
then these conditions are equivalent.
\end{lm}

\proof{
In the commutative diagram
(\ref{eqprbc}),
the left vertical arrow is
an isomorphism
by the projection formula,
the proper base change theorem
and the assumption that
$g'^*h^!\Lambda\to h'^!\Lambda$ 
{\rm (\ref{eqbcdef})} is
an isomorphism.
The right vertical arrow is
an isomorphism
by the proper base change theorem.
If $h'$ is ${\cal F}'$-transversal,
then the lower horizontal arrow
is an isomorphism and
hence $h$ is $g_*{\cal F}'$-transversal.

Conversely,
if $h$ is $g_*{\cal F}'$-transversal,
the horizontal arrows are isomorphism.
Further if $g$ is finite
on the support of ${\cal F}'$,
then the condition that
$g_*(c_{h',{\cal F}'})$ is an isomorphism
implies that
$c_{h',{\cal F}'}$ itself is an isomorphism.
\qed

}

\subsection{${\cal F}$-transversality
and local acyclicity}\label{ssFacyc}

\begin{df}\label{dfla}
{\rm (\cite[Definition 2.12]{TF})}
Let $f\colon X\to Y$ be a morphism
and ${\cal F}$
be a sheaf on $X$.

{\rm 1.}
We say that $f$ is 
locally acyclic relatively
to ${\cal F}$ or
${\cal F}$-acyclic for short
if the following condition
is satisfied:

Let $y$ be any geometric point
of $Y$ and
let $t\to Y_{(y)}$ be
any geometric point of
the strict localization $Y_{(y)}$
of $Y$ at $y$.
Let $i\colon X_y=X\times_Yy
\to X\times_YY_{(y)}$
and $j\colon X_t=X\times_Yt
\to X\times_YY_{(y)}$
be the base changes of the canonical
morphisms $y\to Y_{(y)}\gets t$.
Then,
the morphism
\begin{equation}
{\cal F}|_{X_y}\to 
i^*j_*{\cal F}|_{X_t}
\label{eqdfla}
\end{equation}
is an isomorphism.

Let $Z$ be a closed subset of
$X$, we say that $f$ is ${\cal F}$-acyclic
along $Z$ if {\rm (\ref{eqdfla})}
is an isomorphism on the inverse image of $Z$.

{\rm 2.}
We say that $f$ is universally
locally acyclic
if for any morphism $Y'\to Y$
and the pull-back ${\cal F}'$
of ${\cal F}$ on $X'=X\times_YY'$,
the base change
$f'\colon X'\to Y'$
is ${\cal F}'$-acyclic.
\end{df}

We consider the cartesian diagram
\begin{equation}
\begin{CD}
X\times_YY_{(y)}@<j<< X_t\\
@V{f_{(y)}}VV@VV{f_t}V\\
Y_{(y)}@<{j_Y}<< t.
\end{CD}
\label{eqjYt}
\end{equation}
Then, the morphism
(\ref{eqdfla}) is the pull-back by $i
\colon X_y\to X\times_YY_{(y)}$ of
the morphism
\begin{equation}
{\cal F}|_{X\times_YY_{(y)}}
\otimes f_{(y)}^*j_{Y*}\Lambda
\to
j_*j^*({\cal F}|_{X\times_YY_{(y)}})
\label{equlaj}
\end{equation}
defined as the adjoint of
$j^*{\cal F}|_{X\times_YY_{(y)}}
\otimes j^*f_{(y)}^*j_{Y*}\Lambda
\to
j^*{\cal F}|_{X\times_YY_{(y)}}
\otimes f_t^*j_Y^*j_{Y*}\Lambda
\overset {1\otimes{\rm adj}}\to
j^*({\cal F}|_{X\times_YY_{(y)}}).$

\begin{thm}\label{thmgla}
{\rm (cf.~\cite[Th\'eor\`eme 2.12]{TF})}
Let $k$ be a field
and $X$ be a scheme of finite type
over $k$.
Then the morphism
$X\to {\rm Spec}\, k$ 
is universally ${\cal F}$-acyclic for
any sheaf ${\cal F}$ on $X$.
\end{thm}

\proof{
If ${\cal F}$ is constructible,
then this is a special case of
\cite[Th\'eor\`eme 2.12]{TF}.
The general case follows
by taking the limit.
\qed

}

\begin{lm}\label{lmula}
{\rm (cf.~\cite[Proposition 7.6.3]{Fu})}
In the definition of
universal local acyclicity,
it suffices to consider
$Y'$ smooth over $Y$.
\end{lm}

\proof{
Since the question is local on $Y'$,
we may assume that
$Y$ and $Y'$ are affine
and $Y'$ is a closed subscheme
of ${\mathbf A}^n_Y$.
Since $Y'\to {\mathbf A}^n_Y$
is an immersion,
the local acyclicity of
${\mathbf A}^n_X\to {\mathbf A}^n_Y$
implies
the local acyclicity of
$X'\to Y'$.
\qed

}

\begin{pr}\label{prla0}
{\rm (cf.~\cite[Theorem 7.6.9]{Fu})}
Let $f\colon X\to Y$ be a morphism
of schemes
and ${\cal F}$
be a sheaf on $X$.
The following conditions are equivalent:

{\rm (1)} $f$ is ${\cal F}$-acyclic.

{\rm (2)} 
For every cartesian diagram
\begin{equation}
\begin{CD}
X@<h<<W\\
@VfVV@VVgV\\
Y@<i<<V
\end{CD}
\label{eqfghi}
\end{equation}
such that $i$ is an immersion
and for every sheaf ${\cal G}$ on $V$,
the morphism
\begin{equation}
{\cal F}\otimes f^*i_*{\cal G}
\to 
h_*(h^*{\cal F}\otimes g^*{\cal G})
\label{eqlaG}
\end{equation}
is an isomorphism.

{\rm (2$'$)}
For every cartesian diagram
{\rm (\ref{eqfghi})}
such that $i$ is an immersion
and for every locally constant constructible
sheaf ${\cal L}$ on $V$,
the morphism
\begin{equation}
{\cal F}\otimes f^*i_*{\cal L}
\to 
h_*(h^*{\cal F}\otimes g^*{\cal L})
\label{eqlaL}
\end{equation}
is an isomorphism.

If $Y$ is an excellent
noetherian scheme,
we may assume that $V$ is regular
in {\rm (2$'$)}.
\end{pr}

The morphism (\ref{eqlaG})
is defined as the adjoint of
$1\otimes {\rm adj}\colon
h^*({\cal F}\otimes f^*i_*{\cal G})=
h^*{\cal F}\otimes h^*f^*i_*{\cal G}
=h^*{\cal F}\otimes g^*i^*i_*{\cal G}
\to 
h^*{\cal F}\otimes g^*{\cal G}$.

\proof{
(1)$\Rightarrow$(2):
The immersion $i$ is the composition
of an open immersion
and a closed immersion.
For an open immersion,
the isomorphism (\ref{eqlaG})
is \cite[Proposition 2.10]{App},
\cite[Lemma 7.6.7]{Fu}.
If $i$ is a closed immersion,
the isomorphism (\ref{eqlaG})
follows from the projection
formula and the proper base change theorem.

(2)$\Rightarrow$(2$'$) is clear.

(2$'$)$\Rightarrow$(1):
Let $y$ be a geometric point of $Y$
and $t$ be a geometric point of
$Y_{(y)}$.
We show that
(\ref{eqdfla}) is an isomorphism.
Let $v$ be the image of $t$
by $Y_{(y)}\to Y$
and $V\subset Y$ be the closure of $\{v\}$.
Then, there exists 
a projective system
$p_\lambda\colon 
U_\lambda\to V_\lambda$ of finite
\'etale morphisms of integral schemes
and 
open immersions $j_\lambda
\colon V_\lambda\to V$
and an isomorphism
$\varinjlim (j_\lambda
p_\lambda)_*\Lambda
\to j_*\Lambda$.
If $Y$ is excellent,
we may assume that $V_\lambda$
are regular.
Let $j\colon t\to V$
be the canonical morphism
and consider the cartesian diagram
$$
\begin{CD}
X@<h<<W
@<{h_\lambda}<<W\times_VV_\lambda
@<{k_\lambda}<<W\times_VU_\lambda\\
@VfVV@VVgV@VV{g_\lambda}V
@VVV\\
Y@<i<<V@<{j_\lambda}<<V_\lambda
@<{p_\lambda}<<U_\lambda
\end{CD}
$$
extending (\ref{eqfghi}).
As the isomorphism
(\ref{eqlaL}) for the immersions
$ij_\lambda\colon V_\lambda\to Y$
and the locally constant sheaves
$p_{\lambda*}\Lambda$ on
$V_\lambda$,
we obtain an isomorphism
\begin{equation}
{\cal F}\otimes f^*i_*j_{\lambda*}p_{\lambda*}\Lambda
\to 
(hh_\lambda)_*((hh_\lambda)^*{\cal F}
\otimes g_\lambda^*p_{\lambda*}\Lambda)
\label{eqjt}
\end{equation}
on $X$.
The target of
(\ref{eqjt}) is identified with
$(hh_\lambda k_\lambda)_*
(hh_\lambda k_\lambda)^*{\cal F}$
by the projection formula.

We consider the cartesian diagram
$$\begin{CD}
Y@<i<< V@<{j_\lambda}<< 
V_\lambda @<{p_\lambda}<<
U_\lambda @<<<t\\
@AAA@AAA@AAA@AAA@AAA\\
Y_{(y)}@<{i_{(y)}}<<
V_{(y)}@<{j_{\lambda(y)}}<<
V_{\lambda(y)}@<{p_{\lambda(y)}}<<
U_{\lambda(y)}@<<<
t_{(y)}
\end{CD}$$
and let $_{(y)}$ denote the base change
by $Y\gets Y_{(y)}$.
Then the pull-back of (\ref{eqjt})
to $X\times_YY_{(y)}$ gives an isomorphism
\begin{equation}
{\cal F}|_{X\times_YY_{(y)}}
\otimes f_{(y)}^*(i_{(y)}
p_{\lambda(y)}j_{\lambda(y)})_*\Lambda
\to 
(hh_\lambda k_\lambda)_{(y)*}
(hh_\lambda k_\lambda)_{(y)}^*
{\cal F}|_{X\times_YY_{(y)}}.
\label{eqjty}
\end{equation}
Let $T_{\lambda(y)}\subset U_{\lambda(y)}$
be the connected component containing
the image of $t\to U_{\lambda(y)}$
induced by $t\to Y_{(y)}$
and consider the cartesian diagram
$$
\begin{CD}
X\times_YY_{(y)}@<{h_{T\lambda(y)}}<< X\times_YT_{\lambda(y)}\\
@V{f_{(y)}}VV@VV{f_{T_{\lambda(y)}}}V\\
Y_{(y)}@<{j_{T\lambda(y)}}<< T_{\lambda(y)}
\end{CD}$$
whose limit gives (\ref{eqjYt}).
Then, we obtain
an isomorphism
\begin{equation}
{\cal F}|_{X\times_YY_{(y)}}
\otimes f_{(y)}^*j_{T\lambda(y)*}\Lambda
\to 
h_{T\lambda(y)*}*{h_{T\lambda(y)}}^*
{\cal F}|_{X\times_YY_{(y)}}.
\label{eqjtT}
\end{equation}
as a direct summand of
(\ref{eqjty}).
Since $t=\varprojlim T_{\lambda(y)}$,
by taking the limit of (\ref{eqjtT}),
we obtain
%
%
the isomorphism (\ref{equlaj}).
\qed

}

\begin{cor}\label{corgfla}
Let $f\colon X\to Y$ and 
$g\colon Y\to Z$ be morphisms
of schemes and 
let ${\cal F}$ be a sheaf on $X$.

{\rm 1.
(\cite[Proposition 1.4]{Ext}, 
\cite[Corollaire 2.7]{App})}
Let ${\cal G}$ be a sheaf on $Y$
and 
assume that $f$ is ${\cal F}$-acyclic and
that $g$ is ${\cal G}$-acyclic.
Then, $gf$ is
${\cal F}\otimes f^*{\cal G}$-acyclic.

{\rm 2. (\cite[Lemma 3.9 (i)]{SS})}
Assume that $f$ is proper and
that $gf$ is
${\cal F}$-acyclic.
Then, $g$ is
$f_*{\cal F}$-acyclic.
\end{cor}

\proof{
Let $i\colon V\to Z$
be an immersion
and consider the cartesian diagram
$$
\begin{CD}
W@>{f'}>> W'@>{g'}>> V\\
@VhVV@V{h'}VV@VViV\\
X@>f>>Y@>g>>Z.
\end{CD}$$

1.
For any sheaf ${\cal H}$ on $V$,
we have a commutative diagram
$$
\xymatrix{
{\cal F}\otimes f^*{\cal G}
\otimes (gf)^*i_*{\cal H}
\ar[rr]\ar[d]
&&
h_*(h^*({\cal F}\otimes  f^*{\cal G})
\otimes (g'f')^*{\cal H}))
\ar[d]
\\
{\cal F}\otimes f^*
({\cal G}
\otimes g^*i_*{\cal H})
\ar[r]
&
{\cal F}\otimes f^*
(h'_*(h'^*{\cal G}
\otimes g'^*{\cal H}))
\ar[r]
&
h_*(h^*{\cal F}\otimes f'^*
(h'^*{\cal G}
\otimes g'^*{\cal H})).
}$$
The lower arrows are isomorphisms
by the assumption and
Proposition \ref{prla0} (1)$\Rightarrow$(2).
The vertical arrows
are the canonical isomorphisms.
Hence the upper line is an isomorphism
and the assertion follows
by
Proposition \ref{prla0} (2)$\Rightarrow$(1).

2.
For any sheaf ${\cal H}$ on $V$,
we have a commutative diagram
$$
\xymatrix{
f_*{\cal F}\otimes g^*i_*{\cal H}
\ar[rr]\ar[d]
&&
h'_*(h'^*f_*{\cal F}\otimes g'^*{\cal H})
\ar[d]
\\
f_*({\cal F}\otimes (gf)^*i_*{\cal H})
\ar[r]
&
f_*(h_*(h^*{\cal F}
\otimes (g'f')^*{\cal H}))
\ar[r]
&
h'_*(f'_*(h^*{\cal F}
\otimes (g'f')^*{\cal H})).
}$$
The lower left horizontal arrow is an isomorphism
by the assumption and
Proposition \ref{prla0} (1)$\Rightarrow$(2).
The lower right horizontal arrow is 
the canonical isomorphism.
The vertical arrows
are the isomorphisms of projection formula
together with the proper base change theorem.
Hence the upper line is an isomorphism
and the assertion follows
by
Proposition \ref{prla0} (2)$\Rightarrow$(1).
\qed

}

\begin{cor}\label{corFacyctr}
Let 
$$\begin{CD}
X@<h<< W\\
@VfVV@VVV\\
Y@<i<<V
\end{CD}$$ be a
cartesian diagram of schemes
such 
that the vertical arrows are smooth
and that the horizontal arrows are immersions.
Let ${\cal F}$ be a sheaf on $X$ and
${\cal G}$ be a sheaf on $Y$.
Assume that
$f\colon X\to Y$ is ${\cal F}$-acyclic.
If $i\colon V\to Y$ is ${\cal G}$-transversal,
then
$h\colon W\to X$ is ${\cal F}\otimes
f^*{\cal G}$-transversal.

In particular, 
for any immersion $i\colon V\to Y$,
the morphism
$h\colon W\to X$ is ${\cal F}$-transversal.
\end{cor}

\proof{
Since the assertion is local on $V$,
after replacing $Y$ by an open neighborhood
of $V$, we may assume that
$i$ is a closed immersion.
Let $U_0=Y\sm V$ be the complement
and consider the cartesian diagram
$$\begin{CD}
W@>h>> X@<j<< U\\
@VgVV @VfVV @VV{f_U}V\\
V@>i>> Y@<{j_0}<< U_0.
\end{CD}$$
We consider the commutative
diagram
$$\xymatrix{
{\cal F}\otimes f^*
({\cal G}
\otimes
j_{0*}\Lambda)
\ar[d]\ar[r]&
{\cal F}\otimes f^*
j_{0*}j_0^*{\cal G}
\ar[r]&
j_*(j^*{\cal F}\otimes f_U^*j_0^*{\cal G})
\ar[d]
\\
{\cal F}\otimes f^*{\cal G}
\otimes
j_*\Lambda
\ar[rr]
&& j_*j^*({\cal F}\otimes f^*{\cal G})
}$$
Since $f$ is smooth,
$f$ is $\Lambda$-acyclic
by the local acyclicity of
smooth morphisms.
Since
$f$ is ${\cal F}$-acyclic
and $\Lambda$-acyclic,
the upper right horizontal arrow
and the left vertical arrow
are isomorphisms
by Proposition \ref{prla0}
(1)$\Rightarrow$(2).
The right vertical arrow is a canonical isomorphism.
Since $i$ is ${\cal G}$-transversal,
the upper left horizontal arrow
 is an isomorphism
by Lemma \ref{lmFtrans}.5
(1)$\Rightarrow$(2).
Hence the lower horizontal arrow
is an isomorphism
and the assertion follows
by Lemma \ref{lmFtrans}.5
(2)$\Rightarrow$(1).
\qed

}

\begin{cor}\label{corLZ}
{\rm (\cite[Proposition 2.2]{LZ})}
Assume that
$X\to S$ is universally
${\cal F}$-acyclic.
Then, 
for any immersion
$i\colon V\to Y$
of schemes over $S$
and for any sheaf ${\cal G}$
on $V$,
the morphism
${\cal F}\boxtimes
i_*{\cal G}
\to 
(1\times i)_*({\cal F}\boxtimes
{\cal G})$
is an isomorphism.
\end{cor}

\proof{
We apply
Proposition \ref{prla0}
(1)$\Rightarrow$(2)
to the cartesian square
$$\begin{CD}
X\times_S Y@<{1\times i}<< X\times_S V\\
@V{{\rm pr}_2}VV@VVV\\
Y@<i<<V
\end{CD}$$
and ${\rm pr}_1^*{\cal F}$
on $X\times_S Y$ and
${\cal G}$ on $V$.
Then since
${\rm pr}_2\colon
X\times_SY\to Y$ is
${\rm pr}_1^*{\cal F}$-acyclic,
the morphism
$${\cal F}\boxtimes i_*{\cal G}
={\rm pr}_1^*{\cal F}
\otimes
{\rm pr}_2^*i_*{\cal G}
\to 
(1\times i)_*(
(1\times i)^*{\rm pr}_1^*{\cal F}
\otimes
{\rm pr}_2^*{\cal G})
=
(1\times i)_*({\cal F}\boxtimes {\cal G})$$
is an isomorphism.
\qed

}

\begin{pr}[{cf. \cite[Corollary 8.12]{CC}
\cite[Proposition 1.2.4]{ANT})}]\label{prFtracyc}
Let $f\colon X\to Y$ be a smooth morphism
of excellent noetherian schemes
and ${\cal F}$ be a sheaf on $X$.
Let $Z\subset X$ be a closed subset.
Assume that for every cartesian diagram
\begin{equation}
\begin{CD}
X@<<< X'@<{h'}<< W\\
@VfVV@VVV@VVV\\
Y@<<<Y'@<<<V\\
\end{CD}
\label{eqFtracyc}
\end{equation}
where $V\to Y'$ is an immersion
of regular schemes of finite type
over $Y$,
the morphism $h'$
is ${\cal F}'$-transversal
for the pull-back of ${\cal F}'$ on $X'$
on a neighborhood of the inverse image of $Z$.
Then
$f$ is ${\cal F}$-acyclic along $Z$.
\end{pr}

The same proof as loc.~cit.~works
with a form of alteration
due to Gabber
\cite[Expos\'e VII, Th\'eor\`eme 1.1]{Gabber}.

\begin{pr}\label{prBG}
{\rm (cf.~\cite[Appendix Theorem B.2]{Eis})}
Let $f\colon X\to Y$ be a 
morphism of schemes
over an excellent regular
noetherian
scheme $S$
and let ${\cal F}$
be a sheaf on $X$.
If $X\to S$
is universally ${\cal F}$-acyclic and if
$Y$ is smooth over $S$,
then
the following conditions are equivalent:

{\rm (1)} $f$ is 
${\cal F}$-acyclic.

{\rm (2)} 
For any constructible
sheaf ${\cal G}$ on $Y$,
the morphism
$\gamma=(1_X,f)\colon X\to X\times_S Y$ is
${\cal F}\boxtimes {\cal G}$-transversal.
\end{pr}

The assumption that
$X\to S$ is universally ${\cal F}$-acyclic
is satisfied if $S={\rm Spec}\, k$
for a field $k$
by Theorem \ref{thmgla}.

\proof{
%
%
%
%
%
%
By devissage,
the condition (2) is equivalent to 
the condition where
${\cal G}$
is restricted to
${\cal G}=i_*{\cal L}$
for immersions
$i\colon V\to Y$
of regular subschemes
and locally constant constructible
sheaves ${\cal L}$ on $V$
as in Proposition \ref{prla0} (2$'$)
since $Y$ is excellent.

Let $h\colon W=X\times_YV\to X$
be the base change of $i\colon V\to Y$
and consider
the cartesian diagram
\begin{equation}
\begin{CD}
W@>{\gamma'}>> X\times_S V
@>{{\rm pr}_2}>> V\\
@VhVV@VV{1\times i}V
@VViV\\
X@>{\gamma}>> X\times_S Y
@>{{\rm pr}_2}>> Y.
\end{CD}
\label{eqBG1}
\end{equation}
The condition
(2$'$) in Proposition \ref{prla0}
means that the morphism
\begin{equation*}
\gamma^*({\cal F}\boxtimes i_*{\cal L})
={\cal F}\otimes f^*i_*{\cal L}
\to
h_*(h^*{\cal F}\otimes g^*{\cal L})
=
h_*\gamma'^*({\cal F}\boxtimes {\cal L})
\end{equation*}
is an isomorphism.
Hence this for every 
$V$ and ${\cal L}$ as above
is equivalent to (1)
by Proposition \ref{prla0}
(2$'$)$\Leftrightarrow$(1).
Thus it suffices to show that
the following conditions are equivalent:

(1$'$)
For every immersion $i\colon V\to Y$
of regular subscheme
and every locally constant constructible sheaf
${\cal L}$ on $V$, the morphism
\begin{equation}
\gamma^*({\cal F}\boxtimes i_*{\cal L})
\to
h_*\gamma'^*({\cal F}\boxtimes {\cal L})
\label{eqBG0}
\end{equation}
is an isomorphism.

(2$'$)
For every immersion $i\colon V\to Y$
of regular subscheme
and every locally constant constructible sheaf
${\cal L}$ on $V$, the morphism
\begin{equation}
c_{\gamma, {\cal F}\boxtimes i_*{\cal L}}\colon
\gamma^*({\cal F}\boxtimes i_*{\cal L})
\otimes \gamma^!\Lambda
\to
\gamma^!({\cal F}\boxtimes i_*{\cal L})
\label{eqBG00}
\end{equation}
is an isomorphism.

We construct
a commutative diagram
\begin{equation}
\begin{CD}
\gamma^*({\cal F}\boxtimes i_*{\cal L})
\otimes \gamma^!\Lambda
@>{c_{\gamma, {\cal F}\boxtimes i_*{\cal L}}}>>
\gamma^!({\cal F}\boxtimes i_*{\cal L})\\
@VVV@VVV\\
h_*(\gamma'^*({\cal F}\boxtimes {\cal L})
\otimes h^*\gamma^!\Lambda)
@>>>
h_*\gamma'^!({\cal F}\boxtimes {\cal L})
\end{CD}
\label{eqBG}
\end{equation}
similarly as (\ref{eqprbc}).
The lower horizontal
arrow is obtained by
applying $h_*$
to the composition 
\begin{equation}
\gamma'^*({\cal F}\boxtimes {\cal L})
\otimes h^*\gamma^!\Lambda
\to
\gamma'^!({\cal F}\boxtimes {\cal L})
\label{eqBG2}
\end{equation}
of
$c_{\gamma',
{\cal F}\boxtimes {\cal L}}$
with the morphism
induced by
$h^*\gamma^!\Lambda
\to \gamma'^!\Lambda$
(\ref{eqbcdef}).
The right vertical arrow
is the composition of
$\gamma^!({\cal F}\boxtimes i_*{\cal L})
\to 
\gamma^!(1\times i)_*
({\cal F}\boxtimes {\cal L})$
and the base change isomorphism
$\gamma^!(1\times i)_*\to 
h_*\gamma'^!$
(\ref{eqbcdef2}).
By the assumption that
the morphism
$X\to S$
is universally ${\cal F}$-acyclic
and 
by Corollary \ref{corLZ},
the morphism
${\cal F}\boxtimes i_*{\cal L}
\to 
(1\times i)_*
({\cal F}\boxtimes {\cal L})$
is an isomorphism.
Hence the right vertical arrow is an isomorphism.
The left vertical arrow
is the adjoint of the isomorphism
$h^*\gamma^*({\cal F}\boxtimes i_*{\cal L})
\otimes h^*\gamma^!\Lambda
\to
\gamma'^*({\cal F}\boxtimes {\cal L})
\otimes h^*\gamma^!\Lambda$.
By the assumption that
$Y\to S$ is smooth,
we have an isomorphism
$\Lambda(-d)[-2d]\to \gamma^!\Lambda$
for the relative dimension
$d$ of $Y$ over $S$.
Hence
the left vertical arrow is an isomorphism
if and only if (\ref{eqBG0})
is an isomorphism.
By applying Lemma \ref{lmbc}.2
to the left square of 
(\ref{eqBG1}),
we see that 
(\ref{eqBG}) is commutative.

We show that the following conditions are equivalent:

(1$^\circ$)
(\ref{eqBG2})
and the left vertical arrow
in (\ref{eqBG})
are isomorphisms.

(2$^\circ$)
(\ref{eqBG00}) is an isomorphism,

(3$^\circ$)
All arrows in
(\ref{eqBG})
are isomorphisms.

(1$^\circ$)$\Rightarrow$(3$^\circ$):
If (\ref{eqBG2})
is an isomorphism, then
the lower horizontal arrow
in (\ref{eqBG})
is an isomorphism.
Further if the left vertical arrow
in (\ref{eqBG})
is an isomorphism,
(3$^\circ$) is satisfied.

(2$^\circ$)$\Rightarrow$(1$^\circ$):
If (\ref{eqBG00}) is an isomorphism,
then the composition of
(\ref{eqBG}) via upper right
is an isomorphism.
Since $h$ is an immersion,
this implies that (\ref{eqBG2})
is an isomorphism.
The rest is similar to the proof of
(1$^\circ$)$\Rightarrow$(3$^\circ$).

(3$^\circ$)$\Rightarrow$(2$^\circ$)
is clear.

Thus it suffices
to show that
one of 
(\ref{eqBG2}) and  (\ref{eqBG00}) 
is an isomorphism
assuming (1$'$).
First, we prove (\ref{eqBG2}) in
the case where
$i$ is an open immersion and ${\cal L}=\Lambda$.
Since $Y$ is smooth over $S$,
the morphism $\gamma$
is a section of a smooth morphism
${\rm pr}_1\colon X\times_S Y\to X$.
Hence $\gamma$
is ${\cal F}\boxtimes
\Lambda$-transversal
by Lemma \ref{lmFtrans}.1 and Lemma \ref{lmFtrtr}.
Further, its restriction
$\gamma'\colon
W\to X\times_S V$
is also ${\cal F}\boxtimes
\Lambda$-transversal.
This means that
$
\gamma'^*({\cal F}\boxtimes
\Lambda)
\otimes \gamma'^!\Lambda
\to
\gamma'^!
({\cal F}\boxtimes
\Lambda)
$
(\ref{eqBG2}) 
is an isomorphism
and implies that $\gamma$
is ${\cal F}\boxtimes
i_*\Lambda$-transversal
as we have shown above.

Next, we prove (\ref{eqBG00}) in the case
where $i$ is a closed immersion
and ${\cal L}=i^!\Lambda$.
Note that if $i$ is a regular immersion of codimension
$c$, then
$i^!\Lambda=\Lambda(-c)[-2c]$ is locally constant
constructible by purity \cite[Expos\'e XVI, Th\'eor\`eme 3.1.1]{Gabber}.
Let $j\colon Y\sm V\to Y$ be the open immersion 
of the complement
and consider the distinguished triangle
$
{\cal F}\boxtimes
i_*i^!\Lambda
\to
{\cal F}\boxtimes
\Lambda
\to
{\cal F}\boxtimes
j_*\Lambda\to.$
As
we have already shown above,
$\gamma\colon X\to X\times_S Y$
is transversal for the last two terms.
Hence by Lemma \ref{lmFtrans}.2, 
the morphism $\gamma$
is also ${\cal F}\boxtimes
i_*i^!\Lambda$-transversal
and 
(\ref{eqBG00}) 
for ${\cal L}=i^!\Lambda$
is an isomorphism.

We show (\ref{eqBG2}) in the general case.
Since the isomorphism
(\ref{eqBG2}) is an
\'etale local condition on $V$
and 
since
${\cal L}$ is assumed locally constant,
it is reduced to the case
where
$i\colon V\to Y$ is a closed immersion
by replacing $Y$ by an open neighborhood of $V$
and 
${\cal L}= i^!\Lambda
=\Lambda(-c)[-2c]$.
\qed

}

\subsection{$C$-transversality}\label{ssCtr}

\begin{df}\label{dfCtrans0}
Let $X$ be a scheme.
Let
$V$ be a vector bundle on $X$
and let ${\mathbf P}(V)
=(V\sm X)/{\mathbf G}_m$
be the associated projective space bundle.

{\rm 1.}
We say that a closed subset
$C\subset V$
is conical if it is stable under the
${\mathbf G}_m$-action.

{\rm 2.}
For a closed conical subset $C\subset V$,
we call
the intersection $C\cap X$ with the $0$-section
$X\subset V$
the base of $C$
and the image ${\mathbf P}(C)$ of
$C\sm C\cap X$ by the projection
$V\sm X\to {\mathbf P}(V)$ the
projectivization of $C$.

{\rm 3.}
Let $f\colon V\to W$ be a linear morphism
of vector bundles on $X$ and
 $C\subset V$ be a closed conical
subset.
For $x\in X$,
we say that $f$ is $C$-transversal at $x$
if the fiber at $x$ of the intersection 
$$C\cap {\rm Ker}(f\colon V\to W)
=C\times_WX$$
is a subset of $\{0\}$.
We say that $f$ is $C$-transversal
if we have the inclusion everywhere on $X$.
\end{df}

The base $C\cap X\subset X$
and the projectivization 
${\mathbf P}(C)\subset {\mathbf P}(V)$
are closed subsets.
The base of $C$ equals the image
of $C$ by the projection $V\to X$.

\begin{lm}\label{lmCC'}
Let $C$ and $C'$ be closed conical subsets
of a vector bundle $V$.
Then the following conditions are equivalent:

{\rm (1)}
$C\subset C'$.

{\rm (2)}
$C\cap X\subset C'\cap X$
and 
${\mathbf P}(C)\subset {\mathbf P}(C')$.
\end{lm}

\proof{
(1)$\Rightarrow$(2) is clear.
Since $C$ is the union of
$C\cap X$
and the inverse image
of ${\mathbf P}(C)$
by $V\sm X\to {\mathbf P}(V)$,
we have
(2)$\Rightarrow$(1).
\qed

}

\begin{lm}\label{lmCtrans}
Let $X$ be a scheme and
$f\colon V\to W$ be a linear morphism
of vector bundles on $X$.
Let $C\subset V$ be a closed conical
subset.

{\rm 1.}
The subset $U=\{x\in X\mid
f$ is $C$-transversal at $x\}\subset X$
is an open subset.

{\rm 2.}
Endow $C\subset V$ with a closed 
subscheme structure
stable under the action of ${\mathbf G}_m$.
Then the following conditions are equivalent:

{\rm (1)}
$f$ is $C$-transversal.

{\rm (2)}
The restriction of $f$ on $C$ is finite.

{\rm (3)}
The fiber product
$C\times_WX$ is finite over $X$.

\noindent
If these equivalent conditions are satisfied,
then the image
$f(C)\subset W$ 
is a closed conical subset.
\end{lm}

\proof{
1.
The intersection
$C'=C\cap {\rm Ker}(f\colon V\to W)$
is a closed conical subset of $V$
and its projectivization
${\mathbf P}(C')$
is a closed subset of the
associated projective space bundle
${\mathbf P}(V)$.
Since $U\subset X$ is the complement
of the image of
${\mathbf P}(C')\subset{\mathbf P}(V)$,
it is an open subset.

2.
(1)$\Leftrightarrow$(3):
Condition (1) means that the underlying set
of $C\times_WX\subset V$ is a subset of 
the $0$-section $X$.
Since the ideal defining the closed
immersion $X\to V$ is locally of finite type,
this implies (3).
Conversely,
since $C$ is stable under the ${\mathbf G}_m$-action,
condition (3) implies that the underlying set
of $C\times_WX\subset V$ is a subset of 
the $0$-section $X$.

(2)$\Leftrightarrow$(3):
Since the assertion is local on $X$,
we may assume that
$W={\rm Spec}\, A$,
$V={\rm Spec}\, B$ are affine
for graded rings $A=\bigoplus_{n\in{\mathbf N}}A_n, 
B=\bigoplus_{n\in{\mathbf N}}B_n$
and $C\subset V$ is defined
by a graded ideal $J=\bigoplus_{n\in{\mathbf N}}J_n
\subset B$.
Hence, the assertion follows from Lemma 
\ref{lmGrAB} below.
\qed

}

\begin{lm}\label{lmGrAB}
Let $A=\bigoplus_{n\in{\mathbf N}}A_n$ be a graded ring.

{\rm 1.}
Let $M=\bigoplus_{n\in{\mathbf N}}M_n$
be a graded $A$-module
and
let $S$ be a subset of $M$
consisting of homogeneous elements.
Then the following are equivalent:

{\rm (1)} The $A$-module
$M$ is generated by $S$.

{\rm (2)} For every $n\geqq 0$,
the $A/A_{\geqslant n}$-module
$M/M_{\geqslant n}$ is generated by
the image of $S$.

{\rm 2.}
Let $A=\bigoplus_{n\in{\mathbf N}}A_n\to 
B=\bigoplus_{n\in{\mathbf N}}B_n$ be a morphism of graded rings
and $J=\bigoplus_{n\in{\mathbf N}}J_n\subset B$ be a graded ideal.
Let $S\subset B$ be a set of homogeneous elements.
Then the following conditions are equivalent:

{\rm (1)} 
The $A$-module $B/J$ is generated by the image
of $S$.

{\rm (2)} 
The $A/A_{\geqslant 1}$-module
$B/(J+A_{\geqslant 1}B)$ is generated by the image
of $S$.
\end{lm}

\proof{
1.
By replacing $M$ by the quotient
generated by $S$,
we may assume that $S=\varnothing$.
Then each condition means that
$M=0$.

2.
(2)$\Rightarrow$(1):
By induction on $n$,
the $A/A_{\geqslant n}$-module
$B/(J+A_{\geqslant n}B)$ is generated by the image
of $S$.
Since $A_{\geqslant n}B\subset B_{\geqslant n}$,
the $A/A_{\geqslant n}$-module
$B/(J+B_{\geqslant n})$ is also generated by the image
of $S$.
Hence it suffices to apply 1
to $M=B/J$.

(1)$\Rightarrow$(2) is clear.
\qed

}

\medskip

Let $p$ be a prime number
and we fix a regular noetherian scheme
$S$ over ${\mathbf Z}_{(p)}$.
satisfying the following finiteness condition:
\medskip

\noindent
(F)
The reduced part of
$S_{{\mathbf F}_p}
=S\times_{{\rm Spec}\, {\mathbf Z}_{(p)}}
{\rm Spec}\, {\mathbf F}_p$
is a scheme of finite type
over a field of characteristic $p$
with finite $p$-basis.

\medskip
\noindent
Let $X$ be
a regular scheme of finite type
over $S$.
Then the sheaf $F\Omega_X$
of Frobenius--Witt differentials
is a locally free 
${\cal O}_{X_{{\mathbf F}_p}}$-module 
of finite rank
\cite[Theorem 3.4]{FW}
and 
define a vector bundle
$FT^*X$ on $X_{{\mathbf F}_p}$.
For every point $x$ of $X_{{\mathbf F}_p}$,
we have a short exact sequence
$$0\to
F^*({\mathfrak m}_x/
{\mathfrak m}_x^2)
\to
F\Omega_{X,x}
\otimes_{{\cal O}_{X,x}}k(x)
\to
F^*\Omega_{k(x)/{\mathbf F}_p}
\to 0$$
where the functor $F^*$ denotes
the extension of scalars
by the absolute Frobenius
$k(x)\to k(x)$
[loc.~cit.~Proposition 2.4].
For a morphism $X\to Y$
of regular schemes of finite type
over $S$,
if $X\to Y$ is smooth,
then 
$FT^*Y\times_YX\to
FT^*X$ is injective.
Conversely, 
if $FT^*Y\times_YX\to
FT^*X$ is injective,
then $X\to Y$ is smooth
on a neighborhood of
$X_{{\mathbf F}_p}$ by
[loc.~cit.~Proposition 2.10].
Further, in this case,
if $F^*T^*X/Y|_{X_{{\mathbf F}_p}}$ denote
the Frobenius pull-back of
the restriction to $X_{{\mathbf F}_p}$
of 
the relative cotangent bundle
defined by $\Omega^1_{X/Y}$,
we have a short exact sequence
\begin{equation}
0\to FT^*Y\times_YX\to
FT^*X\to F^*T^*X/Y|_{X_{{\mathbf F}_p}}\to 0
\label{eqFT}
\end{equation}
of vector bundles on $X_{{\mathbf F}_p}$.

\begin{df}\label{dfCtrans}
Let $h\colon W\to X$ be a 
morphism 
of regular schemes of
finite type over $S$.
Let $C\subset FT^*X$ be a
closed conical subset
and let $h^*C$
denote the closed conical subset
$C\times_XW
\subset
FT^*X\times_XW$.

{\rm 1.}
We say that $h$ is
$C$-transversal
if $FT^*X\times_XW\to FT^*W$
is $h^*C$-transversal.

{\rm 2.}
If $h$ is $C$-transversal,
we define a closed conical subset
$h^\circ C\subset FT^*W$
to be the image of
$h^*C$ by $FT^*X\times_XW\to FT^*W$.
\end{df}

The terminology is explained by
Lemma \ref{lmXtr}.
If $h$ is $C$-transversal
and if $C'\subset C$ is another closed
conical subset,
then $h$ is $C'$-transversal.
If $h$ is $C$-transversal,
then $h^\circ C$ is a closed
conical subset by 
Lemma \ref{lmCtrans}.2.

\begin{lm}\label{lmCacyc0}
Let $h\colon W\to X$ be a 
morphism 
of regular schemes of
finite type over $S$.
Let $C\subset FT^*X$ be a
closed conical subset.

{\rm 1.}
If $h$ is smooth,
then $h$ is $C$-transversal.

{\rm 2.}
Assume that $C$ is a closed subset of
the $0$-section $FT^*_XX$.
Then, $h$ is $C$-transversal.

{\rm 3.}
Assume that $C=FT^*X$
and that $h$ is $C$-transversal.
Then $h$ is smooth on
a neighborhood of
$W_{{\mathbf F}_p}$.
\end{lm}

\proof{
1.
Since $FT^*X\times_XW\to FT^*W$
is an injection,
the intersection
$h^*C\cap {\rm Ker}
(FT^*X\times_XW\to FT^*W)
\subset {\rm Ker}
(FT^*X\times_XW\to FT^*W)$
is a subset of the $0$-section.

2.
If 
$C$ is a subset of
the $0$-section $FT^*_XX$,
then $h^*C\cap {\rm Ker}
(FT^*X\times_XW\to FT^*W)
\subset h^*C$
is a subset of the $0$-section.

3.
Since $FT^*X\times_XW\to FT^*W$
is an injection,
$h$ is smooth on
a neighborhood of
$W_{{\mathbf F}_p}$
by \cite[Proposition 2.10]{FW}.
\qed

}

\begin{lm}\label{lmXtr}
Let $h\colon W\to X$ be a morphism
of finite type of regular noetherian schemes.
Let $Z\subset X$ be a regular closed subscheme
and let
$C\subset FT^*X$ be the image
of the Frobenius pull-back
of the conormal bundle
$F^*T^*_ZX|_{Z_{{\mathbf F}_p}}$.
Assume that $X$ is excellent.
Then, the following conditions
are equivalent:

{\rm (1)}
$h$ is $C$-transversal.

{\rm (2)}
$h$ is transversal to the immersion
$Z\to X$
on a neighborhood of $Z_{{\mathbf F}_p}$.
\end{lm}
\proof{
Similarly as in the proof of Lemma \ref{lmpurity},
it is reduced to the case where $h$ is a
closed immersion.
This case follows from the regularity criterion
\cite[Corollary 2.8]{FW}.
\qed

}

\begin{lm}\label{lmCtrtr}
Let $h\colon W\to X'$
and $g\colon X'\to X$ be morphisms of
regular schemes of
finite type over $S$ and
let $C$ be a closed conical subset
of $FT^*X$.
Assume that $g$ is $C$-transversal.
Then, the following conditions are equivalent:

{\rm (1)} $gh$ is $C$-transversal.

{\rm (2)} $h$ is $g^\circ C$-transversal.
\end{lm}

\proof{
Since $g^\circ C$ is defined as the image
of $g^*C$ and $g^*C\to g^\circ C$ is finite,
condition (2) is equivalent to (1)
$h^*g^*C\cap {\rm Ker}(FT^*X\times_XW
\to FT^*W)$ is a subset of the $0$-section.
\qed

}

\begin{df}\label{dff0C}
Let $g\colon X\to Y$ be
a morphism of regular schemes of
finite type over $S$
and let $C\subset FT^*X$
be a closed conical subset.
Assume that $g$ is proper on
the base $B$ of $C$.
Then, we define a closed conical
subset $g_\circ C
\subset FT^*Y$
to be the image 
by $FT^*Y\times_YX\to FT^*Y$
of the inverse image of $C$
by $FT^*Y\times_YX\to FT^*X$.
\end{df}

\begin{lm}\label{lmg0C}
Let $g\colon X'\to X$ be
a morphism of regular schemes of
finite type over $S$
and let $C'\subset FT^*X'$
be a closed conical subset.
Assume that $g$ is proper on
the base $B'$ of $C'$.
Let $h\colon W\to X$
be a morphism of regular schemes of
finite type over $S$
and assume that
$h$ is $g_\circ C'$-transversal.
Let $h'\colon  W'\to X'$
be the base change of $h$.
Assume that $S$ is excellent.

{\rm 1.}
The morphism
$h$ is transversal to $g$ on a neighborhood
of $h'^{-1}(B')$.

{\rm 2.}
There exists a regular neighborhood $U'
\subset W'=W\times_XX'$ of $h'^{-1}(B')$
on which $h'|_{U'}\colon U'\to X'$ is
$C'$-transversal.
\end{lm}

\proof{1.
Since $h\colon W\to X$
is $C$-transversal, on $h'^{-1}(B')$,
the intersection of
the kernel of $(FT^*X\times_XX'\to FT^*X')
\times_{X'}W'$
with 
the kernel $(FT^*X\times_XW\to FT^*W)
\times_WW'$
is a subset of the $0$-section. 
Namely,
the morphism 
$FT^*X\times_XW'\to 
FT^*X'\times_{X'}W'\oplus
FT^*W\times_WW'$
is an injection on $h'^{-1}(B')$.
This implies that
$W'$ is regular on a neighborhood $U'$ of
$h'^{-1}(B')$
by \cite[Theorem 3.4]{FW}
and ${\rm rank}\, FT^*U'-{\rm rank}\, FT^*X'
={\rm rank}\, FT^*W-{\rm rank}\, FT^*X$.

2.
Since the question is local on $W$,
we may decompose $h$ as the composition
of a smooth morphism and an immersion.
Hence by Lemma \ref{lmCtrtr},
we may assume that $h$ is an immersion.
By 1,
$g\colon X'\to X$ is transversal
to the immersion $h\colon W'\to W$
on a neighborhood $U'\subset W'$ of $h'^{-1}(B')$.
Then ${\rm Ker}(FT^*X'\times_{X'}U'\to FT^*U')$
is identified with
${\rm Ker}(FT^*X\times_XW\to FT^*W)\times_WU'$.
Since $g_\circ C'$ is defined as the image
of the intersection of
the inverse image of $C'$ with the image
${\rm Im}(FT^*X\times_XX'\to FT^*X')$
of the injection,
the $C'$-transversality follows.
\qed

}

\subsection{Micro support}\label{ssms}

We keep the notation that
$S$ is a regular noetherian scheme
over ${\mathbf Z}_{(p)}$
satisfying the finiteness condition:

\noindent
(F)
The reduced part of
$S_{{\mathbf F}_p}$
is of finite type over a field of characteristic $p$
with finite $p$-basis.

\begin{df}\label{dfms}
{\rm (\cite[Definition 3.1.1]{ANT})}
Let $X$ be a 
regular scheme of
finite type over $S$.
Let ${\cal F}$ be a constructible
sheaf on $X$
and
$C$ be a closed conical subset
of $FT^*X$.
We say that ${\cal F}$ is micro supported
on $C$ if the following conditions 
{\rm (1)} and {\rm (2)}
are satisfied:

{\rm (1)} The intersection of
the support of ${\cal F}$ and
$X_{{\mathbf F}_p}$ is
a subset of the base of $C$.

{\rm (2)}
Every $C$-transversal
separated morphism
$h\colon W\to X$
of regular schemes
of finite type over $S$
 is ${\cal F}$-transversal
on a neighborhood of
$W_{{\mathbf F}_p}$.
\end{df}

Note that the conditions (1) and (2) are
on neighborhoods of $X_{{\mathbf F}_p}$
and of $W_{{\mathbf F}_p}$.
This reflects the fact that
the FW cotangent bundle 
$FT^*X$ is a vector bundle on $X_{{\mathbf F}_p}$.

\begin{lm}\label{lmmslc}
Let $X$ be a regular scheme of
finite type over $S$
and ${\cal F}$ be a
constructible sheaf on $X$.

{\rm 1.}
${\cal F}$ is micro supported
on $FT^*X$.

{\rm 2.}
We consider the following conditions:

{\rm (1)}
${\cal F}$ is locally constant
on a neighborhood of
$X_{{\mathbf F}_p}$.

{\rm (2)}
${\cal F}$ is micro supported
on the $0$-section $FT^*_XX$.

\noindent
We have
{\rm (1)}$\Rightarrow${\rm (2)}.
Conversely if $X$ is excellent,
we have
{\rm (2)}$\Rightarrow${\rm (1)}.

{\rm 3.}
The following conditions are equivalent:

{\rm (1)}
${\cal F}=0$ 
on a neighborhood of
$X_{{\mathbf F}_p}$.

{\rm (2)}
${\cal F}$ is micro supported
on $\varnothing$.
\end{lm}

\proof{
1.
Let $h\colon W\to X$ be 
a separated morphism of 
regular schemes of finite
type over $S$
and assume that $h$ is 
$FT^*X$-transversal.
Then $h$ is smooth on 
a neighborhood of
$W_{{\mathbf F}_p}$ by Lemma \ref{lmCacyc0}.3
and is ${\cal F}$-transversal
a neighborhood of
$W_{{\mathbf F}_p}$ by Lemma \ref{lmFtrans}.1.
Condition (1) in Definition \ref{dfms} is clearly satisfied.

2.
(1)$\Rightarrow$(2):
If ${\cal F}$ is locally constant
on a neighborhood of
$X_{{\mathbf F}_p}$,
every separated
morphism $h\colon W\to X$
of regular schemes
of finite type over $S$
is ${\cal F}$-transversal 
on a neighborhood of
$W_{{\mathbf F}_p}$
by Lemma \ref{lmFtrans}.3.
Condition (1) in Definition \ref{dfms} is clearly satisfied.

(2)$\Rightarrow$(1):
Assume that
${\cal F}$ is micro supported
on the $0$-section $FT^*_XX$.
Then, every separated
morphism $h\colon W\to X$
of regular schemes
of finite type over $S$
is $FT^*_XX$-transversal 
by Lemma \ref{lmCacyc0}.2
and
hence is ${\cal F}$-transversal
on a neighborhood of
$W_{{\mathbf F}_p}$.
Thus the assertion follows from
Proposition \ref{prlcc}.

3.
By condition (1) in Definition \ref{dfms},
condition (2) is equivalent to that
the support of ${\cal F}$ is
a subset of $\varnothing$ on
a neighborhood of
$X_{{\mathbf F}_p}$.
\qed

}

\begin{lm}\label{lmmstr}
Let $h\colon W\to X$ be a 
separated morphism of
regular schemes of
finite type over $S$.
Let ${\cal F}$ be a constructible
sheaf on $X$
and
$C$ be a closed conical subset
of $FT^*X$.
If ${\cal F}$ is micro supported
on $C$ and if
$h\colon W\to X$ is $C$-transversal,
then $h^*{\cal F}$ is micro supported
on $h^\circ C$.
\end{lm}

\proof{
It suffices
to show that
condition (2) in Definition \ref{dfms} is satisfied.
Since
$h\colon W\to X$ is $C$-transversal,
it is  ${\cal F}$-transversal
on a neighborhood of $W_{{\mathbf F}_p}$.
Let $g\colon V\to W$
be a separated 
$h^\circ C$-transversal
morphism
of regular schemes of finite type
over $S$.
By Lemma \ref{lmCtrtr}
(2)$\Rightarrow$(1),
$hg\colon V\to X$
is $C$-transversal and
hence ${\cal F}$-transversal
on a neighborhood of $V_{{\mathbf F}_p}$.
By Lemma \ref{lmFtrtr}
(1)$\Rightarrow$(2), 
$g\colon V\to W$
is $h^*{\cal F}$-transversal
on a neighborhood of $V_{{\mathbf F}_p}$
since the assumption is satisfied 
by purity \cite[Expos\'e XVI, Th\'eor\`eme 3.1.1]{Gabber}.
\qed

}

\begin{lm}\label{lmmsloc}
Let $X$ be a 
regular scheme of
finite type over $S$.
Let ${\cal F}$ be a constructible
sheaf on $X$
and $C\subset FT^*X$
be a closed conical subset.

{\rm 1.}
Let $X=\bigcup U_i$
be an open covering of $X$.
Then the following conditions are equivalent:

{\rm (1)}
${\cal F}$ is micro supported on $C$.

{\rm (2)}
${\cal F}|_{U_i}$ is micro supported on 
$C|_{U_i}$ for every $i\in I$.

{\rm 2.}
Let $U\subset X$ be an open subscheme
and let $Z=X\sm U$ be the complement.
Let $C'\subset FT^*U$ be
a closed subset.
Assume that ${\cal F}$
is micro supported on $C$
and that ${\cal F}|_U$
is micro supported on $C'$.
Then,  ${\cal F}$
is micro supported on
the union $C_1=\overline {C'}\cup C|_Z$.
\end{lm}

\proof{
1. (1)$\Rightarrow$(2):
Since the open immersion
$j\colon U\to X$
is $C$-transversal
by Lemma \ref{lmCacyc0}.1,
this follows from Lemma \ref{lmmstr}.

(2)$\Rightarrow$(1):
Let $h\colon W\to X$
be a separated $C$-transversal
morphism of
regular schemes of finite type over
$S$.
Then for every $i$,
the base change
$h_i\colon W_i=W\times_XU_i\to U_i$
is $C|_{U_i}$-transversal.
Since ${\cal F}|_{U_i}$
is micro supported on
$C|_{U_i}$,
$h_i$ is  ${\cal F}|_{U_i}$-transversal
on a neighborhood of
$W_{i,{\mathbf F}_p}$.
Hence
$h$ is ${\cal F}$-transversal on
a neighborhood of $W_{{\mathbf F}_p}$.

2.
Let $h\colon W\to X$
be a separated $C_1$-transversal
morphism of
regular schemes of finite type over
$S$.
Then, 
the restriction
$h_U\colon h^{-1}(U)\to U$
is $C'$-transversal
and 
$h\colon W\to X$
is $C$-transversal
on a neighborhood of $h^{-1}(Z)$
by Lemma \ref{lmCtrans}.1.
Hence by 1,
$h_U$ is ${\cal F}|_U$-transversal
on a neighborhood of
$h^{-1}(U)_{{\mathbf F}_p}$
and
$h$ is ${\cal F}$-transversal
on a neighborhood of
$h^{-1}(Z)_{{\mathbf F}_p}$.
Since 
$W=h^{-1}(U)\cup h^{-1}(Z)$,
$h$ is ${\cal F}$-transversal
on a neighborhood of
$W_{{\mathbf F}_p}$.
\qed

}

\begin{lm}\label{lmg0F}
Let $g\colon X'\to X$ be
a morphism of regular schemes of
finite type over $S$
and let $C'\subset FT^*X'$
be a closed conical subset.
Assume that $g$ is proper on
the base $B'$ of $C'$.
Let ${\cal F}'$ be a constructible sheaf on $X'$
and 
assume that ${\cal F}'$ is micro supported on 
$C'\subset FT^*X'$.
Then, $g_*{\cal F}'$ is micro supported on $g_\circ C'$.
\end{lm}

\proof{
It suffices
to show that
condition (2) in Definition \ref{dfms} is satisfied.
Let $h\colon W\to X$ be a morphism
of regular schemes of finite type over $S$
and assume that $h$ is $g_\circ C'$-transversal.
By Lemma \ref{lmg0C}, 
the morphism
$h$ is transversal to $g$ on a neighborhood
$U'\subset W'=W\times_XX'$ of $h'^{-1}(B')$
and $h'\colon U'\to X'$ is
$C'$-transversal.
Since the assumption that
the base change morphism
$g'^*h^!\Lambda\to h'^!\Lambda$ 
{\rm (\ref{eqbcdef})} is
an isomorphism on $U'$
by Lemma \ref{lmpurity}.2,
the morphism $h$ is $g_*{\cal F}'$-transversal
on a neighborhood of $W_{{\mathbf F}_p}$ by
Lemma \ref{lmFtrpr} (2)$\Rightarrow$(1).
\qed

}

\begin{lm}\label{lmSSXP}
Let $i\colon X\to P$ be a closed immersion of
regular schemes of
finite type over $S$
and
let ${\cal F}$ be a 
constructible sheaf on $X$.
Let $C_P$ be a closed conical subset
of $FT^*P|_X
\subset FT^*P$
and define $C\subset FT^*X$ 
to be the closure of
the image of $C_P$
by the surjection
$FT^*P|_X\to FT^*X$.
If $i_*{\cal F}$ is micro supported on $C_P$,
then
${\cal F}$ is micro supported on $C$.
\end{lm}
\proof{
%
It suffices
to show that
condition (2) in Definition \ref{dfms} is satisfied.
Let $h\colon W\to X$ be a 
separated $C$-transversal morphism
of regular schemes of finite type over $S$.
We show that
$h$ is ${\cal F}$-transversal
on a neighborhood of $W_{{\mathbf F}_p}$.
Since the question is local on $W$,
we may assume that there exists
an immersion $W\to {\mathbf A}^n_X$.
By replacing $W\to X\to P$ by 
$W\to {\mathbf A}^n_X\to {\mathbf A}^n_P$
and by Lemma \ref{lmFtrans}.1,
Lemma \ref{lmFtrtr},
Lemma \ref{lmCacyc0}.1
and Lemma \ref{lmCtrtr},
we may assume that $h$ is an immersion.

Since the assertion is local on $W$,
by taking a lifting of a local basis of
the conormal bundle $T^*_WX$,
we may extend the immersion $h\colon W\to X$
to an immersion $g\colon
V\to P$ of regular closed subscheme
transversal with the immersion $X\to P$.
Then,
the $C$-transversality of $h\colon W\to X$
implies the $C_P$-transversality of $V\to P$.
Since $i_*{\cal F}$ is micro supported on $C_P$,
this implies that $g\colon V\to P$ is $i_*{\cal F}$-transversal
on
a neighborhood of $V_{{\mathbf F}_p}$.
Since the intersection $W=V\cap X$ is transversal,
if $i'\colon W\to V$ denote the immersion,
the morphism
$i'^*g^!\Lambda\to h^!\Lambda$
is an isomorphism.
Hence by
Lemma \ref{lmFtrpr} (1)$\Rightarrow$(2),
the $i_*{\cal F}$-transversality of $g$
implies the ${\cal F}$-transversality of $h$
on a neighborhood of $W_{{\mathbf F}_p}$.
Hence ${\cal F}$ is micro supported on $C$.
\qed

}

\medskip
The following lemma will be used in 
the reduction of 
the proof of the existence of
singular support to the case
of projective spaces,
to produce enough sections \'etale locally.

\begin{lm}\label{lmSS'}
Let $X$ be a regular scheme of
finite type over $S$
and ${\cal F}$ be a constructible sheaf on $X$.
Let $S'$ be a $G$-torsor over $S$
for a finite group $G$.
Let $g\colon X'=X\times_SS'\to X$ be the
base change morphism
and let ${\cal F}'=g^*{\cal F}$.

{\rm 1.}
If ${\cal F}$ is micro supported on 
a closed conical subset
$C
\subset FT^*X$,
then ${\cal F}'$ is micro supported on $C'=g^\circ C
\subset FT^*X'$.

{\rm 2.}
Assume that ${\cal F}'$ is micro supported on 
a closed conical subset
$C'
\subset FT^*X'$
and that $C'$ is stable by the action of $G$.
Then, ${\cal F}$ is micro supported on $C=g_\circ C'
\subset FT^*X$.
\end{lm}

\proof{
It suffices
to show that
condition (2) in Definition \ref{dfms} is satisfied.

1. 
Let $h\colon W\to X'$ be a separated morphism
of regular schemes of finite type over $S'$.
Then $h\colon W\to X'$ is a separated morphism
of regular schemes of finite type over $S$
and the condition is satisfied.

2. 
Let $h\colon W\to X$ be a separated morphism
of regular schemes of finite type over $S$.
Then the base change 
$h'\colon W'\to X'$ is a separated morphism
of regular schemes of finite type over $S'$.
Since $C'$ is assumed stable under the $G$-action,
$h$ is $C$-transversal
if and only if
$h'$ is $C'$-transversal.
Hence the condition is satisfied.
\qed

}

\medskip 

If $k$ is a perfect field of
characteristic $p>0$ and
$X$ is a smooth scheme over $k$,
the exact sequence (\ref{eqFT}) means
a canonical isomorphism
$FT^*X\to F^*T^*X$.

\begin{pr}\label{prtensor}
Let $k$ be a perfect field of
characteristic $p>0$ and
let $X$ be a smooth scheme over $k$.
Let ${\cal F}$ be
a constructible sheaf on $X$
and
$C\subset T^*X$
be a closed conical subset.

{\rm 1.}
The following conditions are equivalent:

{\rm (1)}
${\cal F}$ is micro supported on $C
\subset T^*X$ in the sense of Beilinson
{\rm \cite{SS}}.

{\rm (2)}
${\cal F}$ is micro supported on 
$F^*C\subset FT^*X$.

{\rm 2.}
Let ${\cal G}$ be another
constructible sheaf and
assume that ${\cal F}$ and ${\cal G}$ 
are micro supported on 
closed conical subsets
$C$ and on $C'$ of $T^*X$
respectively.
If $C\cap C'$ is a subset of the
$0$-section,
then ${\cal F}\otimes{\cal G}$ is
micro supported on $C+C'
\subset T^*X$.
\end{pr}

\proof{
1.
Since $F^*T^*X=T^*X\times_XX\to T^*X$
is a homeomorphism,
the assertion follows from
\cite[Proposition 8.13]{CC}.

2.
The external tensor product
${\cal F}\boxtimes{\cal G}=
{\rm pr}_1^*{\cal F}\otimes
{\rm pr}_2^*{\cal G}$ on $X\times X$ is
micro supported on $C\times C'
=T^*X\times T^*X=T^*(X\times X)$
by \cite[Theorem 2.2.1]{Ext}.
By the assumption that
$C\cap C'$ is a subset of the
$0$-section,
the diagonal morphism
$\delta\colon X\to X\times X$
is  $C\times C'$-transversal.
Hence
${\cal F}\otimes{\cal G}
=\delta^*({\cal F}\boxtimes{\cal G})$ is
micro supported on $\delta^\circ(
C\times C')=C+C'$
by Lemma \ref{lmmstr}.
\qed

}

\subsection{Singular support}\label{ssSS}

We keep the notation that
$S$ is a regular scheme over
${\mathbf Z}_{(p)}$
satisfying the finiteness condition (F).
For a constructible sheaf ${\cal F}$ on 
a regular scheme $X$ of finite type over $S$,
we call the smallest closed conical subset
$C\subset FT^*X$ on which
${\cal F}$ is micro supported the singular support of
${\cal F}$ if it exists.

\begin{df}\label{dfss}
Let $X$ be a 
regular scheme of
finite type over $S$ and
let ${\cal F}$ be a constructible
sheaf on $X$.
We say that 
a closed conical subset 
$C\subset FT^*X$
is the singular support $SS{\cal F}$ of ${\cal F}$
if for any closed conical subset 
$C'\subset FT^*X$,
the inclusion $C\subset C'$
is equivalent to the condition
that ${\cal F}$ is micro supported on $C'$.
\end{df}

We don't know the existence of singular support
in general because the condition
that ${\cal F}$ is micro supported on $C$ and on
$C'$ does not a priori imply
that ${\cal F}$ is micro supported on 
the intersection $C\cap C'$.
Similarly as in the proof of
\cite[Theorem 1.3 (i)]{SS},
the proof of the existence of
singular support is reduced
to the case where $X={\mathbf P}^n_S$
as Corollary \ref{corSS} below.

\begin{lm}\label{lmSSloc}
Let $X$ be a 
regular schemes of
finite type over $S$ and
let ${\cal F}$ be a constructible
sheaf on $X$.

{\rm 1.}
Let $U\subset X$ be an open subscheme.
Assume that $C\subset FT^*X$ is the singular support
of ${\cal F}$.
Then, $C|_U$ is the singular support
of ${\cal F}|_U$.

{\rm 2.}
Let $(U_i)$ be an open covering of
$X$ and $C_i
\subset FT^*U_i$ be the singular support
of ${\cal F}|_{U_i}$.
Then,
$C=\bigcup_iC_i
\subset FT^*X$ is the singular support
of ${\cal F}$.
\end{lm}

\proof{
1.
The restriction ${\cal F}|_U$
is micro supported on
$C|_U$ by Lemma \ref{lmmsloc}.1.
We show that $C|_U$ is the smallest.
Let $C'\subset FT^*U$ be a closed conical
subset on which ${\cal F}|_U$
is micro supported
and $Z=X\sm U$ be the complement.
Then ${\cal F}$ is micro supported
on $\overline {C'}\cup C|_Z$
by Lemma \ref{lmmsloc}.2.
Since $C$ is the smallest,
we have
$C\subset 
\overline {C'}\cup C|_Z$
and
$C|_U\subset 
(\overline {C'}\cup C|_Z)|_U=C'$.

2.
For every $i,j$,
the restrictions $C_i|_{U_i\cap U_j}$ and
$C_j|_{U_i\cap U_j}$
are the singular support of 
${\cal F}|_{U_i\cap U_j}$
and are the same by 1.
Hence the union 
$C=\bigcup_iC_i$ is a closed conical subset of $FT^*X$.
Since $C_i=C|_{U_i}$,
${\cal F}$ is micro supported on $C$
by Lemma \ref{lmmsloc}.1.

We show that $C$ is the smallest.
Let $C'\subset
FT^*X$ be a closed conical subset 
on which
${\cal F}$ is micro supported.
Then, for each $i$,
we have $C_i\subset C'|_{U_i}$.
Hence we have $C\subset C'$.
\qed
}

\begin{pr}\label{prSSXP}
Let $i\colon X\to P$ be a closed immersion
of regular schemes of
finite type over $S$ and
let ${\cal F}$ be a constructible
sheaf on $X$.
Let $C_P\subset FT^*P$ 
be a closed conical subset
and
assume that $C_P$ is the singular support
of $i_*{\cal F}$.
Then the following holds.

{\rm 1.}
$C_P$ is a subset of
$FT^*P|_X$.

{\rm 2.}
Define $C\subset FT^*X$ 
to be the closure of
the image of $C_P$
by the surjection
$FT^*P|_X\to FT^*X$.
Then $C$ is the singular support $SS{\cal F}$.
\end{pr}

\proof{
1.
Let $U=P\sm X$ be the complement.
Since ${\cal F}|_U=0$
is micro supported on $\varnothing$
by Lemma \ref{lmmslc}.3,
${\cal F}$ is micro supported on
$C_P|_X\subset FT^*P$ by Lemma \ref{lmmsloc}.2.
Since $C_P$ is the smallest,
we have $C_P=C_P|_X
\subset FT^*P|_X$.

2.
By Lemma \ref{lmSSXP},
${\cal F}$ is micro supported on $C$.
We show that $C$ is the smallest.
Assume that ${\cal F}$ is micro supported on 
a closed conical subset
$C'\subset FT^*X$.
Then, since $i_*{\cal F}$ is micro supported on $i_\circ C'$
by  Lemma \ref{lmSSXP}.1,
we have $C_P\subset i_\circ C'$.
By taking the closure of the image 
by the surjection
$FT^*P|_X\to FT^*X$,
we obtain $C\subset C'$.
\qed

}

\begin{cor}\label{corSS}
Assume that the singular support
exists for every constructible sheaf
on every projective space
${\mathbf P}$ 
over $S$.
Then, 
the singular support
$SS{\cal F}$ exists for every constructible sheaf
${\cal F}$ on every regular scheme
$X$ of finite type over $S$.
\end{cor}

\proof{
By Lemma \ref{lmSSloc},
the assertion is local on $X$.
Hence, we may assume that $X$ is affine
and is a closed subscheme of
${\mathbf A}^n_S$.
By the assumption and Lemma \ref{lmSSloc},
the assertion holds for
${\mathbf A}^n_S$.
Hence the assertion follows from Proposition 
\ref{prSSXP}.
\qed

}

%
%
%

\subsection{$S$-acyclicity}

We keep the notation that
$S$ is a regular noetherian scheme over
${\mathbf Z}_{(p)}$
satisfying the finiteness condition (F).
We will introduce two notions
the $S$-acyclicity 
(Definition \ref{dfSacyc} and \ref{dfFf0})
and 
the smoothness modulo $S$ 
(Definition \ref{dff0}) on
a scheme $Y$ smooth over $S$
as a preparation of the definitions 
(Definitions \ref{dfCacyc}
and \ref{dfFacyc})
and a property (Lemma \ref{lmuniCacy}.3) of
acyclicity over $S$
in the next chapter.
If $Y=Y_0\times_kS$ is a base change
of a smooth scheme $Y_0$ over a perfect
field $k$ of characteristic $p>0$,
the two notions should be conditions on $Y_0$.
Because we don't have such a situation
in mixed characteristic,
we formulate the conditions using the injection
$FT^*S\times_SY\to FT^*Y$.

\begin{df}\label{dfSacyc}
Let $Y\to S$ be a smooth scheme
over $S$.
We say that a closed conical subset
$C'\subset FT^*Y$ is $S$-acyclic
if the inverse image of
$C'$ by the injection
$FT^*S\times_SY\to FT^*Y$
is a subset of the $0$-section.
\end{df}

We may define the $S$-acyclicity without
assuming that $Y$ is smooth over $S$.
However, since the condition implies
that $Y$ is smooth on a neighborhood
of the base of $C'$ if $S$ is assumed excellent,
we impose the smoothness.

\begin{lm}\label{lmCSacyc}
Let $g\colon Y\to S$ 
be a smooth scheme over $S$
and
let $C'\subset FT^*Y$
be an $S$-acyclic closed conical subset.
Let $X$ be a regular scheme of
finite type over $S$.

{\rm 1.}
The second projection
${\rm pr}_2\colon X\times_SY\to Y$
is $C'$-transversal and
${\rm pr}_2^\circ C'
\subset FT^*(X\times_SY)$
is $X$-acyclic.

{\rm 2.}
Let $C\subset FT^*X$ be a closed
conical subset. Then,
the intersection ${\rm pr}_1^\circ C\cap 
{\rm pr}_2^\circ C'\subset FT^*(X\times_SY)$
is a subset of the $0$-section.

{\rm 3.}
Let $C'_1\subset FT^*Y$
be a closed conical subset
satisfying $C'+
(FT^*S\times_SY)\subset C'_1$.
Let $f\colon X\to Y$ be a morphism over
$S$ and assume that $f$ is
$C'_1$-transversal.
Then, 
on a neighborhood of the inverse image 
$f^{-1}(B')$ of the base $B'$ of $C'$,
the scheme $X$ is smooth over $S$
and $f^\circ C'$ is $S$-acyclic.

Further if $S'\to S$
is a morphism of regular schemes of
finite type,
the base change $X'=X\times_SS'
\to Y'=Y\times_SS'$
is ${\rm pr}_1^\circ C'$-transversal
on a neighborhood of the inverse image of
$f^{-1}(B')$.
\end{lm}

\proof{1.
Since ${\rm Ker}(FT^*Y\times_Y(X\times_SY)
\to FT^*(X\times_SY))$
is a subset of
\begin{align*}
&{\rm Ker}(FT^*Y\times_Y(X\times_SY)
\to F^*T^*(X\times_SY)/X|_{(X\times_SY)
_{{\mathbf F}_p}})\\
&=
{\rm Ker}(FT^*Y
\to F^*T^*Y/S
|_{Y_{{\mathbf F}_p}})\times_Y(X\times_SY),
\end{align*}
the projection
${\rm pr}_2\colon X\times_SY\to Y$
is $C'$-transversal.

Since
$F^*T^*Y/S|_{Y_{{\mathbf F}_p}}
\times_Y(X\times_SY)
\to  F^*T^*(X\times_SY)/X|_{(X\times_SY)
_{{\mathbf F}_p}}$
is an isomorphism,
the intersection
${\rm pr}_2^\circ C'
\cap {\rm Im}
(FT^*X\times_X(X\times_SY)\to
FT^*(X\times_SY))$
equals the image of
$\bigl(C'
\cap {\rm Im}
(FT^*S\times_SY\to
FT^*Y)\bigr)
\times_Y(X\times_SY)$.
Hence
the $S$-acyclicity of $C'$
implies
the $X$-acyclicity of
${\rm pr}_2^\circ C'$.

2.
The intersection ${\rm pr}_1^\circ C\cap 
{\rm pr}_2^\circ C'
\subset
{\rm pr}_2^\circ C'\cap {\rm Im}(FT^*X\times_X
(X\times_SY)\to FT^*(X\times_SY))$
is a subset of the $0$-section by 1.

3.
Since $f$ is
$C'_1$-transversal
and $(FT^*S\times_SY)|_{B'}\subset C'_1$,
the morphism
$FT^*S\times_SX\to FT^*X$
is an injection on the inverse image
$f^{-1}(B')$.
Hence $X$ is smooth over $S$
on a neighborhood of $f^{-1}(B')$
by \cite[Proposition 2.10]{FW}.

To show the rest, we may assume $X$ is smooth.
Since $f$ is $C'_1$-transversal,
the intersection 
$f^\circ C'\cap {\rm Im}(FT^*S\times_SX
\to FT^*X)$
is a subset of the image of
$f^*C'\cap {\rm Im}(FT^*S\times_SX
\to FT^*Y\times_YX)$.
This is a subset of the $0$-section
since $C'$ is assumed $S$-acyclic.

By 1,
${\rm pr}_1\colon Y'\to Y$ is $C'$-transversal
and ${\rm pr}_1^\circ C'$ is defined.
Since $f^\circ C'$ is $S$-acyclic,
${\rm pr}_1\colon X'\to X$ is $f^\circ C'$-transversal.
Since ${\rm pr}_1f'=f\circ {\rm pr}_1$,
the base change $f'\colon X'\to Y'$
is ${\rm pr}_1^\circ C'$-transversal by
Lemma \ref{lmCtrtr}.
\qed

}

%
%
%
%
%

\begin{df}\label{dfFf0}
Let $Y$ be a smooth
scheme of finite type over $S$.
We say that a constructible sheaf ${\cal G}$ on $Y$
is $S$-acyclic 
if there exists an $S$-acyclic 
closed conical subset
$C'\subset FT^*Y$ on which
${\cal G}$ is micro supported.
\end{df}

If $S={\rm Spec}\, k$ for a perfect field
of characteristic $p$,
then every $C'$ is $S$-acyclic
and hence every ${\cal G}$ is $S$-acyclic.
If $Y=S$ is excellent, the condition that
${\cal F}$ is $S$-acyclic means
that ${\cal F}$ is locally constant
on a neighborhood of $S_{{\mathbf F}_p}$
by Lemma \ref{lmmslc}.

\begin{lm}\label{lmSacyc}
Assume that $S$ is excellent.
Let $Y$ be a smooth 
scheme of finite type over $S$
and let
${\cal G}$ be an $S$-acyclic constructible sheaf
on $Y$.
Then, $g\colon Y\to S$
is universally ${\cal G}$-acyclic
along $Y_{{\mathbf F}_p}$.
\end{lm}

\proof{
Let $C'\subset FT^*Y$
be an $S$-acyclic closed conical subset
on which ${\cal G}$ is
micro supported.
First, we show that
it suffices to show that $g$ is ${\cal G}$-acyclic at 
$Y_{{\mathbf F}_p}$.
By Lemma \ref{lmCSacyc}.1,
for every morphism $S'\to S$ 
of finite type
of regular noetherian schemes,
the base change $f\colon Y'=Y\times_SS'\to Y$
is $C'$-transversal
and $f^\circ C'$ is $S'$-acyclic.
Hence the pull-back
${\cal G}'=f^*{\cal G}$ on $Y'$
is micro supported on $f^\circ C'$
by Lemma \ref{lmmstr}.
Thus by Lemma \ref{lmula},
it suffices to show that $g$ is ${\cal G}$-acyclic 
along $Y_{{\mathbf F}_p}$.

Let
$$\begin{CD}
Y@<f<< Y'@<h<< W\\
@VgVV@VVV@VVV\\
S@<<< S'@<<< V
\end{CD}$$
be a cartesian diagram
where $V\to S'$ is an immersion
of regular schemes of finite type over $S$
as in (\ref{eqFtracyc}) in
Proposition \ref{prFtracyc}.
By Lemma \ref{lmCSacyc}.1,
$f$ and $fh$ are $C'$-transversal
and 
$f^\circ C'$ is $S'$-acyclic.
Hence $h$ is
$f^\circ C'$-transversal
by Lemma \ref{lmCtrtr}.
Hence the pull-back
${\cal G}'=f^*{\cal G}$ on $Y'$
is micro supported on $f^\circ C'$
by Lemma \ref{lmmstr}
and $h$ is ${\cal G}'$-transversal
on a neighborhood of $W_{{\mathbf F}_p}$.
Thus by Proposition \ref{prFtracyc},
$f$ is ${\cal G}$-acyclic along $Y_{{\mathbf F}_p}$.
\qed

}

%
%

%

\begin{lm}\label{lmSacyfun}
{\rm 1.}
Let $i\colon Z\to Y$ be a closed immersion
of smooth schemes over $S$
and let ${\cal G}$
be an $S$-acyclic constructible sheaf on $Z$.
Then $i_*{\cal G}$ on $Y$
is $S$-acyclic.

{\rm 2.}
Let $f\colon Y\to Y'$ be a morphism
of smooth schemes over $S$
and ${\cal G}'$ be an $S$-acyclic
constructible sheaf on $Y'$
micro supported on an $S$-acyclic
closed conical subset 
$C'\subset FT^*Y'$.
Let $C'_1\subset FT^*Y$
be a closed conical subset
satisfying $C'+
(FT^*S\times_SY)\subset C'_1$
and assume that $f$ is
$C'_1$-transversal.
Then the pull-back
${\cal G}=f^*{\cal G}'$ on $Y$ is $S$-acyclic.

{\rm 3.}
Let $S'\to S$ be a separated morphism
of finite type of regular noetherian schemes.
Let $Y\to S$ be a smooth morphism
and $h\colon Y'=Y\times_SS'\to Y$ be the base change.
If a constructible sheaf 
${\cal G}$ on $Y$ is $S$-acyclic,
then $h$ is ${\cal G}$-transversal
on a neighborhood of $Y'_{{\mathbf F}_p}$
and the pull-back
${\cal G}'=h^*{\cal G}$ on $Y'$ is $S'$-acyclic.
\end{lm}

\proof{
1.
Let $C'\subset FT^*Z$ be an $S$-acyclic 
closed conical subset on which
${\cal G}$ is micro supported.
Then $i_*{\cal G}$ is micro supported
on $i_\circ C'\subset FT^*Y$ 
by Lemma \ref{lmg0F}
and $i_\circ C'$ is $S$-acyclic.

2.
The pull-back
${\cal G}$ on $Y$
of ${\cal G}'$ is micro supported
on $C=g^\circ C'\subset FT^*Y$.
Since $C$ is $S$-acyclic
by Lemma \ref{lmCSacyc}.3,
${\cal G}$ is $S$-acyclic.

3.
Assume that ${\cal G}$ is 
micro supported on $S$-acyclic $C'$.
By Lemma \ref{lmCSacyc}.1,
$h\colon
Y'\to Y$ 
is $C'$-transversal
and $h^\circ C'$ is $S'$-acyclic.
Hence $h$
is ${\cal G}$-transversal
and ${\cal G}'=h^*{\cal G}$ is micro supported
on $h^\circ C'$
by Lemma \ref{lmmstr}.
Hence ${\cal G}'$ is $S'$-acyclic.
\qed

}

\begin{df}\label{dff0}
Let $f\colon W\to Y$ be a morphism
of regular scheme of finite type
over $S$.
Assume that $Y$ is smooth over $S$.
We say that $f$ is smooth modulo $S$
if we have an inclusion
\begin{equation}
{\rm Ker}(FT^*Y\times_YW \to FT^*W)
\subset
{\rm Im}(FT^*S\times_SW \to FT^*Y\times_YW).
\label{eqSsm}
\end{equation}

For a closed subset $Z\subset W$,
we say that $f$ is smooth modulo $S$ on $Z$
if we have an inclusion {\rm (\ref{eqSsm})}
on the inverse image of $Z$.
\end{df}

If $S$ is smooth over
a perfect field $k$
of characteristic $p$
and if $Y=Y_0\times_kS$
for a scheme $Y_0$ smooth over $k$,
then the smoothness 
of $W\to Y$ modulo $S$
implies the smoothness
of the composition $W\to Y_0$.

\begin{lm}\label{lmf0}
Let $f\colon W\to Y$
be a morphism of regular schemes
of finite type over $S$
and assume that $Y$ is smooth over $S$.

{\rm 1.}
If $Y=S$,
then every morphism $W\to S$
is smooth modulo $S$.

{\rm 2.}
Let $g\colon Y'\to Y$ be a morphism of
smooth schemes over $S$
and set $W'=W\times_YY'$.
If $f$ is smooth modulo $S$
and if $S$ is excellent,
then there exists a neighborhood $U'\subset W'$ of 
$W'_{{\mathbf F}_p}$ 
such that
$g$ is transversal to $f$
on $U'$  and
the base change
$f'\colon W'\to Y'$
of $f$ is smooth modulo $S$ on $U'$.
In particular,
the projection 
${\rm pr}_2\colon W\times_SY\to Y$
is smooth modulo $S$.

{\rm 3.}
Let $S'\to S$ be a smooth morphism
and let 
$f'\colon W\to Y'=Y\times_SS'$
be a morphism over $S'$.
If $f'$ is smooth modulo $S'$, 
then $f={\rm pr}_1\circ f'$ is smooth modulo $S$.

{\rm 4.}
Assume that $S$ is smooth over 
a perfect field $k$ of characteristic $p$
and $Y=Y_0\times_kS$
for a smooth scheme $Y_0$ over $k$.
Then $f$ is smooth modulo $S$
if and only if
the composition ${\rm pr}_1\circ f\colon
W\to Y_0$ is smooth.
In particular if $S={\rm Spec}\, k$,
then $f$ is smooth modulo $k$
means
that $f$ is smooth.
\end{lm}

\proof{1.
Clear from (\ref{eqSsm}).

2.
By replacing $W\to Y$ by 
$W\times_SY'\to Y\times_SY'$,
we may assume that $Y'\to Y$ is a
regular immersion of smooth schemes
over $S$.
Then, the morphism
$T^*_{Y'}Y\to T^*Y/S$
is an injection.
Hence the intersection
$F^*T^*_{Y'}Y
\cap {\rm Im}(FT^*S\times_SY\to FT^*Y)$
is a subset of the $0$-section.
Thus the morphism
$F^*T^*_{Y'}Y\times_YW
\to FT^*W$
is also an injection and
$W'$ is regular on a neighborhood of 
$W'_{{\mathbf F}_p}$
by \cite[Theorem 3.4]{FW}.
Further 
${\rm Ker}(FT^*Y\times_YW \to FT^*W)
\to 
{\rm Ker}(FT^*Y'\times_{Y'}W' \to FT^*W')$
is an isomorphism and the assertion follows.

3.
Since the inverse image of
${\rm Im}(FT^*S'\times_{S'}W \to FT^*Y'\times_{Y'}W)$
by 
$FT^*Y\times_YW \to FT^*Y'\times_{Y'}W$
equals
${\rm Im}(FT^*S\times_SW \to FT^*Y\times_YW)$,
the assertion follows.

4.
The assertion follows from the exact
sequence
$0\to FT^*S\times_SW \to FT^*Y\times_YW
\to FT^*Y_0\times_{Y_0}W\to 0$.
\qed

}

%

\begin{lm}\label{lmsmacyc}
Let $f\colon W\to Y$ be a morphism
of regular schemes of finite type over $S$.
Assume that $Y$ is smooth over $S$
and that $f$ is smooth modulo $S$.

{\rm 1.}
Let $C'\subset FT^*Y$
be a closed conical subset.
If $C'$ is $S$-acyclic,
then 
$f$ is $C'$-transversal.

{\rm 2.}
Let
${\cal G}$ be an $S$-acyclic constructible sheaf
on $Y$.
Then $f$ is ${\cal G}$-transversal
on a neighborhood of
$W_{{\mathbf F}_p}$.
\end{lm}

\proof{
1.
Clear from the inclusion
(\ref{eqSsm}).

2.
Clear from 1.
\qed

}

\begin{lm}\label{lmcn}
Let $S$ be a smooth scheme
over a perfect field $k$.
Let $Y$ be a smooth scheme of finite type over $S$
and let
${\cal G}$ be a sheaf
on $Y$
micro supported
on an $S$-acyclic 
$C'\subset FT^*Y$.
Let $X$ be a regular scheme
of finite type over $S$
and 
let ${\cal F}$ be a sheaf on $X$
micro supported on $C
\subset FT^*X$.
Then,  
${\cal F}\boxtimes {\cal G}$
on $X\times_SY$ is
micro supported on ${\rm pr}_1^\circ C
+{\rm pr}_2^\circ C'$.
\end{lm}

\proof{
Since the first projection
${\rm pr}_1\colon X\times_SY\to X$
is smooth,
${\rm pr}_1$ is $C$-transversal
and 
${\rm pr}_1^*{\cal F}$ is 
micro supported on ${\rm pr}_1^\circ C$
by Lemma \ref{lmmstr}.
Since the second projection
${\rm pr}_2\colon X\times_SY\to Y$
is smooth modulo $S$ by Lemma \ref{lmf0}.2,
${\rm pr}_2$ is $C'$-transversal
by Lemma \ref{lmsmacyc}.1
and 
${\rm pr}_2^*{\cal G}$ is 
micro supported on ${\rm pr}_2^\circ C'$
by Lemma \ref{lmmstr}.
The intersection
${\rm pr}_1^\circ C\cap 
{\rm pr}_2^\circ C'$ is a subset of the $0$-section
by Lemma \ref{lmCSacyc}.2.
Hence by Proposition \ref{prtensor},
${\cal F}
\boxtimes{\cal G}$
is micro supported on
${\rm pr}_1^\circ C+ 
{\rm pr}_2^\circ C'$.
\qed

}

\medskip
We prove a weaker assertion in general.

\begin{pr}
Let $g\colon Y\to S$ be a smooth scheme over $S$
and let ${\cal G}$ be an $S$-acyclic constructible
sheaf on $Y$ micro supported on an
$S$-acyclic closed conical subset $C'\subset FT^*Y$.
Assume that $S$ is excellent.
Let $C'_1\subset FT^*Y$
be a closed conical subset
satisfying $C'+
(FT^*S\times_SY)\subset C'_1$.
Then,
for every constructible sheaf ${\cal H}$
on $S$,
the tensor product
${\cal G}\otimes g^*{\cal H}$ is
micro supported on $C'_1$.
\end{pr}

\proof{
Let 
\begin{equation}
\xymatrix{
W\ar[r]^f\ar[rd]_q&Y\ar[d]^g\\
&S}
\label{eqWYS}
\end{equation}
 be a
commutative diagram of
separated morphism of regular schemes of finite
type over $S$
and assume that 
$f$ is $C'_1$-transversal.
We show that $f$ is
${\cal G}\otimes g^*{\cal H}$-transversal
on a neighborhood of $W_{{\mathbf F}_p}$.
Since the base $B'$ of $C'$ contains 
the intersection of the support of ${\cal G}$
with $Y_{{\mathbf F}_p}$,
we may assume that $W$ is smooth over $S$
by Lemma \ref{lmCSacyc}.3.

Since $f\colon W\to Y$ is $C'$-transversal,
it is ${\cal G}$-transversal 
on a neighborhood of $W_{{\mathbf F}_p}$.
Hence by Lemma \ref{lmcGL},
it suffices to show that the morphism
$c_{f,{\cal G},g^*{\cal H}}
\colon f^!{\cal G}\otimes q^*{\cal H}\to
f^!({\cal G}\otimes g^*{\cal H})$
is an isomorphism
on a neighborhood of $W_{{\mathbf F}_p}$.
By devissage,
we may assume that there
exists an immersion $j\colon V\to S$
of regular schemes,
a locally constant constructible sheaf
${\cal L}$ on $V$ and
and 
${\cal H}=j_*{\cal L}$.

Let $$\xymatrix{
W'\ar[r]^{f'}\ar[rd]_{q'}&Y'\ar[d]^{g'}\\
&V}
$$ be the base change of
(\ref{eqWYS}) by $j\colon V\to S$
and let $j_Y\colon Y'\to Y$ and $j_W\colon W'\to W$
be the induced morphisms.
Since $g$ is ${\cal G}$-acyclic
by Lemma \ref{lmSacyc}
along $Y_{{\mathbf F}_p}$,
the canonical morphism
${\cal G}\otimes g^*j_*{\cal L}
\to j_{Y*}(j_Y^*{\cal G}\otimes g'^*{\cal L})$
is an isomorphism 
on a neighborhood of $Y_{{\mathbf F}_p}$
by Proposition \ref{prla0}.
Similarly, since 
$f^*{\cal G}$ is also $S$-acyclic by Lemma \ref{lmSacyfun}.2,
the morphism
$f^!{\cal G}\otimes q^*j_*{\cal L}
\to j_{W*}(j_W^*f^!{\cal G}\otimes q'^*{\cal L})$
is an isomorphism
on a neighborhood of $W_{{\mathbf F}_p}$.
Thus, it suffices to show that
the morphism
$f^! j_{Y*}(j_Y^*{\cal G}\otimes g'^*{\cal L})
\to  j_{W*}(j_W^*f^!{\cal G}\otimes q'^*{\cal L})$
is an isomorphism
on a neighborhood of $W_{{\mathbf F}_p}$.

Since the morphism
$f^! j_{Y*}\to j_{W*}f'^!$ (\ref{eqbcdef2})
is an isomorphism,
it suffices to show that
$f'^! (j_Y^*{\cal G}\otimes g'^*{\cal L})
\to  j_W^*f^!{\cal G}\otimes q'^*{\cal L}$
is an isomorphism
on a neighborhood of $W_{{\mathbf F}_p}$.
Since ${\cal L}$ is locally constant,
it is reduced to showing that
$f'^! j_Y^*{\cal G}
\to  j_W^*f^!{\cal G}$
is an isomorphism
on a neighborhood of $W'_{{\mathbf F}_p}$.
By Lemma \ref{lmSacyfun}.3 and Lemma \ref{lmCSacyc}.1,
on a neighborhood of $Y'_{{\mathbf F}_p}$,
$j_Y$ is ${\cal G}$-transversal 
and $j_Y^*{\cal G}$ is $V$-acyclic
and is micro supported on $j_Y^\circ C'$.
By Lemma \ref{lmCSacyc}.3,
$f'$ is
$j_Y^\circ C'$-transversal and
is $j_Y^*{\cal G}$-transversal
on a neighborhood of $W'_{{\mathbf F}_p}$.
Since $f$ is ${\cal G}$-transversal
on a neighborhood of $W_{{\mathbf F}_p}$
and is transversal to $j_Y$,
the morphism
$f'^! j_Y^*{\cal G}
\to  j_W^*f^!{\cal G}$
is an isomorphism
on a neighborhood of $W'_{{\mathbf F}_p}$
by Lemma \ref{lmbc}.1.
\qed

}

\section{$S$-micro support}\label{sSms}

Let $S$ be an excellent regular noetherian scheme
over ${\mathbf Z}_{(p)}$
satisfying the finiteness condition:

\noindent
(F)
The reduced part of
$S_{{\mathbf F}_p}$
is of finite type over a field of characteristic $p$
with finite $p$-basis.

\subsection{$C$-acyclicity over $S$}\label{ssCacycS}

\begin{df}\label{dfCacyc}
Let $h\colon W\to X$ and
$f\colon W\to Y$ be morphisms
of regular schemes of finite type over $S$.
Assume that $Y$ is smooth over $S$.
Let $C\subset FT^*X$ be a closed conical subset.

{\rm 1.}
Let $w\in W_{{\mathbf F}_p}$.
We say that $(h,f)$ is $C$-acyclic
over $S$
at $w$ if we have an inclusion
\begin{align}
(h^*C\times_W&\,  (FT^*Y\times_YW))
\cap {\rm Ker}((FT^*X\times_XW)\times_W(FT^*Y\times_YW)
\to FT^*W)
\label{eqCacyc}
\\
&
\subset
{\rm Ker}((FT^*X\times_XW)\times_W(FT^*Y\times_YW)\to
FT^*(X\times_SY)\times_{X\times_SY}W).
\nonumber
\end{align}
on the fiber at $w$.

We say that $(h,f)$ is $C$-acyclic
over $S$
if $(h,f)$ is $C$-acyclic
over $S$
at every point of $W_{{\mathbf F}_p}$.

{\rm 2.}
We say $(h,f)$ is  universally $C$-acyclic
over $S$
if for every morphism
$g\colon Y'\to Y$ of smooth schemes over $S$
and the commutative diagram
\begin{equation}
\xymatrix{
&W'\ar[r]^{f'}\ar[ld]_{h'}\ar[d]^{g'}&Y'\ar[d]^g\\
X&W\ar[l]_h\ar[r]^f &Y}
\label{eqCuniv}
\end{equation}
with cartesian square,
there exists a neighborhood $U'\subset W'$
of the inverse image
$h'^{-1}(B)$ of
the base $B$ of $C$
such that $g$ is transversal to $f$
on $U'$
and the pair $(h',f')$ is $C$-acyclic
on $U'$.
\end{df}

Unlike the $C$-transversality,
the $C$-acyclicity over $S$ may not
be an open condition in general.
Since
$(FT^*Y\times_YW))
\cap {\rm Ker}((FT^*X\times_XW)\times_W(FT^*Y\times_YW)
\to FT^*W$
factors through
$FT^*(X\times_SY)\times_{
X\times_SY}W$,
the inclusion
(\ref{eqCacyc})
is equivalent to the inclusion
$${\rm Im}(h^*C\times_W
(FT^*Y\times_YW))
\cap 
{\rm Ker}
(FT^*(X\times_SY)\times_{
X\times_SY}W
\to FT^*W)
\subset {\text {0-section}}
\leqno{\rm (\ref{eqCacyc})'}$$
\noindent
\!\!\!\!
in $FT^*(X\times_SY)\times_{
X\times_SY}W$.
Using this interpretation,
we will show in Lemma \ref{lmCCbar} that  
the $C$-acyclicity over $S$ 
is an open condition for $S$-saturated $C$
defined in Definition \ref{dfsat}.1.

If $S={\rm Spec}\, k$ for a perfect field $k$
of characteristic $p$,
then the $C$-acyclicity is the same
as the $C$-transversality for a test pair
defined by Beilinson \cite[1.2]{SS}
by Lemma \ref{lmuniCacy}.4,
Lemma \ref{lmCacycq}
and Lemma \ref{lmCacychs}.
We show in Lemma \ref{lmuniCacy}.1 below
that the $C$-acyclicity implies 
the universal $C$-acyclicity. 
If $C'\subset FT^*X$
is another closed conical subset
satisfying 
$C\subset C'\subset
C+FT^*S\times_SX,$
then the $C$-acyclicity over $S$
and 
the $C'$-acyclicity over $S$
are equivalent to each other
since the kernels in 
(\ref{eqCacyc}) are stable
under the actions of
$FT^*S\times_SW$.



\begin{lm}\label{lmuniCacy}
Let $h\colon W\to X$ and
$f\colon W\to Y$ be morphisms
of regular scheme of finite type
over $S$.
Let $C\subset FT^*X$ be a closed conical
subset and let $B\subset X$ be
the base of $C$.
Assume that $Y$ is smooth over $S$
and that 
$(h,f)$ is $C$-acyclic over $S$.

{\rm 1.}
The pair $(h,f)$ is universally $C$-acyclic over $S$.

{\rm 2.}
On the inverse image $h^{-1}(B)\subset W$,
the kernel
${\rm Ker}(FT^*Y\times_YW\to FT^*W)$
is a subset of the image of
${\rm Ker}(FT^*S\times_SW\to FT^*X\times_XW)$.

{\rm 3.}
The morphism $f\colon
W\to Y$ is smooth modulo $S$
on $h^{-1}(B)$.

If $X$ is smooth over $S$,
then
$W\to Y$ is smooth
on a neighborhood of $h^{-1}(B)$.

{\rm 4.}
The morphism
$h$ is $C$-transversal.
If $q\colon W\to S$ denote
the canonical morphism,
then $(h,q)$ is $C$-acyclic over $S$.

{\rm 5.}
Let $C'\subset FT^*Y$
be a closed conical subset
and assume that $C'$ is $S$-acyclic.
Then the morphism
$(h,f)\colon W\to X\times_SY$
is ${\rm pr}_1^\circ C
+{\rm pr}_2^\circ C'$-transversal.
\end{lm}

\proof{
1.
Let $g\colon Y'\to Y$
be a morphism of smooth schemes
over $S$ and
consider the commutative diagram (\ref{eqCuniv}).
By Lemma \ref{lmuniCacy}.3 and Lemma \ref{lmf0}.2,
the morphism $g$ is transversal to $f$
on a neighborhood $U'\subset W'$ of $h'^{-1}(B)$.
To show that $(h',f')$ is $C$-acyclic over $S$,
it suffices to show
that 
\begin{align*}
&{\rm Ker}((FT^*X\times_XW)
\times_W
(FT^*Y\times_YW)
\to FT^*W)\times_WU'\\
&\to 
{\rm Ker}((FT^*X\times_XU')
\times_{U'}
(FT^*Y'\times_{Y'}U')
\to FT^*U')
\end{align*}
is an isomorphism.
By the factorization
$Y'\to Y\times_SY'\to Y$,
it suffices to prove the cases
where $Y'\to Y$ is smooth
and $Y'\to Y$ is an immersion
respectively.
Since the morphism $f\colon W\to Y$
is transversal to $g\colon Y'\to Y$ 
in the second case, the assertion follows.

2.
The intersection
$$((FT^*_XX\times_XW)\times_W
(FT^*Y\times_YW))
\cap
{\rm Ker}((FT^*X\times_XW)\times_W
(FT^*Y\times_YW)\to
FT^*W)$$
is identified
with
${\rm Ker}(FT^*Y\times_YW\to
FT^*W)$.
In the commutative diagram 
$$\begin{CD}
0@>>>
FT^*S\times_SW
@>>>
FT^*Y\times_YW
@>>>
F^*T^*Y/S\times_YW
@>>>0\\
@.@VVV@VVV@VVV@.\\
0@>>>
FT^*X\times_XW
@>>>
FT^*(X\times_SY)\times_{X\times_SY}W
@>>>
(F^*T^*(X\times_SY)/X)\times
_{X\times_SY}W@>>>0
\end{CD}$$
of exact sequences,
the right vertical arrow is an isomorphism.
Hence the intersection
\begin{align*}
&((FT^*_XX\times_XW)\times_W
(FT^*Y\times_YW))
\\&\
\cap
{\rm Ker}((FT^*X\times_XW)\times_W
(FT^*Y\times_YW)\to
FT^*(X\times_SY)_{X\times_SY}W)
\end{align*}
is identified
with the image of
${\rm Ker}(FT^*S\times_SW\to
FT^*X\times_XW)$.
Thus, 
we obtain the inclusion.

3.
By 2, we have an inclusion (\ref{eqSsm})
on the inverse image $h^{-1}(B)$.

If $X$ is smooth over $S$,
then $FT^*S\times_SX\to FT^*X$
is an injection. Hence by 2,
$FT^*Y\times_YW\to FT^*W$
is also an injection on  $h^{-1}(B)$.

4.
Assume that we have an inclusion 
(\ref{eqCacyc}).
The intersection of (\ref{eqCacyc})
with the first factor
$FT^*X\times_XW$ gives an inclusion
$$h^*C\cap {\rm Ker}(FT^*X\times_XW
\to FT^*W)
\subset {\rm Ker}(FT^*X\times_XW
\to FT^*(X\times_SY)_{X\times_SY}W).$$
Since $Y$ is smooth over $S$,
the last morphism is an injection
and hence $h$ is $C$-transversal.

Since 
the morphism
$FT^*S\times_SY
\to
FT^*Y$ is an injection and
the morphism
${\rm Ker}
((FT^*X\times_XW)\times_W
(FT^*S\times_SW)
\to
FT^*X\times_XW)
\to
{\rm Ker}
((FT^*X\times_XW)\times_W
(FT^*Y\times_YW)
\to
FT^*(X\times_SY)\times_{X\times_SY}W)
$
is an isomorphism,
the $C$-transversality of $(h,f)$ over $S$
implies that of $(h,q)$.

5.
The assumption that $C'$ is $S$-acyclic 
implies that that the
intersection
$(C\times C')
\cap {\rm Ker}(FT^*X\times_SFT^*Y
\to FT^*(X\times_SY))$ is a subset
of the $0$-section
and that
${\rm pr}_1^\circ C+{\rm pr}_2^\circ C'
\subset T^*(X\times_SY)$
is the image of $C\times C'$
and is a closed subset.
If $(h,f)$ is $C$-acyclic,
then the intersection
$(h^*C\times_Wf^*C')
\cap {\rm Ker}((FT^*X\times_XW)\times_W(FT^*Y\times_YW)
\to FT^*W)$
is a subset of
${\rm Ker}((FT^*X\times_XW)\times_W(FT^*Y\times_YW)\to
FT^*(X\times_SY)\times_{X\times_SY}W)$.
Since ${\rm pr}_1^\circ C
+{\rm pr}_2^\circ C'$ is the image of
$h^*C\times_Wf^*C'$
by the surjection
$(FT^*X\times_XW)\times_W(FT^*Y\times_YW)\to
FT^*(X\times_SY)\times_{X\times_SY}W$,
its intersection with the kernel of
$FT^*(X\times_SY)\times_{X\times_SY}W
\to FT^*W$
is a subset of the $0$-section.
\qed

}
\begin{lm}\label{lmuniCacy6}
Let $f\colon X\to Y$ be morphisms
of regular scheme of finite type
over $S$.
Let $C\subset FT^*X$ be a closed conical
subset.
Assume that $Y$ is smooth over $S$.
Then, the following conditions are equivalent:

{\rm (1)}
$(1,f)$ is $C$-acyclic over $S$.

{\rm (2)}
The inverse image of
$C$ by $FT^*Y\times_YX\to FT^*X$
is a subset of the image
of the injection $FT^*S\times_SX\to FT^*Y\times_YX$.
\end{lm}

\proof{
Since the kernels
${\rm Ker}(FT^*X\times_X(FT^*Y\times_YX)
\to FT^*X)$
and
${\rm Ker}(FT^*X\times_X(FT^*Y\times_YX)
\to FT^*(X\times_SY)\times_{X\times_SY}X)$
are identified with the images of 
$FT^*Y\times_YX$
and 
$FT^*S\times_SX$,
the assertion follows.
\qed

}

\begin{lm}\label{lmCacycsm}
Let $h\colon W\to X$ and
$f\colon W\to Y$ be morphisms
of regular schemes of finite type over $S$.
Assume that $Y$ is smooth over $S$.

{\rm 1.}
For 
%
$C=FT^*X$,
the following conditions
are equivalent:

{\rm (1)}
The pair $(h,f)$ is $C$-acyclic over $S$.

{\rm (2)}
The morphism
$(h,f)\colon W\to X\times_SY$ is smooth
on a neighborhood of $W_{{\mathbf F}_p}$.

{\rm 2.}
If $C$ is the $0$-section $FT^*_XX$,
the following conditions
are equivalent:

{\rm (1)}
The pair $(h,f)$ is $C$-acyclic over $S$.

{\rm (2)}
The kernel
${\rm Ker}(FT^*Y\times_YW\to FT^*W)$
is a subset of the image of
${\rm Ker}(FT^*S\times_SW\to FT^*X\times_XW)$
on a neighborhood of $W_{{\mathbf F}_p}$.
\end{lm}

\proof{
1.
The condition (1)
holds if
and only if
$FT^*(X\times_SY)
\times_{X\times_SY}W\to FT^*W$
is an injection.
This is equivalent to (2)
by \cite[Proposition 2.8]{FW}.

2.
This is proved by the same argument
as in the proof of Lemma \ref{lmuniCacy}.2.
\qed

}

\begin{lm}\label{lmCacycq}
Let 
$$\xymatrix{
X'\ar[d]_g&&\\
X&W\ar[l]_h\ar[ul]_{h'}\ar[r]^f&Y}$$ 
be a commutative diagram
of regular schemes
of finite type over $S$.
Assume that $Y$ is smooth over $S$.
Let $C\subset FT^*X$ be a closed conical subset.
Let $q\colon X'\to S$ be the canonical
morphism and assume that
$(g,q)$ is $C$-acyclic over $S$.
Then, the following conditions are equivalent:

{\rm (1)} $(h,f)$ is $C$-acyclic over $S$.

{\rm (2)}
$(h',f)$ is $g^\circ C$-acyclic over $S$.
\end{lm}


\proof{
(1)$\Rightarrow$(2):
Since
${\rm Ker}((FT^*X\times_XW)\times_W(FT^*Y\times_YW)
\to FT^*W)$ 
equals the inverse image of
${\rm Ker}((FT^*X'\times_{X'}W)
\times_W(FT^*Y\times_YW)
\to FT^*W)$,
by taking the image of (\ref{eqCacyc})
by 
$(FT^*X\times_XW)\times_W(FT^*Y\times_YW)
\to (FT^*X'\times_{X'}W)\times_W(FT^*Y\times_YW)$,
we obtain an inclusion
\begin{align}
(g^\circ C\times_W&\,  (FT^*Y\times_YW))
\cap {\rm Ker}(
(FT^*X'\times_{X'}W)\times_W(FT^*Y\times_YW)
\to FT^*W)
\label{eqCacycW}
\\
&
\subset
{\rm Ker}((FT^*X'\times_{X'}W)
\times_W(FT^*Y\times_YW)\to
FT^*(X'\times_SY)\times_{X'\times_SY}W).
\nonumber
\end{align}

(2)$\Rightarrow$(1):
Since the kernel
of $(FT^*X'\times_{X'}W)\times_W
(FT^*Y\times_YW)\to 
FT^*(X'\times_SY)\times
_{X'\times_SY}W$
equals the image of
$FT^*S\times_SW\to 
(FT^*X'\times_{X'}W)
\times_W
(FT^*Y\times_YW)$,
the intersection of
its inverse image with
the kernel
of $(FT^*X\times_XW)\times_W
(FT^*Y\times_YW)\to FT^*W$
is a subset of the image of
the kernel of
$(FT^*X\times_XW)\times_W
(FT^*S\times_SW)\to 
(FT^*X'\times_{X'}W)$.
Hence the inclusion (\ref{eqCacycW})
implies that
the left hand side of (\ref{eqCacyc})
is a subset of the intersection
$(h^*C\times_W(FT^*S\times_SW))
\cap
{\rm Ker}((FT^*X\times_XW)\times_W
(FT^*S\times_SW)\to (FT^*X'\times_{X'}W))$.
Hence by the assumption
that $(g,q)$ is $C$-acyclic over $S$,
the assertion follows.
\qed

}

\begin{lm}\label{lmCacych}
Let $$\xymatrix{
X&W\ar[l]_h\ar[r]^f\ar[rd]_{f'}&Y\ar[d]^g\\
&&Y'}$$ 
be a commutative diagram
of regular scheme of finite type over $S$.
Assume that $Y$ and $Y'$ are smooth over $S$
and that $g$ is smooth.
Let $C\subset FT^*X$ be a closed conical subset.
Assume that $(h,f)$ is $C$-acyclic
over $S$.
Then, 
$(h,f')$ is $C$-acyclic
over $S$.
\end{lm}

\proof{
Since $FT^*Y'\times_{Y'}Y\to FT^*Y$
is an injection and
the inverse image of the kernel
${\rm Ker}((FT^*X\times_XW)\times_W(FT^*Y\times_YW)\to
FT^*(X\times_SY)\times_{X\times_SY}W)$
by 
$(FT^*X\times_XW)\times_W(FT^*Y'\times_{Y'}W)\to
(FT^*X\times_XW)\times_W(FT^*Y\times_YW)$
is 
${\rm Ker}((FT^*X\times_XW)\times_W(FT^*Y'\times_{Y'}W)\to
FT^*(X\times_SY')\times_{X\times_SY'}W)$,
the assertion follows.
\qed

}

%
%
%
%
%
%

\begin{lm}\label{lmCacycpr}
Let $g\colon X'\to X$ be a morphism
of regular schemes of finite type over $S$.
Let $C'\subset FT^*X'$ be a closed
conical subset and assume that $g$ is
proper on the base $B'$ of  $C'$.
Let $C=g_\circ C'$.
Let $$
\xymatrix{
X'\ar[d]_g&W'\ar[l]_{h'}\ar[rd]^{f'}\ar[d]_{g'}&\\
X&W\ar[l]_h\ar[r]^f &Y}
$$ be a commutative
diagram of separated morphisms
of schemes of finite
type over $S$ with cartesian square.
Assume that $W$ is regular and 
that $Y$ is smooth over
$S$.
Assume that
$(h,f)$ is $C$-acyclic over $S$
and that $S$ is excellent.

{\rm 1.}
The morphism
$h$ is transversal to $g$ on a neighborhood
of $h'^{-1}(B')$.

{\rm 2.}
There exists a regular neighborhood $U'
\subset W'$ of $h'^{-1}(B')$
on which
$(h',f')$ is $C'$-acyclic over $S$.
\end{lm}

\proof{1.
The morphism $h\colon W\to X$
is $C$-transversal by Lemma \ref{lmuniCacy}.4.
Hence the assertion follows from
Lemma \ref{lmg0C}.1.

2.
By 1,
there exists a regular neighborhood $U'
\subset W'$ of $h'^{-1}(B')$ such that
the sequence $0\to FT^*X\times_XU'\to 
FT^*X'\times_{X'}U'\oplus
FT^*W\times_WU'
\to FT^*U'$ is exact.
Then,
the canonical morphism
${\rm Ker}((FT^*X\times_XU')
\times_{U'}(FT^*Y\times_YU')
\to FT^*W\times_WU')
\to 
{\rm Ker}((FT^*X'\times_{X'}U')
\times_{U'}(FT^*Y\times_YU')
\to FT^*U')$ is an isomorphism.
Since
${\rm Ker}((FT^*X\times_XU')
\times_{U'}(FT^*Y\times_YU')
\to FT^*(X\times_SY)\times_{X\times_SY}U')
\to 
{\rm Ker}((FT^*X'\times_{X'}U')
\times_{U'}(FT^*Y\times_YU')
\to FT^*(X'\times_SY)\times_{X'\times_SY}U')$ is an isomorphism,
the assertion follows from the definition
of $C=g_\circ C'$.
\qed

}

%
%
%
%
%
%
%

\begin{df}\label{dfsat}
Let $C\subset FT^*X$ be a closed conical subset.

{\rm 1.}
We say that $C$ is $S$-saturated
if $C$ is stable under the action of
$FT^*S\times_SX$.

{\rm 2.}
Let $C^{\rm sat}\subset FT^*X$ denote the
smallest $S$-saturated
closed conical subset
containing $C$
as a subset.
\end{df}

If the morphism $FT^*S\times_SX\to FT^*X$ is 0,
then every closed conical subset $C\subset FT^*X$ is
$S$-saturated and hence $C=C^{\rm sat}$.
In particular, if $S={\rm Spec}\, k$
for a perfect field $k$,
every $C$ is $S$-saturated.
The saturation
$C^{\rm sat}$ is given by the closure
$\overline{C+FT^*S\times_SX}$.
If
$C^{\rm sat}=
C+FT^*S\times_SX,$
the $C$-acyclicity over $S$
and 
the $C^{\rm sat}$-acyclicity over $S$
are equivalent to each other
since the kernels in 
(\ref{eqCacyc}) are stable
under the actions of
$FT^*S\times_SW$.

\begin{lm}\label{lmCCbar}
Let $h\colon W\to X$ and $f\colon W\to Y$
be morphisms of regular schemes of
finite type over $S$
and assume that $Y$ is smooth over $S$.
Let $C\subset FT^*X$ be an $S$-saturated
closed conical subset.
%
%
%
%
Then the subset
$\{w_{{\mathbf F}_p}\in W\mid
\text{\rm 
$(h,f)$ is $C$-acyclic over $S$
at $w$}\}\cup
W_{{\mathbf Z}[\frac 1p]}\subset W$
is an open subset.
\end{lm}

\proof{
%
%
Since $C$ is $S$-saturated,
by the split short exact sequence
$0\to 
FT^*S\times_S(X\times_SY)
\to
FT^*X\times_SFT^*Y
\to
FT^*(X\times_SY)\to 0$,
the image
${\rm Im}(h^*C\times_W
(FT^*Y\times_YW))
\subset 
FT^*(X\times_SY)$
is a closed subset.
Hence the assertion follows from
(\ref{eqCacyc})$'$ and
Lemma \ref{lmCtrans}.1.
\qed

}

\begin{lm}\label{lmCacychs}
Let $h\colon W\to X$
be a morphism of regular schemes
of finite type over $S$ and
$q\colon W\to S$ be the canonical
morphism.
Let $C\subset FT^*X$ be a closed conical subset.
%
Assume
that 
$h$ is $C$-transversal and
$C$ is $S$-saturated.
Then
$(h,q)$ is $C$-acyclic over $S$.
%
%
%
\end{lm}

\proof{
Assume that $(a,b)
\in  h^*C\times_W(FT^*S\times_SW)$ 
is in the kernel
of $(FT^*X\times_XW)\times_W
(FT^*S\times_SW)\to FT^*W$.
Since $C$ is $S$-stable,
we have $a+b
\in h^*C\subset FT^*X\times_XW$
and this is contained in the kernel of
$FT^*X\times_XW\to FT^*W$.
Since $h$ is assumed $C$-transversal,
we have $a+b=0$ 
and hence 
$(a,b)$ is in the kernel of
$(FT^*X\times_XW)\times_W
(FT^*S\times_SW)\to FT^*X\times_XW$.
\qed

}

\begin{lm}\label{lmCacycS0}
Let $$\xymatrix{
X&W\ar[l]_h\ar[r]^f\ar[rd]_{f_0}&Y\ar[d]^g\ar[r]&S\ar[d]\\
&&Y_0\ar[r]&S_0}$$ 
be a commutative diagram
of regular noetherian schemes
satisfying the finiteness
condition {\rm (F)}
with cartesian square.
Assume that $X,W$ and $Y$ are of finite type over $S$
and that $S$ and $Y_0$ are smooth over $S_0$.
We consider the following conditions:

{\rm (1)}
$(h,f)$
is $C$-acyclic over $S$.

{\rm (2)}
$(h,f_0)$
is $C$-acyclic over $S_0$.

\noindent
We have
{\rm (1)}$\Rightarrow${\rm (2)}.
Conversely, if $C$ is $S$-saturated,
we have
{\rm (2)}$\Rightarrow${\rm (1)}.
\end{lm}

\proof{
{\rm (1)}$\Rightarrow${\rm (2)}:
Since $FT^*Y_0\times_{Y_0}Y\to FT^*Y$
is an injection and
since $X\times_{S_0}Y_0=X\times_SY$,
the assertion follows.

{\rm (2)}$\Rightarrow${\rm (1)}:
Let $q\colon W\to S$ be
the canonical morphism.
By Lemma \ref{lmCacycq},
it suffices to show that
$(h,q)$ is $C$-acyclic over $S$
and that $(1,f)$ is $h^\circ C$-acyclic over $S$.
Since $C$ is assumed $S$-saturated,
$(h,q)$ is $C$-acyclic by
Lemma \ref{lmuniCacy}.4 and
Lemma \ref{lmCacychs}.
Hence by
replacing $X$ and $C$
by $W$ and $h^\circ C$,
we may assume $X=W$ and $h=1_X$.
By Lemma \ref{lmuniCacy6},
it suffices to show that the condition (2)
there for $X\to Y_0$ over $S_0$ implies that
for $X\to Y$ over $S$.
Since $FT^*Y$ is generated by the images
of $FT^*Y_0$ and $FT^*S$,
this follows from the assumption
that $C$ is $S$-saturated.
%
%
\qed

}

\begin{lm}\label{lmsat}
Assume that $X$ is smooth
over $S$.
Then,
the mapping
$$\{
\text{\rm closed conical subsets of }
T^*X/S|_{X_{{\mathbf F}_p}}\}
\to 
\{ \text{\rm closed conical subsets of } FT^*X\}$$
sending $\overline C$ to
its inverse image is an injection
and the image consists of
saturated closed conical subsets.
\end{lm}

\proof{
Since the base change
$F^*(T^*X/S)|_{X_{{\mathbf F}_p}}
\to 
T^*X/S|_{X_{{\mathbf F}_p}}$
is a homeomorphism,
the assertion follows from the exact sequence
\begin{equation*}
0\to FT^*S\times_SX\to
FT^*X\to F^*T^*X/S|_{X_{{\mathbf F}_p}}\to 0.
\leqno{\rm (\ref{eqFT})}
\end{equation*}
\qed

}

\begin{lm}\label{lmCCbarsm}
Let $h\colon W\to X$ and $f\colon W\to Y$
be morphisms of smooth schemes over $S$.
Let $C\subset FT^*X$ be an $S$-saturated
closed conical subset
and let $\overline C\subset T^*X/S$
be the corresponding closed conical subset.

For $w\in W_{{\mathbf F}_p}$,
the following conditions are equivalent:

{\rm (1)}
$(h,f)$ is $C$-acyclic over $S$
at $w$.

{\rm (2)}
The fiber at $w$ of the intersection
$$((\overline C\times_XW)\times_W
(T^*Y/S\times_YW))
\cap 
{\rm Ker}
((T^*X/S\times_XW)\times_W
(T^*Y/S\times_YW)
\to T^*W/S)$$
is a subset of $\{0\}$.
\end{lm}

\proof{
The kernel ${\rm Ker}
((FT^*X\times_XW)\times_W
(FT^*Y\times_YW)
\to FT^*W)$
is the extension of
${\rm Ker}
((F^*T^*X/S|_{X_{{\mathbf F}_p}}\times_XW)\times_W
(F^*T^*Y/S|_{Y_{{\mathbf F}_p}}\times_YW)
\to F^*T^*W/S|_{W_{{\mathbf F}_p}})$
by the diagonal image of
$FT^*S\times_SW$
and
$T^*X/S
\times_S
T^*Y/S
\to
T^*(X\times_SY)/S$ is an isomorphism.
Hence the assertion follows.
\qed

}

\subsection{${\cal F}$-acyclicity over $S$}\label{ssFacycS}

Lemma \ref{lmcn}
suggests the following definition.

\begin{df}\label{dfFacyc}
Let $h\colon W\to X$ and
$f\colon W\to Y$ be separated morphisms
of regular 
schemes of finite type over $S$.
Assume that $Y$ is smooth over $S$.
Let ${\cal F}$ be a constructible sheaf on $X$.

{\rm 1.}
We say that $(h,f)$ is ${\cal F}$-acyclic over $S$
if the following condition is satisfied:

For every $S$-acyclic constructible sheaf ${\cal G}$ on $Y$,
the morphism
$(h,f)\colon W\to X\times_SY$ 
is ${\cal F}\boxtimes {\cal G}$-transversal
on an open neighborhood of $W_{{\mathbf F}_p}$.

{\rm 2.}
Let $B\subset X$ be the support of ${\cal F}$.
We say that $(h,f)$ is universally 
${\cal F}$-acyclic over $S$
if the following condition is satisfied:

For every morphism $g\colon Y'\to Y$
of smooth schemes over $S$
and the commutative diagram
$$\begin{CD}
\xymatrix{
&W'\ar[r]^{f'}\ar[ld]_{h'}\ar[d]^{g'}&Y'\ar[d]^g\\
X&W\ar[l]_h\ar[r]^f &Y}
\end{CD}$$
with cartesian square,
there exists an open neighborhood $U'
\subset W'$
of $h'^{-1}(B)_{{\mathbf F}_p}$
such that $g$ is transversal to
$f$ on $U'$  
and 
the pair $(h',f')$ 
is ${\cal F}$-acyclic over $S$
on $U'$.

If $Y=S$, we say that
$h$ is universally ${\cal F}$-transversal
over $S$.
\end{df}


In Definition \ref{dfFacyc}.1,
it is more natural to define a Zariski local condition
on $Y$. However, since
the $C$-acyclicity over $S$
is a local condition on $Y$,
this does not affect the definition
of $S$-micro support in
Definition \ref{dfSms}.
Thus, we adopt a simpler formulation here.
If $Y=S$, 
the condition that
$(h,f)$ is ${\cal F}$-acyclic over $S$
means that $h$ is ${\cal F}$-transversal
on a neighborhood of $W_{{\mathbf F}_p}$
since 
$S$-acyclicity of ${\cal G}$ on $Y=S$ means
that ${\cal G}$ is locally constant
on a neighborhood of $S_{{\mathbf F}_p}$
by Lemma \ref{lmmslc}.
In particular, 
for every ${\cal F}$ on $X$
and the morphism $f\colon X\to S$,
the pair of $(1_X,f)$
is ${\cal F}$-acyclic
over $S$.
The condition that
$h\colon W\to X$ 
is universally ${\cal F}$-transversal
over $S$ means that
for every smooth scheme $Y$ over $S$
and for every $S$-acyclic sheaf ${\cal G}$
on $Y$,
the base change
$h'\colon W\times_SY\to X\times_SY$
is ${\cal F}\boxtimes {\cal G}$-transversal
on a neighborhood of $(W\times_SY)_{{\mathbf F}_p}$.

\begin{lm}\label{lmk}
Let $f\colon W\to Y$ be a 
smooth morphism of smooth schemes over $S$.
Let ${\cal F}$ be a constructible sheaf on $W$.
We consider the following conditions:

{\rm (1)} $(1_W,f)$ is universally 
${\cal F}$-acyclic over $S$.

{\rm (2)} $f$ is universally ${\cal F}$-acyclic.

{\rm 1.}
We have {\rm (2)}$\Rightarrow${\rm (1)}.

{\rm 2.}
If $S={\rm Spec}\, k$ for a perfect field
$k$ of characteristic $p$,
the conditions {\rm (1)} and {\rm (2)} are
equivalent.
\end{lm}

The condition that $f$ is smooth is satisfied
on a neighborhood of the intersection
of the support of $h^*{\cal F}$
with $W_{{\mathbf F}_p}$
if 
$(h,f)$ is $C$-acyclic over $S$
for a closed conical subset $C\subset 
FT^*X$ such that the intersection
of the support of ${\cal F}$ 
with $X_{{\mathbf F}_p}$ is
a subset of the base of $C$
and 
if $X$ is smooth over $S$
by Lemma \ref{lmuniCacy}.3.

\proof{
1.
By replacing $W\to Y$
by the base change by
a morphism $Y'\to Y$ of smooth schemes over $Y$,
it suffices to show that $(1_W,f)$ is 
${\cal F}$-acyclic over $S$.
Let ${\cal G}$ be an $S$-acyclic constructible sheaf
on $Y$.
We show that
$\gamma=(1,f)\colon W\to W\times_SY$
is ${\cal F}\boxtimes {\cal G}$-transversal.
We consider the cartesian diagram
$$\begin{CD}
W\times_SY@<\gamma<< W\\
@V{f\times 1_Y}VV@VVfV\\
Y\times_SY@<\delta<< Y.
\end{CD}$$
Since $f$ is assumed universally ${\cal F}$-acyclic,
the morphism
$f\times 1_Y$ is
${\cal F}\boxtimes\Lambda$-acyclic.
Since $\delta$ is a morphism of smooth schemes
over $Y$,
the diagonal $\delta$
is $\Lambda\boxtimes {\cal G}$-transversal.
Since $f$ is assumed smooth,
by Corollary \ref{corFacyctr},
the morphism
$\gamma$
is ${\cal F}\boxtimes{\cal G}$-transversal.

2.
It suffices to show that
$f$ is ${\cal F}$-acyclic.
Since every constructible sheaf on $Y$
is $S$-acyclic,
the assertion follows from Proposition \ref{prBG}.
\qed

}

\begin{lm}\label{lmFacyc0}
Let $X$ be a regular scheme of finite type over
$S$ and ${\cal F}$ be a constructible sheaf on $X$.
Let $h\colon W\to X$ be an open immersion
and define $f\colon W\to Y={\mathbf A}^1_S$
to be the composition of the canonical
morphism $W\to S$ with
the $0$-section $i\colon S\to Y$.
If $(h,f)$ is ${\cal F}$-acyclic over $S$,
on a neighborhood of $W_{{\mathbf F}_p}$,
we have $h^*{\cal F}=0$.
\end{lm}

\proof{
The constructible
sheaf ${\cal G}=i_*\Lambda_S$
on $Y$ is $S$-acyclic by Lemma \ref{lmSacyfun}.1.
Hence the $0$-section $i_W\colon 
W\to {\mathbf A}^1_X=X\times_SY$
is ${\cal F}\boxtimes i_*\Lambda_S$-transversal
on a neighborhood of $W_{{\mathbf F}_p}$.
By Lemma \ref{lmFtrans}.4,
we have
${\cal F}=0$
on a neighborhood of $W_{{\mathbf F}_p}$.
\qed

}

\begin{lm}\label{lmFacyc}
Let 
$$\xymatrix{
X'\ar[d]_g&&\\
X&W\ar[l]_h\ar[ul]_{h'}\ar[r]^f&Y}$$ 
be a commutative diagram
of separated morphisms regular schemes
of finite type over $S$.
Assume that $Y$ is smooth over $S$.
Let ${\cal F}$ be a constructible sheaf on $X$.
Assume that 
$g$ is universally ${\cal F}$-transversal.
Then, the following conditions
are equivalent:

{\rm (1)} $(h,f)$ is ${\cal F}$-acyclic over $S$.

{\rm (2)} $(h',f)$ is $g^*{\cal F}$-acyclic over $S$.
\end{lm}

\proof{
Let ${\cal G}$ be an $S$-acyclic constructible sheaf on $Y$.
Since $g$ is universally ${\cal F}$-transversal,
the morphism 
$g\times 1_Y\colon X'\times_SY \to X\times_SY$ is 
${\cal F}\boxtimes {\cal G}$-transversal
on a neighborhood of $(X'\times_SY)_{{\mathbf F}_p}$.
Since $Y$ is smooth over $S$,
$(g\times 1)^!\Lambda$ is 
isomorphic to $\Lambda(d)[2d]$
for the relative dimension $d$ of $g$
by Lemma \ref{lmpurity}.1.
Hence by Lemma \ref{lmFtrtr},
$W\to  X\times_SY$
is ${\cal F}\boxtimes {\cal G}$-transversal
on a neighborhood of $W_{{\mathbf F}_p}$
if and only if
$W\to  X'\times_SY$
is $g^*{\cal F}\boxtimes {\cal G}$-transversal
on a neighborhood of $W_{{\mathbf F}_p}$.
\qed

}

\begin{lm}\label{lmFacycfun}
Let $$\xymatrix{
X&W\ar[l]_h\ar[r]^f\ar[rd]_{f'}&Y\ar[d]^g\\
&&Y'}$$ 
be a commutative diagram
of separated morphisms of
regular scheme of finite type over $S$.
Assume that $Y$ is smooth over $S$.
Let ${\cal F}$ be a constructible sheaf on $X$.
Assume that
$(h,f)$ is ${\cal F}$-acyclic over $S$
and let $g\colon Y\to Y'$
be a smooth morphism of
smooth schemes over $S$.
Then $(h,f')$ is ${\cal F}$-acyclic over $S$.
\end{lm}

\proof{
Let ${\cal G}'$ be an $S$-acyclic
constructible sheaf on $Y'$.
Then, the pull-back
${\cal G}$ on $Y$
of ${\cal G}'$ is $S$-acyclic
by Lemma \ref{lmSacyfun}.2.
Hence $(h,f)\colon W\to X\times_SY$
is ${\cal F}\boxtimes {\cal G}$-transversal
on a neighborhood of $W_{{\mathbf F}_p}$.
Since $X\times_SY\to X\times_SY'$
is smooth,
$(h,f')\colon W\to X\times_SY'$
is also
${\cal F}\boxtimes {\cal G}'$-transversal
on a neighborhood of $W_{{\mathbf F}_p}$
by  Lemma \ref{lmFtrans}.1 and Lemma \ref{lmFtrtr}.
\qed

}

\begin{lm}\label{lmFacycg}
Let $h\colon W\to X$ and $f\colon W\to Y$ 
be separated morphisms of regular schemes
of finite type over $S$.
%
Let
$$\xymatrix{
X'\ar[d]_g&W'\ar[l]_{h'}\ar[rd]^{f'}\ar[d]_{g'}&\\
X&W\ar[l]_h\ar[r]^f &Y}
$$
be a commutative diagram
of schemes of finite type over $S$
with cartesian square.
Assume that $X,W$ and $X'$ are regular.
Let ${\cal F}'$ be
a constructible sheaf ${\cal F}'$ on $X'$
and $B'\subset X'$ be the support.
Assume that $g$ is proper 
on $B'$
and let ${\cal F}=g_*{\cal F}'$.
Assume further that $h$ 
is transversal to $g$ on
a neighborhood of $h'^{-1}(B')_{{\mathbf F}_p}
\subset W'$.
We consider the following conditions:

{\rm (1)} $(h,f)$ is ${\cal F}$-acyclic over $S$.

{\rm (2)} 
$(h',f')$ is ${\cal F}'$-acyclic over $S$,

We have the implication {\rm (2)}$\Rightarrow${\rm (1)}.
If $g$ is finite
on $B$, 
then these conditions are equivalent.
\end{lm}


\proof{
%
Consider the cartesian diagram
$$\begin{CD}
W'@>{h'_Y}>>X'\times_SY\\
@V{g'}VV@VV{g\times 1_Y}V\\
W@>{h_Y}>> X\times_SY.
\end{CD}$$
Since $h$ is transversal to $g$
on a neighborhood $U'$
of $h'^{-1}(B')_{{\mathbf F}_p}
\subset W'$,
the base change
$h_Y$ is transversal to $g\times 1_Y$
on $U'$.
Hence 
the base change morphism
$g'^*h_Y^!\Lambda\to h'^!_Y\Lambda$ 
is an isomorphism by Lemma \ref{lmpurity}.2.
Let ${\cal G}$ be an $S$-acyclic constructible sheaf
on $Y$.
Since $g\times 1_Y$ is proper
on the support of
${\cal F}'\boxtimes{\cal G}$,
the ${\cal F}'\boxtimes{\cal G}$-transversality
on a neighborhood of $W'_{{\mathbf F}_p}$
implies
the $g_*{\cal F}'\boxtimes{\cal G}$-transversality
on a neighborhood of $W_{{\mathbf F}_p}$
by Lemma \ref{lmFtrpr} {\rm (2)}$\Rightarrow${\rm (1)}.
Hence we have {\rm (2)}$\Rightarrow${\rm (1)}.
Conversely, if $g$ is finite on $B$,
we have {\rm (1)}$\Rightarrow${\rm (2)}
by Lemma \ref{lmFtrpr} {\rm (1)}$\Rightarrow${\rm (2)}.
\qed

}

\begin{lm}\label{cnacyc}
Let $$\xymatrix{
X&W\ar[l]_h\ar[r]^f\ar[rd]_{f_0}&Y\ar[d]^g\ar[r]&S\ar[d]\\
&&Y_0\ar[r]&S_0}$$ 
be a commutative diagram
of separated morphisms of 
regular noetherian schemes
satisfying the finiteness
condition {\rm (F)}
with cartesian square.
Assume that $X,W$ and $Y$ are of finite type over $S$
and that $S$ and $Y_0$ are smooth over $S_0$.
Let ${\cal F}$ be a constructible sheaf on $X$.
If $(h,f)$ is ${\cal F}$-acyclic
over $S$,
then $(h,f_0)$ is ${\cal F}$-acyclic
over $S_0$.
\end{lm}

\proof{
For an $S_0$-acyclic constructible sheaf
${\cal G}_0$
on $Y_0$,
its pull-back ${\cal G}$ on $Y$ 
is $S$-acyclic by Lemma \ref{lmSacyfun}.3.
Hence 
if $(h,f)$ is ${\cal F}$-acyclic over $S$,
then
the morphism
$W\to X\times_SY=X\times_{S_0}Y_0$
is ${\cal F}\boxtimes{\cal G}=
{\cal F}\boxtimes{\cal G}_0$-transversal
on a neighborhood of $W_{{\mathbf F}_p}$.
\qed

}

\subsection{$S$-micro support}\label{ssSms}

Let $S$ be a regular excellent noetherian
scheme over ${\mathbf Z}_{(p)}$
satisfying the finiteness condition (F).

\begin{df}\label{dfSms}
Let $X$ be a regular scheme
of finite type over $S$
and let ${\cal F}$ be a constructible sheaf on $X$.
Let $C\subset FT^*X$ be
a closed conical subset.
We say that ${\cal F}$ is $S$-micro supported
on $C$ if the following condition is satisfied:

For every pair of separated 
morphisms $h\colon W\to X$ and
$f\colon W\to Y$ 
of regular schemes of finite type over $S$
such that $Y$ is smooth over $S$,
the $C$-acyclicity of $(h,f)$ over $S$
implies
the ${\cal F}$-acyclicity of $(h,f)$
over $S$.
\end{df}

We may replace
the ${\cal F}$-acyclicity over $S$ in the definition
by the universal
${\cal F}$-acyclicity over $S$
by Lemma \ref{lmSbase} below.
If $S={\rm Spec}\, k$ for a perfect field
of characteristic $p$,
then being $S$-micro supported is
equivalent to being micro supported
by \cite[Proposition 8.13]{CC}.

\begin{lm}\label{lmSbase}
Let $X$ be a regular scheme
of finite type over $S$
and let ${\cal F}$ be a constructible sheaf on $X$.
Let $C\subset FT^*X$ be
a closed conical subset
and assume that
${\cal F}$ is $S$-micro supported on $C$.

{\rm 1.}
The intersection of
$X_{{\mathbf F}_p}$ with
the support of ${\cal F}$
is a subset of the base $B$ of $C$.

{\rm 2.}
Let $h\colon W\to X$ and
$f\colon W\to Y$ be morphisms 
of regular schemes of finite type over $S$.
If $(h,f)$ is $C$-acyclic over $S$,
then
$(h,f)$ is universally ${\cal F}$-acyclic over $S$
and $h$ is 
universally ${\cal F}$-transversal
over $S$.

\end{lm}

\proof{
1.
Let $W=X\sm B$ be the complement.
Let $h\colon W\to X$ be the immersion
and $f\colon W\to Y={\mathbf A}^1_S$
be the composition of the canonical
morphism $W\to S$ with
the $0$-section $i\colon S\to Y$
as in Lemma \ref{lmFacyc0}.
Then, $(h,f)$ is $C$-acyclic over $S$ and
hence ${\cal F}$-acyclic over $S$.
Thus, we have $h^*{\cal F}=0$
on a neighborhood of $W_{{\mathbf F}_p}$
by Lemma \ref{lmFacyc0}.

2.
Assume that $(h,f)$ is
$C$-acyclic over $S$.
By 1 and Lemma \ref{lmuniCacy}.1,
on a neighborhood of
the intersection of
$W_{{\mathbf F}_p}$ with the support
of the pull-back of ${\cal F}$,
$(h,f)$ is universally 
$C$-acyclic over $S$
and hence $(h,f)$ is universally
${\cal F}$-acyclic over $S$.

Let $q\colon W\to S$ be
the canonical morphism.
Then, $(h,q)$
is $C$-acyclic over $S$
by Lemma \ref{lmuniCacy}.4.
Hence $(h,q)$
is universally ${\cal F}$-acyclic over $S$.
This means that
the morphism $h$ is 
universally ${\cal F}$-transversal
over $S$.
\qed

}

\begin{lm}\label{lmS}
Let $X$ be a regular scheme
of finite type over $S$.
Let ${\cal F}$ be a constructible sheaf on $X$
and $C\subset FT^*X$
be a closed conical subset.

{\rm 1.}
If ${\cal F}$ is $S$-micro supported on $C$
and if $C$ is $S$-saturated,
then ${\cal F}$ is micro supported on $C$.

{\rm 2.}
Assume that
${\cal F}$ is micro supported on $C$.
Assume further that the following condition
is satisfied:

\noindent
{\rm (T)}
For every smooth scheme $Y$
over $S$,
for every $S$-acyclic $C'\subset FT^*Y$
closed conical subset
and for every constructible
sheaf ${\cal G}$ on $Y$ micro supported
on $C'$,
the external tensor product
${\cal F}\boxtimes {\cal G}$
is micro supported on ${\rm pr}_1^\circ C
+{\rm pr}_2^\circ C'$.

\noindent
Then ${\cal F}$ is $S$-micro supported on $C$.

{\rm 3.}
Let $S\to S_0$ be a smooth morphism
of regular noetherian schemes
over ${\mathbf Z}_{(p)}$.
If ${\cal F}$ is $S$-micro supported on $C$
and 
if $C$ is $S$-saturated,
then
${\cal F}$ is $S_0$-micro supported on $C$.
\end{lm}

The condition (T) in Lemma \ref{lmS}.2
is satisfied if $S$ is smooth over a perfect field $k$
of characteristic $p>0$
by By Lemma \ref{lmcn}.

\proof{
1.
Let $h\colon W\to X$
be a separated morphism of
regular schemes
of finite type over $S$.
Assume that
$h$ is $C$-transversal
and let $q\colon W\to S$ be 
the canonical morphism.
Then, $(h,q)$ is $C$-acyclic over $S$
by Lemma \ref{lmCacychs}.
Hence $h\colon W\to X=X\times_SS$
is ${\cal F}$-transversal
on a neighborhood of $W_{{\mathbf F}_p}$.

2.
Let $h\colon W\to X$
and $f\colon W\to Y$
be morphisms over $S$ such that
$(h,f)$ is
$C$-acyclic over $S$.
Let 
${\cal G}$ be an $S$-acyclic constructible sheaf
on $Y$.
Assume that ${\cal G}$ is micro supported
on an $S$-acyclic closed conical subset
$C'\subset FT^*Y$.
By assumption,
the tensor product
${\cal F}\boxtimes{\cal G}$
on $X\times_SY$
is micro supported on ${\rm pr}_1^\circ C+
{\rm pr}_2^\circ C'$.
By Lemma \ref{lmuniCacy}.5,
the morphism
$h\colon W\to X\times_SY$
is 
${\rm pr}_1^\circ C+
{\rm pr}_2^\circ C'$-transversal.
Hence
$h$ is
${\cal F}\boxtimes{\cal G}$-transversal
on a neighborhood of $W_{{\mathbf F}_p}$.

3.
Let $h\colon W\to X$ and
$f_0\colon W\to Y_0$ be morphisms
over $S_0$.
Assume that 
$Y_0\to S_0$ is smooth and that
$(h,f_0)$
is $C$-acyclic over $S_0$.
Consider $W$ as a scheme over $S$
by the composition $W\to X\to S$
and let $f\colon W\to Y=Y_0\times_{S_0}S$
be the induced morphism over $S$.
Then, $(h,f)$
is $C$-acyclic over $S$
by Lemma \ref{lmCacycS0}
and hence $(h,f)$
is universally ${\cal F}$-acyclic over $S$.
By Lemma \ref{cnacyc},
$(h,f_0)$
is universally ${\cal F}$-acyclic over $S_0$.
\qed

}

\begin{lm}\label{lmFperp}
Let $X$ be a regular scheme of
finite type over $S$.

{\rm 1.}
Every constructible sheaf ${\cal F}$
on $X$ is $S$-micro supported on $FT^*X$.

{\rm 2.}
Every locally constant constructible sheaf ${\cal F}$
on $X$ is $S$-micro supported
on the $0$-section $FT^*_XX$.

{\rm 3.}
Let ${\cal F}$ be a constructible sheaf on $X$.
Assume ${\cal F}$ is 
$S$-micro supported on the $0$-section
and that $FT^*S\times_SX\to FT^*X$
is the $0$-morphism.
Then ${\cal F}$ is locally constant
on a neighborhood of $X_{{\mathbf F}_p}$.

{\rm 4.}
Assume $X=S$. Then every constructible sheaf
on $S$ is $S$-micro supported on
the $0$-section $FT^*_SS\subset FT^*S$.
\end{lm}

Lemma \ref{lmFperp}.3 and 4.~imply
that the notion of $S$-micro support
can be a reasonable one
only for $S$-saturated conical closed subset.

\proof{
1.
Assume that
$(h,f)$ is $FT^*X$-acyclic over $S$.
Then,
by Lemma \ref{lmCacycsm}.1,
$(h,f)\colon W\to X\times_S Y$ is smooth
on a neighborhood of $W_{{\mathbf F}_p}$.
Hence for every constructible sheaf
${\cal F}$ on $X$ and 
every constructible sheaf ${\cal G}$ on $Y'$,
the pair $(h,f)$ is
${\cal F}\boxtimes{\cal G}$-transversal
on a neighborhood of $W_{{\mathbf F}_p}$
by Lemma \ref{lmFtrans}.1.

2.
We may assume ${\cal F}=\Lambda$.
Assume that
$(h,f)$ is $FT^*_XX$-acyclic over $S$.
Then the base change
$f\colon W\to Y$
is smooth modulo $S$
by Lemma \ref{lmuniCacy}.3
and ${\rm pr}_2\colon X\times_SY\to Y$
is smooth modulo $S$
by Lemma \ref{lmf0}.2.
Hence for every 
$S$-acyclic sheaf ${\cal G}$ on $Y$,
on a neighborhood of $W_{{\mathbf F}_p}$,
the morphism $f$ is ${\cal G}$-transversal 
and ${\rm pr}_2\colon X\times_SY\to Y$
is also ${\cal G}$-transversal
by Lemma \ref{lmsmacyc}.2.
Hence 
$W\to X\times_SY$
is $\Lambda\boxtimes{\cal G}$-transversal
on a neighborhood of $W_{{\mathbf F}_p}$
by Lemma \ref{lmFtrtr}.

3.
Assume that $FT^*S\times_SX\to FT^*X$
is the $0$-morphism.
By Lemma \ref{lmCacycsm}.2,
for every morphism
$h\colon W\to X$ of regular
schemes of finite type
and the canonical morphism
$q\colon W\to S$,
the pair $(h,q)$ is $C$-acyclic over
$S$ for the $0$-section
$C=FT^*_XX$.
Hence $h\colon W\to X=X\times_SS$
is ${\cal F}={\cal F}\boxtimes \Lambda$-transversal
on a neighborhood of $W_{{\mathbf F}_p}$.
This implies that ${\cal F}$ is locally constant
on a neighborhood of $X_{{\mathbf F}_p}$ by 
Proposition \ref{prlcc}.

4.
For $C=FT^*_SS$,
we have
$C\subset FT^*S\subset C+FT^*S$.
Hence the $C$-acyclicity over $S$
is equivalent to
the $FT^*S$-acyclicity over $S$.
Thus it follows from 1.
\qed

}

\begin{lm}\label{lmCbarms}
{\rm 1.}
Let $g\colon X'\to X$ be a separated morphism
of regular schemes of finite type over $S$
and $q\colon X'\to S$ be the canonical morphism.
Let ${\cal F}$ be a constructible sheaf on $X$.
If ${\cal F}$ is $S$-micro supported on $C$
and if $(g,q)$ is $C$-acyclic over $S$,
then $g^*{\cal F}$ is $S$-micro supported
on $g^\circ C$.

{\rm 2.}
Let $g\colon X'\to X$ be a morphism
of regular schemes of finite type over $S$.
Assume that $g$ is proper on the support of $C'
\subset FT^*X'$.
Let ${\cal F}'$ be a constructible sheaf on $X'$.
If ${\cal F}'$ is $S$-micro supported on $C'$,
then $Rg_*{\cal F}'$ is $S$-micro supported
on $C=g_\circ C'$.
\end{lm}

\proof{
1. 
Let $h\colon W\to X'$ and
$f\colon W\to Y$ be separated morphisms
of regular schemes of finite type over $S$
and assume that $Y$ is smooth over $S$.
Assume that $(h,f)$ is $g^\circ C$-acyclic over $S$.
Then by 
Lemma \ref{lmCacycq},
$(gh,f)$ is $C$-acyclic over $S$
and 
is ${\cal F}$-acyclic over $S$.
Since $(g,q)$ is $C$-acyclic over $S$,
the morphism $g$ is universally 
${\cal F}$-transversal over $S$
by Lemma \ref{lmSbase}.2.
Hence 
$(h,f)$ is $g^*{\cal F}$-acyclic over $S$
by Lemma \ref{lmFacyc}.

2.
Let $h\colon W\to X$ and $f\colon W\to Y$
be separated morphisms of regular schemes
of finite type over $S$
such that $Y$ is smooth over $S$
and $(h,f)$ is $C$-acyclic over $S$.
Let $$
\begin{CD}
X'@<{h'}<<W'\\
@VgVV@VV{g'}V\\
X@<h<<W
\end{CD}$$
be a cartesian diagram.
Then, by Lemma \ref{lmCacycpr},
there exists an open regular neighborhood
$U'\subset W'=W\times_XX'$
of the inverse image of $B'$
where $h$ is transversal to $g$ and
$(h',fg')$ is $C'$-acyclic over $S$.
Since $B'$ contains the intersection of
the support of ${\cal F}'$ and $X'_{{\mathbf F}_p}$,
$(h',fg')$ is ${\cal F}'$-acyclic over $S$ on 
a neighborhood of $U'_{{\mathbf F}_p}$.
Hence the assertion follows from
Lemma \ref{lmFacycg}.
\qed

}

\begin{lm}\label{lmSSXPS}
Let $i\colon X\to P$ be a closed immersion of
regular schemes of
finite type over $S$
and
let ${\cal F}$ be a 
constructible sheaf on $X$.
Let $C_P$ be a closed conical subset
of $FT^*P|_X
\subset FT^*P$
and define $C\subset FT^*X$ 
to be the closure of
the image of $C_P$
by the surjection
$FT^*P|_X\to FT^*X$.
If $i_*{\cal F}$ is $S$-micro supported on $C_P$,
then
${\cal F}$ is $S$-micro supported on $C$.
\end{lm}
\proof{
%
Let $h\colon W\to X$ and 
$f\colon W\to Y$ be 
separated morphisms
of regular schemes of finite type over $S$
and assume that
$(h,f)$ is $C$-acyclic over $S$.
We show that
$(h,f)$ is ${\cal F}$-acyclic over $S$.
By the same argument as in
the proof of Lemma \ref{lmSSXP},
we may assume that $h$ is an immersion.
Then $(h,f)\colon W\to X\times_SY$ is also 
an immersion.

Since the assertion is local on $W$,
by taking a lifting of a local basis of
the conormal bundle $T^*_W(X\times_SY)$,
we may extend the immersion $(h,f)$
to an immersion $(g,q)\colon
V\to P\times_SY$ of regular closed subscheme
transversal with the immersion $X\times_SY\to P\times_SY$.
Then,
the $C$-acyclicity over $S$ of $(h,f)$
implies the $C_P$-acyclicity over $S$ of $(g,q)$.
Since $i_*{\cal F}$ is $S$-micro supported on $C_P$,
this implies that $(g,q)$ is $i_*{\cal F}$-acyclic over $S$.
Since the intersection $W=V\cap (X\times_SY)$ 
is transversal,
if $i'\colon W\to V$ denote the immersion,
the morphism
$i'^*g^!\Lambda\to h^!\Lambda$
is an isomorphism.
Hence by
Lemma \ref{lmFacycg} (1)$\Rightarrow$(2),
the $i_*{\cal F}$-acyclicity over $S$ of $(g,q)$
implies the ${\cal F}$-acyclicity over $S$ of $(h,f)$.
Hence ${\cal F}$ is $S$-micro supported on $C$.
\qed

}

\begin{lm}\label{lmCms}
Let $X$ be a 
regular schemes of
finite type over $S$.
Let $C$ be a closed
conical subset of $FT^*X$
and ${\cal F}$ be a constructible sheaf on $X$.

{\rm 1.}
Let $X=\bigcup U_i$
be an open covering of $X$.
Then the following conditions are equivalent:

{\rm (1)}
${\cal F}$ is $S$-micro supported on $C$.

{\rm (2)}
${\cal F}|_{U_i}$ is $S$-micro supported on 
$C|_{U_i}$ for every $i\in I$.

{\rm 2.}
Let $U\subset X$ be an open subscheme
and let $Z=X\sm U$ be the complement.
Let $C'\subset FT^*U$ be
a closed subset.
Assume that ${\cal F}$
is $S$-micro supported on $C$
and that ${\cal F}|_U$
is $S$-micro supported on $C'$.
Assume further that $C$ and $C'$ are $S$-saturated.
Then,  ${\cal F}$
is $S$-micro supported on
the union $C_1=\overline {C'}\cup C|_Z$.
\end{lm}

\proof{
The same argument as in the proof of
Lemma \ref{lmmsloc} works using
Lemma \ref{lmCCbar}.
\qed

}

\medskip

For an $S$-saturated closed conical subset
on a smooth scheme over $S$,
the $S$-micro supported condition can be rephrased
as follows.

\begin{lm}\label{lmmssat}
Let $X$ be a smooth scheme over $S$
and ${\cal F}$ be a constructible sheaf on $X$.
Let $C\subset FT^*X$ be an $S$-saturated
closed conical subset and
let  $\overline C
\subset T^*X/S|_{X_{{\mathbf F}_p}}$
be the corresponding closed conical subset.
Then, the following conditions are equivalent:

{\rm (1)}
${\cal F}$ is $S$-micro supported on $C$.

{\rm (2)}
For every pair of morphisms 
$h\colon W\to X$ and
$f\colon W\to Y$ 
of smooth separated schemes of finite type over $S$
such that
the morphism
$(T^*X/S\times_XW)
\times_W
(T^*Y/S\times_YW)
\to
T^*W/S
$
is $h^*\overline C\times_W
(T^*Y/S\times_YW)|_{W_{{\mathbf F}_p}}
$-transversal,
the pair $(h,f)$ is ${\cal F}$-acyclic over $S$.
\end{lm}

\proof{
By Lemma \ref{lmuniCacy}.3,
it suffices to consider pairs $(h,f)$ such that
$W$ is smooth over $S$.
Hence the assertion follows from
Lemma \ref{lmCCbarsm}.
\qed

}

\begin{lm}\label{lmSSSS}
Let $k$ be a field of characteristic $p>0$
and $X$ be a smooth scheme over $S={\rm Spec}\, k$.
Let ${\cal F}$ be a constructible sheaf on $X$.
Let $C\subset FT^*X$ be an $S$-saturated
closed conical subset and
$\overline C\subset T^*X/k$
be the corresponding closed conical subset.
If ${\cal F}$ is micro supported
on $\overline C\subset T^*X/k$ in the sense of {\rm \cite{SS}},
then
${\cal F}$ is $S$-micro supported
on $C$.
\end{lm}

\proof{
By Lemma \ref{lmmssat},
the assertion follows from 
Lemma \ref{lmk}.
\qed

}
\medskip
The converse of Lemma \ref{lmSSSS} does not hold if $k$
is imperfect. For a concrete example, see Remark after Definition \ref{dfSss}.

\subsection{$S$-singular support}\label{ssSSS}

Let $S$ be a regular excellent noetherian
scheme over ${\mathbf Z}_{(p)}$
satisfying the finiteness condition (F).

\begin{df}\label{dfSss}
Let $X$ be a regular scheme
of finite type over $S$
and let ${\cal F}$ be
a constructible sheaf on $X$.

{\rm 1.}
We say that 
a closed conical subset 
(resp. an $S$-saturated closed conical subset)
$C\subset FT^*X$
is the $S$-singular support $SS_S{\cal F}$
(resp. the $S$-saturated singular support
$SS^{\rm sat}_S{\cal F}$)
of ${\cal F}$
if for every closed conical subset 
(resp. every $S$-saturated closed conical subset)
$C'\subset FT^*X$,
the inclusion $C\subset C'$
is equivalent to the condition
that ${\cal F}$ is $S$-micro supported on $C'$.

{\rm 2.}
If $X$ is smooth over $S$,
we call the closed conical subset
$\overline{SS}_S{\cal F}
\subset T^*X/S|_{{X_{{\mathbf F}_p}}}$ 
corresponding to
$SS^{\rm sat}_S{\cal F}
\subset FT^*X$
the $S$-relative singular support.
\end{df}

If $S={\rm Spec}\, k$ for a perfect field
of characteristic $p$,
the $S$-singular support and
the $S$-saturated singular support is the same
as the singular support.
If $k'$ is not perfect,
the $k'$-relative singular support
may be different from that defined
by Beilinson \cite{SS}
or by Hu--Yang \cite{HY}.
By \cite[Theorem  1.4 (iii)]{SS},
the singular support
defined by Beilinson commutes with base change.
In this sense, Beilinson's singular support is
essentially an invariant for 
smooth schemes over a perfect field.
On the other hand,
the relative singular support defined
in Definition \ref{dfSss} does not commute
with base change in general.

\medskip
\noindent
{\it
Example.}
Let $k$ be a perfect field
and $k'=k(y)$ be
the rational function
field with one variable.
Let ${\cal F}$ be the sheaf on
$X={\mathbf A}^1_{k'}
={\rm Spec}\, k'[x]$
defined by $t^p-t=y/x^p$.
Then, the 
singular support
$SS{\cal F}\subset FT^*X$
is the union of the $0$-section
and the line spanned by $dy$
at the origin,
the $k'$-singular support
$SS_{k'}{\cal F}\subset FT^*X$
and the $k'$-relative singular support
$\overline{SS}_{k'}{\cal F}\subset T^*X/k'$
are the $0$-sections.
On the other hand,
the singular support
defined by Beilinson in \cite{SS}
is the union of the $0$-section
and the line spanned by $dx$
at the origin.
\medskip

We don't know the existence of $SS_S{\cal F}$ in general
since the condition that
${\cal F}$ is $S$-micro supported 
on $C$ and $C'$ does not a priori imply
that 
${\cal F}$ is $S$-micro supported 
on $C\cap C'$.

\begin{lm}\label{SSSS}
Assume that $SS{\cal F}$,
$SS_S{\cal F}$ and
$SS_S^{\rm sat}{\cal F}$ exist.

{\rm 1.}
We have
$$SS_S{\cal F}\subset
SS_S^{\rm sat}{\cal F}$$
and the equality holds
if the morphism
$FT^*S\times_SX\to FT^*X$ is $0$.

{\rm 2.}
We have $$SS{\cal F}\subset
SS_S^{\rm sat}{\cal F}.$$

{\rm 3.}
If the condition {\rm (T)} in Lemma {\rm \ref{lmS}.2} is satisfied,
we have
$$SS_S{\cal F}\subset
SS{\cal F}.$$
In particular, this inclusion is satisfied
if $S$ is a smooth scheme
over a perfect field $k$ of characteristic $p>0$.
\end{lm}

\proof{
1.
Clear.

2.
This follows from Lemma \ref{lmS}.1.

3.
The inclusion follows from Lemma {\rm \ref{lmS}.2}.
In the last case,
condition (T) is satisfied by Lemma \ref{lmcn}.
\qed

}
\medskip

We reduce the proof of the existence of
$SS_S^{\rm sat}{\cal F}$
on a regular scheme $X$ 
of finite type over $S$
to the case
where $X$ is a projective space
${\mathbf P}_S$ over $S$
as in Corollary \ref{corSS}.


\begin{lm}\label{lmSSlocsat}
Let $X$ be a regular scheme of
finite type over $S$ and
let ${\cal F}$ be a constructible
sheaf on $X$.

{\rm 1.}
Let $U\subset X$ be an open subscheme.
Assume that $C\subset FT^*X$ is the 
$S$-saturated singular support
of ${\cal F}$.
Then, $C|_U$ is the $S$-saturated singular support
of ${\cal F}|_U$.

{\rm 2.}
Let $(U_i)$ be an open covering of
$X$ and $C_i
\subset FT^*U_i$ be the $S$-saturated singular support
of ${\cal F}|_{U_i}$.
Then,
$C=\bigcup_iC_i
\subset FT^*X$ is the $S$-saturated singular support
of ${\cal F}$.
\end{lm}

\proof{
The same argument as in the proof
of Lemma \ref{lmSSloc} works using
Lemma \ref{lmCms}.
\qed

}

\begin{pr}\label{prSSXPsat}
Let $i\colon X\to P$ be a closed immersion
of regular schemes of
finite type over $S$ and
let ${\cal F}$ be a constructible
sheaf on $X$.
Let $C_P\subset FT^*P$ 
be a closed conical subset
and
assume that $C_P$ is the $S$-saturated singular support
of $i_*{\cal F}$.
Then the following holds.

{\rm 1.}
$C_P$ is a subset of
$FT^*P|_X$.

{\rm 2.}
Define $C\subset FT^*X$ 
to be the closure of
the image of $C_P$
by the surjection
$FT^*P|_X\to FT^*X$.
Then $C$ is the $S$-saturated singular support $SS_S^{\rm sat}{\cal F}$.
\end{pr}

\proof{
1.
The same argument as in the proof
of Proposition \ref{prSSXP} works using
Lemma \ref{lmCms}.

2.
By Lemma \ref{lmSSXPS},
${\cal F}$ is $S$-micro supported on $C$.
Since $C_P$ is $S$-saturated,
$C$ is also $S$-saturated.
We show that $C$ is the smallest.
Assume that ${\cal F}$ is $S$-micro supported on 
an $S$-saturated closed conical subset
$C'\subset FT^*X$.
Then, since $i_*{\cal F}$ is $S$-micro supported on $i_\circ C'$
by  Lemma \ref{lmCbarms}.2,
we have $C_P\subset i_\circ C'$.
Thus the image of $C_P$
by the surjection
$FT^*P|_X\to FT^*X$
is a subset of $C'$ and
we obtain $C\subset C'$.
\qed

}
\medskip

We will prove the existence of
$SS_S^{\rm sat}{\cal F}$
for the projective space
${\mathbf P}_S$ over $S$
using the Radon transform,
following Beilinson's argument in \cite{SS}.
Let $S$ be a scheme and
$V={\mathbf A}^{n+1}_S$
be the vector bundle of rank $n+1$.
Let ${\mathbf P}={\mathbf P}_S={\mathbf P}(V)$
be the projective space bundle of dimension $n$
parametrizing lines in $V$.
The dual projective space ${\mathbf P}^\vee=
{\mathbf P}(V^\vee)$
is the moduli space of hyperplanes
in ${\mathbf P}$.

By the exact sequence
$0\to \Omega^1_{{\mathbf P}/S}(1)
\to {\cal O}_{\mathbf P}\otimes V^\vee
\to {\cal O}_{\mathbf P}(1)\to 0$
of locally free ${\cal O}_{\mathbf P}$-modules,
we define a closed subscheme
$Q={\mathbf P}(T^*{\mathbf P}/S)
\subset 
{\mathbf P}\times_S {\mathbf  P}(V^\vee)=
{\mathbf P}\times_S {\mathbf  P}^\vee$
of codimension 1.
This equals the universal family of hyperplanes
since it is defined by the tautological section
$\Gamma({\mathbf P}\times_S {\mathbf  P}^\vee,
{\cal O}(1)\boxtimes{\cal O}(1))
=V^\vee\otimes V$ corresponding
to the identity $1\in {\rm End}(V)$.
Let $q\colon Q\to {\mathbf P}$
and $q^\vee\colon Q\to {\mathbf P}^\vee$
be the restrictions of the projections
${\mathbf P}\times_S {\mathbf  P}^\vee
\to {\mathbf P}$
and
${\mathbf P}\times_S {\mathbf  P}^\vee
\to {\mathbf P}^\vee$.
By symmetry, $Q\subset {\mathbf P}\times_S {\mathbf  P}^\vee$ is identified with 
${\mathbf P}(T^*{\mathbf P}^\vee/S)$.

The conormal bundle
$L_Q=T^*_Q({\mathbf P}\times_S {\mathbf P}^\vee)
\subset (T^*{\mathbf P}/S\times_S T^*{\mathbf P}^\vee/S)|_Q$
is a line bundle.
Since $1\in {\rm End}(V)
=V^\vee\otimes V$ regarded as
a global section of
${\cal O}(1)\boxtimes{\cal O}(1)$
is the bilinear form defining $Q
\subset {\mathbf P}\times_S {\mathbf P}^\vee$,
the morphism
$N_{Q/({\mathbf P}\times_S {\mathbf P}^\vee)}
\to \Omega^1_{({\mathbf P}\times_S {\mathbf P}^\vee)/
{\mathbf P}^\vee}
\otimes_{{\cal O}_{{\mathbf P}\times_S {\mathbf P}^\vee}}
{\cal O}_Q
=
\Omega^1_{{\mathbf P}/S}
\otimes_{{\cal O}_{\mathbf P}}
{\cal O}_Q$
defines a tautological sub invertible sheaf
on $Q={\mathbf P}(T^*{\mathbf P}/S)$.
In other words,
the tautological sub line bundle
$L
\subset T^*{\mathbf P}/S\times_{\mathbf P}Q$ 
is the image of $L_Q$ by
the first projection
${\rm pr}_1\colon
(T^*{\mathbf P}/S\times_S T^*{\mathbf P}^\vee/S)|_Q
\to 
T^*{\mathbf P}/S\times_{\mathbf P}Q$.
By symmetry,
the image by
the second projection
equals
the tautological sub line bundle $L^\vee$
on $Q={\mathbf P}(T^*{\mathbf P}^\vee/S)$.

Since the conormal bundle $L_Q$
is the kernel of the surjection
$(T^*{\mathbf P}/S\times_S
T^*{\mathbf P}^\vee/S)|_Q\to T^*Q$,
the intersection
$q^\circ T^*{\mathbf P}/S
\cap q^{\vee \circ}T^*{\mathbf P}^\vee/S
\subset T^*Q$
equals the image of
the tautological bundle
$L\subset
T^*{\mathbf P}/S\times_{\mathbf P}Q$.
By symmetry, the intersection also 
equals the image of
the tautological bundle
$L^\vee\subset
T^*{\mathbf P}^\vee/S\times_{{\mathbf P}^\vee}Q$.

Let $C\subset T^*{\mathbf P}/S$ denote
a closed conical subset.
We define the Legendre transform
$C^\vee\subset T^*{\mathbf P}^\vee/S$
to  be $q^\vee_\circ q^\circ C$.
We consider projectivizations
${\mathbf P}(C)
\subset {\mathbf P}(T^*{\mathbf P}/S)$
and
${\mathbf P}(C^\vee)
\subset {\mathbf P}(T^*{\mathbf P}^\vee/S)$
as closed subsets of
$Q$.
Let $C^+\subset T^*{\mathbf P}/S$ denote
the union of $C$ and the $0$-section.

In the following,
for a smooth scheme $X$ over
$S$ and a closed conical subset $C\subset T^*X/S$,
we say that a constructible sheaf ${\cal F}$
on $X$ is micro supported on $C$,
if it is so on the corresponding
$S$-saturated closed conical subset of $FT^*X$,
and similarly for other properties
by abuse of terminology.

\begin{pr}\label{prL}
Let $C\subset T^*{\mathbf P}/S$
be a closed conical subset.
Let 
$E={\mathbf P}(C)
\subset 
Q={\mathbf P}(T^*{\mathbf P}/S)$
be the projectivization.
Let $L_Q=T^*_Q({\mathbf P}\times_S {\mathbf P}^\vee)
\subset
(T^*{\mathbf P}\times_S T^*{\mathbf P}^\vee)|_Q$
be the conormal line bundle.

{\rm 1.}
The projectivization
$E={\mathbf P}(C)
\subset Q$
is the complement of the largest open subset
where $(q,q^\vee)$ is $C$-acyclic.


{\rm 2.}
The Legendre transform
$C^\vee$
equals the union of the image of
$L|_E
\subset q^\circ T^*{\mathbf P}
\cap q^{\vee\circ}T^*{\mathbf P}^\vee
\subset T^*Q$
and its base.
%
We have ${\mathbf P}(C)
={\mathbf P}(C^\vee)$.

{\rm 3.}
We have $C^{\vee\vee}
\subset C^+$.
\end{pr}

\proof{
1.
The kernel 
${\rm Ker}((T^*{\mathbf P}/S\times_S
T^*{\mathbf P}^\vee/S)|_Q\to T^*Q)$
equals the conormal bundle
$L_Q$ and the first projection
induces an isomorphism
$L_Q\to L\subset T^*{\mathbf P}/S\times_{\mathbf P}Q$
to the tautological bundle.
By this isomorphism,
the intersection
$(C\times_S T^*{\mathbf P}^\vee/S)|_Q
\cap L_Q$ is identified with
$q^*C\cap L$.
Hence $(q,q^\vee)$
is $C$-acyclic on $U\subset Q$
if and only if the restriction
$(q^*C\cap L)|_U$ is a subset of 
the $0$-section.
Since the projectivization
$E={\mathbf P}(C)\subset
{\mathbf P}(T^*{\mathbf P}/S)$
equals
${\mathbf P}(q^* C\cap L)\subset
{\mathbf P}(L)=Q$,
the assertion follows.

2.
The Legendre transform $C^\vee$
is the image of 
the intersection
$q^\circ C\cap
q^{\vee \circ}T^*{\mathbf P}^\vee/S
\subset T^*Q$
by $q^{\vee \circ}T^*{\mathbf P}^\vee/S
\to T^*{\mathbf P}^\vee/S$.
We identify the intersection
$q^\circ T^*{\mathbf P}/S
\cap q^{\vee\circ}T^*{\mathbf P}^\vee/S
\subset T^*Q$
with the tautological line bundle
$L\subset T^*{\mathbf P}/S\times_{\mathbf P}Q$.
Then, the intersection
$q^\circ C\cap 
q^{\vee \circ}T^*{\mathbf P}^\vee/S
\subset T^*Q$
is identified with
$q^* C\cap L$.
Since the projectivization
$E={\mathbf P}(C)\subset
{\mathbf P}(T^*{\mathbf P}/S)$
equals
${\mathbf P}(q^* C\cap L)\subset
{\mathbf P}(L)=Q$,
the closed conical subset $q^* C\cap L$
equals $L|_E$ up to the base.

%
%
Since 
$L|_E$ is identified with
$L^\vee|_E$ inside 
$T^*Q$,
we have
${\mathbf P}(C^\vee)=
{\mathbf P}(L^\vee|_E)
=E\subset
{\mathbf P}(T^*{\mathbf P}^\vee/S)=Q$.

3.
By 2 and symmetry,
we have
${\mathbf P}(C)
={\mathbf P}(C^{\vee})
={\mathbf P}(C^{\vee\vee})$.
Hence we have
$C^{\vee\vee}\subset C^+$.
\qed

}

\begin{lm}\label{lmLeg}
Let $C \subset T^*{\mathbf P}/S$
be a closed conical subset and
$E={\mathbf P}(C)
\subset Q=
{\mathbf P}(T^*{\mathbf P}/S)$
be its projectivization.
Let 
$$\xymatrix{
{\mathbf P} &
Q\ar[l]_q\ar[d]_{q^\vee}&
Q_W\ar[l]_{h'}\ar[d]_{q^\vee_W}
\ar[rd]^{f'}&
\\
&{\mathbf P}^\vee
&W\ar[l]_{h}\ar[r]^f&
Y}$$ be a commutative
diagram of smooth schemes over $S$
with cartesian square.
Let $C^{\vee+}
=C^\vee\cup  T^*_{{\mathbf P}^\vee}
{\mathbf P}^\vee/S
\subset T^*{\mathbf P}^\vee/S$
be the union of the Legendre transform
with the $0$-section and
suppose that $(h,f)$ is $C^{\vee+}$-acyclic.

{\rm 1.}
The morphism $f$ is smooth
on an open neighborhood of
$W_{{\mathbf F}_p}$
and the pair
$(qh',f')$ is
$C$-acyclic on a neighborhood of
$(Q_W\sm E_W)_{{\mathbf F}_p}$.

{\rm 2.}
On a neighborhood of $E_{W,{\mathbf F}_p}$,
the pair $(qh',f')$ is
$T^*{\mathbf P}/S$-acyclic.
\end{lm}

\proof{
1.
Since $(h,f)$ is
$C^{\vee+}$-acyclic over $S$ and 
$C^{\vee+}$ contains the $0$-section,
the morphism $f\colon W\to Y$ is smooth.
By Proposition \ref{prL}.1,
$(q,q^\vee)$ is $C$-acyclic over $S$
outside $E$.
Hence
$(qh',q^\vee_W)$ is $C$-acyclic over $S$
outside $E_W$ by Lemma \ref{lmuniCacy}.1
and 
$(qh',f')$ is $C$-acyclic over $S$
outside $E_W$ by Lemma \ref{lmCacych}.

2.
By the description of $C^\vee$
in Proposition \ref{prL}.3
and by the open condition Lemma \ref{lmCCbar},
the $C^{\vee}$-acyclicity over $S$ of
$(h,f)$ implies
the $T^*{\mathbf P}$-acyclicity over $S$ of
$(qh',f')$
on a neighborhood $U$ of $E_W\subset Q_W$.
\qed

}

\medskip

We define the naive Radon transform
$R{\cal F}$ to be
$q^\vee_*q^*{\cal F}$ and
the naive inverse Radon transform
$R^\vee{\cal G}$ to be
$q_*q^{\vee*}{\cal G}$.
Since we consider only the naive Radon transform,
we drop the adjective `naive' in the sequel.
We will refine in Proposition \ref{prSSR},
after studying the difference between
$R^\vee R{\cal F}$
and ${\cal F}$,
the following elementary property.

\begin{lm}\label{lmL}
Assume that ${\cal F}$ is $S$-micro supported on $C$.

{\rm 1.}
The Radon transform
$R{\cal F}$ is $S$-micro supported on $C^\vee$.

{\rm 2.}
$R^\vee R{\cal F}$ is $S$-micro supported on $C^+$.
\end{lm}

\proof{
1.
By Lemma \ref{lmCbarms},
the Radon transform
$R{\cal F}=q^\vee_*q^*{\cal F}$
is $S$-micro supported
on the Legendre transform
$C^\vee=q^\vee_\circ q^\circ C$.

2.
By 1,
$R^\vee R{\cal F}$ is $S$-micro supported on 
$C^{\vee\vee}\subset C^+$.
\qed
}

\begin{lm}\label{lmRn}
We consider the commutative diagram
$$
\xymatrix{
Q\times_{{\mathbf P}^\vee}Q
\ar[r]^-i\ar[rd]_{q\times q}&
{\mathbf P}\times_S
{\mathbf P}^\vee\times_S {\mathbf P}
\ar[d]^{{\rm pr}_{13}}
\\
&
{\mathbf P}\times_S 
{\mathbf P}
&
{\mathbf P}\ar[l]_-{\delta_{\mathbf P}}
}
$$
where 
$\delta_{\mathbf P}
\colon {\mathbf P}\to {\mathbf P}\times_S
{\mathbf P}$
is the diagonal immersion.

{\rm 1.}
The closed immersion $i=((q,q^\vee),(q^\vee,q))$
induces isomorphisms
\begin{equation}
R^s(q\times q)_*\Lambda_{
Q\times_{{\mathbf P}^\vee}Q}
\to 
\begin{cases}
\Lambda(-t)[-2t]
&\text{if }s=2t
\text{ and }0\leqslant t\leqslant n-2,\\
\delta_{{\mathbf P}*}\Lambda(-(n-1))[-2(n-1)]
&\text{if }s=2(n-1),
\\0&\text{if otherwise.}
\end{cases}
\label{eqPQP}
\end{equation}

{\rm 2.}
Let 
$p\colon {\mathbf P}\to S$
and $p^\vee\colon {\mathbf P}^\vee\to S$
denote the projections.
Then, we have a distinguished triangle
\begin{equation}
\to 
\tau_{\leqq 2(n-2)}
(p\times p)^*
p^\vee_*\Lambda
\to
(q\times q)_*\Lambda_{
Q\times_{{\mathbf P}^\vee}Q}
\to
\delta_{{\mathbf P}*}
\Lambda_{\mathbf P}(n-1)[2(n-1)]\to.
\label{eqppL}
\end{equation}
\end{lm}

\proof{
1.
The immersion $i$ induces a morphism
\begin{equation}
(p\times p)^*
p^\vee_*\Lambda
={\rm pr}_{13*}\Lambda
\to 
(q\times q)_*\Lambda_{
Q\times_{{\mathbf P}^\vee}Q
}
\label{eqQQ}
\end{equation}
and we have isomorphisms
$R^sp^\vee_*\Lambda
\to 
\Lambda(-t)[-2t]$
for $s=2t, 0\leqslant t\leqslant n$
and 
$R^sp^\vee_*\Lambda
=0$ otherwise.
The restriction of
the closed immersion $
i\colon Q\times_{{\mathbf P}^\vee}Q
\to
{\mathbf P}\times_S 
{\mathbf P}^\vee
\times_S {\mathbf P}$
on the diagonal
${\mathbf P}\subset
{\mathbf P}\times_S {\mathbf P}$
is the sub ${\mathbf P}^{n-1}$-bundle $Q
\subset {\mathbf P}\times_S{\mathbf P}^\vee$.
On the complement
${\mathbf P}\times_S {\mathbf P}
\sm {\mathbf P}$,
it is a sub
${\mathbf P}^{n-2}$-bundle.
Hence (\ref{eqQQ}) induces an isomorphism
of cohomology sheaves except
for degree $s$ and
induces an isomorphism
$\delta_{{\mathbf P}*}
p^*R^{2(n-1)}p^\vee_*\Lambda
\to 
R^{2(n-1)}(q\times q)_*\Lambda_{
Q\times_{{\mathbf P}^\vee}Q}$.

2.
This follows from the isomorphisms (\ref{eqPQP}).
\qed

}
\medskip

Next, we consider the diagram
\begin{equation}
\begin{CD}
{\mathbf P}@<{{\rm pr}_1}<<
{\mathbf P}\times_S {\mathbf P}
@<{q\times q}<<
Q\times_{{\mathbf P}^\vee}Q
\\
@.@V{{\rm pr}_2}VV@.\\
@.{\mathbf P}.@.
\end{CD}
\label{eqRRv}
\end{equation}

\begin{pr}\label{prRn}
{\rm 1.}
We have a canonical isomorphism
\begin{equation}
R^\vee R{\cal F}
\to
R{\rm pr}_{2*}\bigl({\rm pr}_1^*{\cal F}
\otimes R(q\times q)_*\Lambda_{
Q\times_{{\mathbf P}^\vee}Q}\bigr).
\label{eqRRF}
\end{equation}

{\rm 2.}
The isomorphism {\rm (\ref{eqRRF})}
induces a
distinguished triangle
$$
\to
\bigoplus_{q=0}^{n-2}
p^*p_*{\cal F}(-q)[-2q]
\to
R^\vee R{\cal F}
\to 
{\cal F}(-(n-1))[-2(n-1)]
\to.
$$
\end{pr}

\proof{
1.
By the cartesian diagram
$$\begin{CD}
{\mathbf P}@<q<<
Q@<{{\rm pr}_1}<<
Q\times_{{\mathbf P}^\vee}Q
\\
@.@V{q^\vee}VV@VV{{\rm pr}_2}V\\
@.{{\mathbf P}^\vee}@<{q^\vee}<<Q\\
@.@.@VVqV\\
@.@.{\mathbf P},
\end{CD}$$
we have 
$R^\vee R{\cal F}=
Rq_*q^{\vee*}Rq^\vee_*q^*{\cal F}$.
By the proper base change theorem,
we have a canonical isomorphism
$Rq_*q^{\vee*}Rq^\vee_*q^*{\cal F}
\to 
R(q\circ {\rm pr}_2)_*
(q\circ {\rm pr}_1)^*{\cal F}$.
In the notation of (\ref{eqRRv}),
the latter is identified with
$R({\rm pr}_2\circ (q\times q))_*
({\rm pr}_1\circ (q\times q))^*{\cal F}$.
This is identified with
$R{\rm pr}_{2*}({\rm pr}_{1*}{\cal F}
\otimes 
R(q\times q)_*\Lambda_{Q\times_{{\mathbf P}^\vee}Q})$ by the projection formula.

2.
This follows from the isomorphism (\ref{eqRRF})
and the distinguished triangle (\ref{eqppL}).
}

\begin{pr}\label{prSSR}
For a constructible sheaf ${\cal F}$ on ${\mathbf P}$
and a closed conical subset
$C\subset T^*{\mathbf P}$,
we have implications
{\rm (1)}$\Rightarrow${\rm (2)}$\Rightarrow${\rm (3)}$\Rightarrow${\rm (4)}.

{\rm (1)}
${\cal F}$ is $S$-micro supported on $C$.

{\rm (2)}
The pair
$(q,q^{\vee})$
is universally ${\cal F}$-acyclic
over $S$ outside $E={\mathbf P}(C)$.

{\rm (3)}
The Radon transform
$R{\cal F}$ is $S$-micro supported on $C^{\vee+}$.

{\rm (4)}
${\cal F}$ is $S$-micro supported on $C^+$.
\end{pr}

%

\proof{
{\rm (1)}$\Rightarrow${\rm (2)}:
The pair $(q,q^\vee)$
of $q\colon Q\to {\mathbf P}$
and $q^\vee\colon Q\to {\mathbf P}^\vee$ is 
$C$-acyclic over $S$
outside $E={\mathbf P}(C)$
by Proposition \ref{prL}.1.
Hence (1) implies that
$(q,q^\vee)$ is universally
${\cal F}$-acyclic over $S$ on the complement  
$Q\sm E$ by Lemma \ref{lmSbase}.2.

{\rm (2)}$\Rightarrow${\rm (3)}:
Assume that a pair of morphisms $h\colon W\to {\mathbf P}^\vee,
f\colon W\to Y$ is $C^{\vee+}$-acyclic over $S$
and show that
$(h,f)$ is $R{\cal F}$-acyclic over $S$.
We consider the commutative diagram
$$\xymatrix{
{\mathbf P} &
Q\ar[l]_q\ar[d]_{q^\vee}&
Q_W\ar[l]_{h'}\ar[d]_{q^\vee_W}
\ar[rd]^{f'}&
\\
&{\mathbf P}^\vee
&W\ar[l]_{h}\ar[r]^f&
Y}$$
with cartesian square as in Lemma \ref{lmLeg}.
Since $q^\vee$ is proper,
it suffices to show that
$(h',f')$ is $q^*{\cal F}$-acyclic over $S$
by Lemma \ref{lmFacycg}.
Since $q$ is smooth,
it suffices to show that
$(qh',f')$ is ${\cal F}$-acyclic over $S$
by Lemma \ref{lmFacyc} (1)$\Rightarrow$(2).

By (2), 
$(qh',q^\vee_W)$ is ${\cal F}$-acyclic over $S$
on a neighborhood of $(Q_W\sm E_W)_{{\mathbf F}_p}$.
By Lemma \ref{lmLeg}.1,
the morphism $f\colon W\to Y$ is smooth
on an open neighborhood of
$W_{{\mathbf F}_p}$.
Hence
$(qh',f')$ is ${\cal F}$-acyclic over $S$
on a neighborhood of $(Q_W\sm E_W)_{{\mathbf F}_p}$
by Lemma \ref{lmFacycfun}.

By Lemma \ref{lmLeg}.2,
the pair $(qh',f')$ is
$T^*{\mathbf P}$-acyclic over $S$
on a neighborhood $U$ of $E_W\subset Q_W$.
Since ${\cal F}$ is $S$-micro supported on 
$FT^*{\mathbf P}$
by Lemma \ref{lmFperp}.1,
$(qh',f')$ is ${\cal F}$-acyclic over $S$
on a neighborhood of $E_{W,{\mathbf F}_p}$.
Thus 
$(qh',f')$ is ${\cal F}$-acyclic over $S$
on a neighborhood of  $Q_{W,{\mathbf F}_p}$
as required.

(3)$\Rightarrow$(4)
By (3) and (1)$\Rightarrow$(3),
$R^\vee R{\cal F}$ is $S$-micro supported on $(C^{\vee+})
^{\vee+}=
C^+$.
Since $p^*p_*{\cal F}$
is $S$-micro supported
on the $FT^*S\times_S{\mathbf P}$,
by the distinguished triangle in
Proposition \ref{prRn}.2,
${\cal F}$  is also $S$-micro supported on $C^+$.
\qed

}

\begin{cor}\label{corSSR}
Let ${\cal F}$ be a constructible
sheaf on ${\mathbf P}$.
Let $E\subset Q={\mathbf P}(T^*{\mathbf P})$ be 
the complement of
the largest
open subset on which
$(q,q^\vee)$
is universally ${\cal F}$-acyclic over $S$.
Then the closed conical
subset $C\subset FT^*{\mathbf P}$
corresponding to the base
$B={\rm supp}\, {\cal F}
\cap {\mathbf P}_{{\mathbf F}_p}$ 
and $E\cap Q_{{\mathbf F}_p}
\subset {\mathbf P}(T^*{\mathbf P}/S)_{{\mathbf F}_p}$
is the $S$-saturated
singular support $SS_S^{\rm sat}{\cal F}$
of ${\cal F}$.
\end{cor}

\proof{
By Proposition \ref{prSSR} (2)$\Rightarrow$(4),
${\cal F}$ is $S$-micro supported on $C^+$.
Hence
${\cal F}$ is $S$-micro supported on $
C=C^+|_B$ by Lemma \ref{lmCms}.1.

Assume that ${\cal F}$ is $S$-micro supported on $C'
\subset T^*{\mathbf P}/S_{{\mathbf F}_p}$.
Then, by Proposition \ref{prSSR} (1)$\Rightarrow$(2),
we have ${\mathbf P}(C')\supset E={\mathbf P}(C)$
since $E$ is the smallest.
Since the base of $C'$ contains
$B={\rm supp}\, {\cal F}$
as a subset,
we have $C'\supset C$
by Lemma \ref{lmCC'}.
\qed

}

\begin{thm}\label{thmSS}
Let $S$ be an excellent regular scheme
over ${\mathbf Z}_{(p)}$
such that the reduced part of
$S_{{\mathbf F}_p}$ is of
finite type over a field $k$ such that $[k:k^p]$ is finite.
Let $X$ be a regular scheme of
finite type over $S$ and ${\cal F}$ be
a constructible sheaf on $X$.
Then, the $S$-saturated singular support
$SS_S^{\rm sat}{\cal F}$ exists.
\end{thm}

\proof{
By the same argument as the proof of
Corollary \ref{corSS}
using Lemma \ref{lmSSlocsat} 
and Proposition \ref{prSSXPsat}
in place of Lemma \ref{lmSSloc} 
and Proposition \ref{prSSXP},
it suffices to prove the case where
$X={\mathbf P}_S$.
This case is proved in Corollary \ref{corSSR}.
\qed

}
\medskip

By the same argument as in Lemma \ref{SSSS},
we have the following.

\begin{cor}
Assume further that the morphism
$FT^*S\times_SX\to FT^*X$
is $0$.
Then the $S$-singular support
$SS_S{\cal F}$ exists and equals
$SS_S^{\rm sat}{\cal F}$.
If we further assume that the condition
{\rm (T)} in Lemma {\rm \ref{lmS}.2}
is satisfied, then
$SS{\cal F}$ exists and equals
$SS_S{\cal F}=SS_S^{\rm sat}{\cal F}$.
\end{cor}

\end{document}